# Analytic Extension of Keiper-Li Coefficients

Krzysztof Dominik Maślanka

Polish Academy of Sciences

Institute for the History of Science

Nowy Świat 72, 00-330 Warsaw, Poland

e-mail krzysiek2357@gmail.com

December 18, 2022

**Abstract**

We construct certain entire function $\lambda(s)$ which for integer $s$ coincides with the well-known Keiper-Li coefficients (see [7], [10]), i.e. $\lambda(n) = \lambda_n$. This is an even function $\lambda(s) = \lambda(-s)$ and has an infinitude of complex zeros exhibiting interesting distribution. Extensive tables of more than 3500 complex zeros of $\lambda(s)$ with precision of 14 significant digits are included. A detailed analysis of the distribution of these zeros may shed some light on the Riemann hypothesis. It turns out that possible violation of the Riemann hypothesis (if such is the case) would be clearly reflected in the specific distribution of these zeros. More specifically, they form complex quadruplets if the Riemann hypothesis is true, and aligned real doublets if it is false.

MSC classes: 11-04, 11M06, 68W30

## 1 Introduction and definitions

One of the objectives of experimental mathematics is to look for hidden, unexpected patterns or structures using modern computational techniques. Of course, these experiments are not as valuable as the rigorous and immutable proofs that have been the goal and attribute of mathematics for centuries. But they can sometimes lead to the right track and suggest specific, useful directions for theoretical research – in the spirit of an ambitious program formulated as a kind of manifesto by Jonathan Borwein, an enthusiast and ardent propagator of experimental mathematics. In particular, Borwein recommends using computers to gain insight and intuition as well as discover new patterns and relationships ([4], page 2). And this is what we intend to do in this paper.

Here we will describe such an approach based on a specific example taken from analytic number theory. But first, let us recall some definitions and the notation adopted.

Function $\xi(s)$ introduced by Riemann in his fundamental paper [17] (see also [5], p.16-21, for a detailed exposition):

$$\xi(s) = 2(s-1)\pi^{-s/2}\Gamma\left(1 + \frac{s}{2}\right)\zeta(s) \qquad (1)$$

is related to the zeta function $\zeta(s)$, having the same complex (nontrivial) zeros $\rho_n$. All the so-called trivial zeros of the zeta function located at the negative even integers have been "annihilated" by the simple poles of the gamma function. $\xi(s)$ is an entire function[1] and satisfies simple functional equation:

$$\xi(s) = \xi(1-s) \qquad (2)$$

---

[1] Strictly speaking, Riemann considered another function which is related to (1), namely $\xi(\frac{1}{2} + it)$. Since it was sometimes confusing, Edmund Landau introduced the present convention.

There are different conventions concerning the trivial factor in the definition of $\xi(s)$. Here we have chosen such that seems the most natural: $\xi(0) = \xi(1) = 1$.

Let us further define:
$$f(s) = \ln \xi(s) \tag{3}$$
This function is of crucial importance since it appears in the definition of certain real coefficients:
$$\lambda_n = \frac{1}{\Gamma(n)} \lim_{s \to 1} \frac{d^n}{ds^n}\left[s^{n-1} f(s)\right] \tag{4}$$

As was outlined by Keiper in 1992 [7] and rigorously proved by Li in 1997 [10] the celebrated Riemann hypothesis (RH) [17] is equivalent to nonnegativity of coefficients (4):
$$\text{RH} \iff \lambda_n \geq 0 \tag{5}$$

Of course, criterion (5), while deceptively simple, at least in its concise notation, is as difficult to settle as the Riemann hypothesis itself. Nevertheless it is worth detailed investigations, see [3], [6], [8], [12], [16], [20]. It should be emphasized that the coefficients $\lambda_n$ in the form (4), although important from theoretical point of view, are practically useless for numerical estimates.

Introducing function $\varphi(s)$:
$$\varphi(s) := \xi\left(\frac{1}{1-s}\right) \tag{6}$$
one can easily show that $\lambda_n$ are power series coefficients of the logarithmic derivative of $\varphi(s)$ (see [10] for details):
$$\frac{\varphi\prime(s)}{\varphi(s)} = \sum_{n=0}^{\infty} \lambda_{n+1}\, s^n \qquad |s| < 1 \tag{7}$$

The third formula for $\lambda_n$ coefficients involves explicitly complex zeros of the Riemann zeta function located within the critical strip:
$$\lambda_n = \sum_{\rho}\left[1 - \left(1 - \frac{1}{\rho}\right)^n\right] \tag{8}$$
where the summation is taken over all these zeros with $\rho$ and $1 - \rho$ being paired together.

A detailed inspection of the formula (8) leads to a rather pessimistic conclusion: any attempt to locate a counterexample for RH by finding a negative $\lambda_n$ is hopeless. Note that if $\rho = 1/2 + i\, y$ (for real $y$) then $|1 - \frac{1}{\rho}| = 1$. However, if the $\rho$ deviates from the critical line then $|1 - \frac{1}{\rho}|$ is slightly different from one. Then $n$ should be extremely great, especially for high $\rho$, to make $\lambda_n$ negative.

Note that (8) can be thought of as an analytical extension of the coefficients $\lambda_n$ to an entire function $\lambda(n)$, since $n$ need not be an integer at all – it can be negative, rational, or even complex. The summation in (8) would still converge over the entire complex plane. Unfortunately, this convergence is very slow – so slow that, for example, numerical calculation of the position of zeros of the entire function $\lambda(n)$ is practically impossible.

Several important results on such an analytical extension of the Keiper-Li coefficients are contained in [8], chapter 7. In this paper Lagarias rigorously proved its existence and also showed some of its properties.

There is yet another representation of these important constants, which we will give here just for the sake of completeness. Let
$$\varphi(s) = \sum_{j=0}^{\infty} a_j\, s^j \tag{9}$$
where coefficients $a_j$ are all positive (see [10] for the proof). Then
$$\lambda_n = n \sum_{i=0}^{\infty} \frac{(-1)^{i-1}}{i} \sum a_{k_1}...a_{k_i} \tag{10}$$

where the last summation is performed over all integer indices $k_1, ..., k_i$ satisfying two constraints: $1 \leq k_1, ..., k_i \leq n$ and $k_1 + ... + k_i = n$.

Needless to say that formula (10), just like (4), while important for theoretical considerations, is rather useless when trying to estimate numerically the values of $\lambda_n$.

## 2 Historical remarks

Coefficients $\lambda_n$ were introduced by Jerry B. Keiper in 1992 [7] who also showed (without rigorously proving) that the Riemann hypothesis *implies* that $\lambda_n > 0$.

In 1997 Xian-Jin Li provided rigorous proof of a stronger statement, namely that the Riemann hypothesis is in fact *equivalent* to $\lambda_n \geqslant 0$ [10]. Li does not refer to Keiper's work. As a historian of mathematics, I wondered if Li knew or was inspired by Keiper's 1992 paper. This issue was settled thanks to first-hand information, i.e. from Li himself [11]. As is sometimes the case in the history of science, the matter is more complicated – and more interesting.

Li was unaware of Keiper's paper. The first version of his work dates back to 1988. At that time, the author was a second year graduate student at Purdue University. He sent his important result for publication in a reputable periodical, but the editors rejected it. Only after several years, encouraged by Jeffrey Lagarias, Li published his result in *Journal of Number Theory* [10].

As for terminology, the following quote from Lagarias' paper ([8], p. 1691) seems to be particularly meaningful:

> We adopt the term *Li coefficients* because Li's work shows the nonnegativity of $\lambda_n$ is equivalent to the Riemann hypothesis and because he generalized them to various other zeta functions.

In my opinion, it would be correct, and this is generally the case in the literature on the subject, to use the name "Keiper-Li coefficients" – to pay proper credit to the late Jerry Keiper (see e.g. [6], [16], [20]), while the criterion (5) for the Riemann hypothesis should bear only Li's name: "the Li criterion for the Riemann hypothesis".

Jerry Keiper was killed in a road accident at the age of 41 on January 18, 1995, after being hit by a car while cycling home after work at Wolfram Research [22]. He was a brilliant computer programmer but also had an intuitive sense of number theory. His contribution to the development of algorithms later implemented in Wolfram Mathematica cannot be overestimated. In particular, he developed effective methods of calculating the values of many special functions, for any complex argument, with any accuracy – a really formidable achievement.

In his important 1992 paper [7], without paying special attention to mathematical precision, he simply wrote: "It is clear that the Riemann hypothesis implies that $\lambda_n > 0$ for all positive $n$", see [7], p. 769. He did not even care to promote this modest sentence, with the manner typical of mathematicians, to the formal theorem: "A new criterion for the Riemann hypothesis". There are a few more such inconspicuous but profound statements in his paper.

Keiper was also an extraordinary personality. As an extremely honest pacifist, he gave up his salary at Wolfram Research. He justified this in a way that was probably incomprehensible and even offensive to the majority of modern greedy people: the salaries are subject to (compulsory) tax, some of which, without the taxpayer's knowledge and consent, is spent on armaments. His former teacher put it like this:

> Jerry Keiper was indeed a most remarkable individual. [...] As a Mennonite [member of the radical Protestant church] he was very pious, but did lead by example rather than by "preaching". As a law-abiding citizen he had to pay taxes but as a pacifist he could not help support the military. His solution was to lead a life of poverty. As a very talented individual whose expertise was tremendously in demand (computer science) this was actually very difficult to do. [...] In a way he was just too good for this world [18].

# 3 Interpolation with nonequidistant nodes

Now our goal is to construct an *interpolating polynomial* convergent to $f(s) \equiv \ln \xi(s)$ in the closed interval $[0, 1]$. To achieve this we choose as interpolation nodes the following set of points:

$$u = \{1, \frac{1}{2}, \frac{1}{3}, \frac{1}{4}, \frac{1}{5}, ...\} \tag{11}$$

together with:

$$v = \{0, \frac{1}{2}, \frac{2}{3}, \frac{3}{4}, \frac{4}{5}, ...\} \tag{12}$$

These nodes are obviously unequally spaced and get denser towards the edges of the interval $[0, 1]$ on which we want to interpolate $f(s)$. Note that with such choice, point $\frac{1}{2}$ is doubled. Also note that each number from the set (11) is equal to one minus the corresponding number from the set (12), and vice versa which, of course, is intentionally related to the functional equation (2).

In the case of interpolation attempts, typical phenomenon is the so-called Runge effect [19], which for *equidistant* nodes occurs even for some perfectly smooth (on the real axis) functions. These are undesirable oscillations between the nodes, growing with increasing order of the interpolating function. The remedy for this is the appropriate choice of *unequally* spaced nodes, e.g. so called Chebyshev nodes, which get denser towards the edges of the interval on which we want to interpolate a given function. It is this uneven distribution of nodal points that guarantees the convergence of the interpolation along the entire interval, and even in its neighborhood, but, of course, not the global convergence.

But first recall the definition and properties of the Pochhammer symbol $(x)_k$ which we shall need soon:

$$(x)_k \equiv \frac{\Gamma(k+x)}{\Gamma(x)} = \prod_{i=0}^{k-1}(x+i) = (-1)^k \sum_{i=0}^{k}(-1)^i S_k^{(i)} x^i = \sum_{i=0}^{k} \left| S_k^{(i)} \right| x^i \tag{13}$$

where $x$ is complex and $k$ is a nonnegative integer. Let us now introduce the following family of polynomials of the variable $s$:

$$\eta_k(s) \equiv \prod_{i=1}^{k}\left(s - \frac{1}{i}\right) = \frac{s^k}{k!}\left(\frac{s-1}{s}\right)_k = \frac{s^k}{k!} \sum_{i=0}^{k} \left|S_k^{(i)}\right| \left(\frac{s-1}{s}\right)^i \tag{14}$$

$$\omega_k(s) = \eta_k(s)\eta_k(1-s) = \frac{1}{(k!)^2} \sum_{i,j=1}^{k} S_k^{(i)} S_k^{(j)} s^{k+i-j}(1-s)^{k+j-i} \tag{15}$$

where $S_k^{(i)}$ are (signed) Stirling numbers of the first kind (see [1] and [25]). The degree of the polynomial $\omega_k(s)$ is $2k$. As usual in such situations, we assume that if in (15) any of the exponents is zero then the corresponding power is taken to be one. It is evident that $\omega_k(s)$ is zero in $k$ initial nodes (11) and (12). Also, by construction, $\omega_k(s) = \omega_k(1-s)$, i.e. for each $k$ the functions $\omega_k(s)$ have the same symmetry (2) as $\xi$ and $f$. Also, by construction the nodes (11), (12) are taken in pairs in (15). This is realized as follows:

| $k^{\text{th}}$ order | nodes |
| --- | --- |
| 0 | $(0, 1)$ |
| 1 | $(0, 1), (\frac{1}{2}, \frac{1}{2})$ |
| 2 | $(0, 1), (\frac{1}{2}, \frac{1}{2}), (\frac{1}{3}, \frac{2}{3})$ |
| 3 | $(0, 1), (\frac{1}{2}, \frac{1}{2}), (\frac{1}{3}, \frac{2}{3}), (\frac{1}{4}, \frac{3}{4})$ |
| ... | ... |

Such an interpolation mechanism is realized by taking a linear combination of polynomials (15):

$$F_m(s) = \sum_{k=1}^{m} (-1)^k \alpha_k\, \omega_k(s) \qquad (16)$$

where $m$ tends to infinity and $\alpha_k$ are some coefficients that have yet to be determined. The condition from which we calculate them is $F_m(s) = f(s)$ if $s$ takes the value of any of the nodes (11) or (12). (Due to symmetry $f(s) = f(1-s)$ nodes (11) and (12) are equivalent.) This condition leads to a system of linear equations:

| $m$ | $F_m(u_{m+1})$ or $F_m(v_{m+1})$ |
|---|---|
| 1 | $-\frac{1}{4}\alpha_1 = f\left(\frac{1}{2}\right)$ |
| 2 | $-\frac{2}{9}\alpha_1 - \frac{1}{162}\alpha_2 = f\left(\frac{2}{3}\right)$ |
| 3 | $-\frac{3}{16}\alpha_1 - \frac{3}{256}\alpha_2 - \frac{5}{12288}\alpha_3 = f\left(\frac{3}{4}\right)$ |
| 4 | $-\frac{4}{25}\alpha_1 - \frac{9}{625}\alpha_2 - \frac{14}{15625}\alpha_3 - \frac{77}{3125000}\alpha_4 = f\left(\frac{4}{5}\right)$ |
| 5 | $-\frac{5}{36}\alpha_1 - \frac{5}{324}\alpha_2 - \frac{5}{3888}\alpha_3 - \frac{35}{559872}\alpha_4 - \frac{133}{100776960}\alpha_5 = f\left(\frac{5}{6}\right)$ |
| ... | ... |

which may readily be solved:

$$\begin{pmatrix} \alpha_1 \\ \alpha_2 \\ \alpha_3 \\ \alpha_4 \\ \alpha_5 \\ ... \end{pmatrix} = \begin{pmatrix} -4 & 0 & 0 & 0 & 0 \\ 144 & -162 & 0 & 0 & 0 \\ -2304 & \frac{23328}{5} & -\frac{12288}{5} & 0 & 0 \\ 25600 & -\frac{524880}{7} & \frac{983040}{11} & -\frac{3125000}{77} & 0 \\ -230400 & \frac{23328}{7} & -\frac{141557760}{77} & \frac{2812500000}{1463} & -\frac{100776960}{133} \\ ... & ... & ... & ... & ... \end{pmatrix} \cdot \begin{pmatrix} f\left(\frac{1}{2}\right) \\ f\left(\frac{2}{3}\right) \\ f\left(\frac{3}{4}\right) \\ f\left(\frac{4}{5}\right) \\ f\left(\frac{5}{6}\right) \\ ... \end{pmatrix} \qquad (17)$$

Rational coefficients $c_{kj}$ which form the triangular matrix in (17) were simply guessed which took some time and also required some intuition[2]:

$$c_{kj} = \binom{k}{j} k!(k+1)^2 \frac{j^{-k}(j+1)^{2k-2}}{\left(\frac{1}{j} - k\right)_{k-j} \left(-\frac{1}{j}\right)_j} \qquad (18)$$

Finally we have:

$$\alpha_k = \sum_{j=1}^{k} c_{kj}\, f\left(\frac{j}{j+1}\right) \qquad (19)$$

As usual (cf. [13]) the most efficient way to obtain numerical values of $\alpha_k$ with high precision consists in tabulating at first appropriate values of $f$ in the nodes with high precision (typically at least 10000 digits). In this way we avoid unnecessary, repetitive and time-consuming calculations of $f$ in the nodes. A small but highly effective program PARI/GP [23] is most appropriate to achieve this. Having these values one can use (19) to calculate $\alpha_k$. It turns out that they decrease very quickly to zero as $k$ grows, and their signs alternate (at least for the values available for numerical experiments, cf. Fig. 1.)

---

[2] Explaining in detail how the formula (18) was obtained would be uninformative and, worse, very boring. There was a lot of guesswork, and there was also trial and error. Not very ambitious, but still effective. First, the formula for the numbers from the main diagonal of (17) was guessed. Frequent consultation with the *On-Line Encyclopedia of Integer Sequences* [21] was very helpful at this stage.

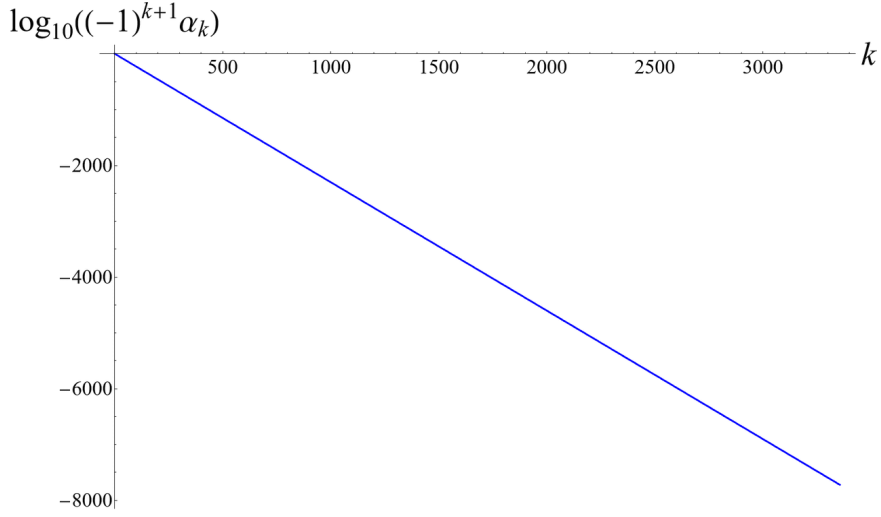

Figure 1: Real coefficients $\alpha_k$ (19) diminish quickly to zero when $k$ grows which guarantees fast convergence of series (16). Numerical calculations suggest that $(-1)^{k+1}\alpha_k$ are all positive, at least if the Riemann hypothesis is true.

## 4 Numerically effective formula for $\lambda_n$

We are now ready to find numerically efficient expression for the Keiper-Li coefficients (4). Substituting the interpolating polynomial (16) together with (15) in place of $f(s)$ in (4) and letting $m \longrightarrow \infty$ we get:

$$\lambda_n = \frac{1}{\Gamma(n)} \lim_{s \to 1} \frac{d^n}{ds^n} \left[ s^{n-1} \sum_{k=1}^{\infty} \frac{(-1)^k \alpha_k}{(k!)^2} \sum_{i,j=0}^{k} S_k^{(i)} S_k^{(j)} s^{k+i-j} (1-s)^{k+j-i} \right] = \qquad (20)$$

$$= \frac{1}{\Gamma(n)} \sum_{k=1}^{\infty} \frac{(-1)^k \alpha_k}{(k!)^2} \sum_{i,j=0}^{k} S_k^{(i)} S_k^{(j)} \lim_{s \to 1} \frac{d^n}{ds^n} \left[ s^{n-1+k+i-j} (1-s)^{k+j-i} \right]$$

Using Leibniz' formula:

$$\frac{d^n}{ds^n}(f(s)g(s)) = \sum_{r=0}^{n} \binom{n}{r} \frac{d^r}{ds^r} f(s) \frac{d^{n-r}}{ds^{n-r}} g(s) \qquad (21)$$

we get:

$$\lambda_n = \frac{1}{\Gamma(n)} \sum_{k=1}^{\infty} \frac{(-1)^k \alpha_k}{(k!)^2} \sum_{i,j=0}^{k} S_k^{(i)} S_k^{(j)} \lim_{s \to 1} \sum_{r=0}^{n} \binom{n}{r} \frac{d^r}{ds^r} \left((1-s)^{k+j-i}\right) \frac{d^{n-r}}{ds^{n-r}} \left(s^{n-1+k+i-j}\right)$$

The $n^{\text{th}}$ derivative of monomial is:

$$\frac{d^n}{ds^n}(s-c)^r = \frac{r!}{(r-n)!}(s-c)^{r-n} \qquad n \leq r \qquad (22)$$

so

$$\lambda_n = \frac{1}{\Gamma(n)} \sum_{k=1}^{\infty} \frac{\alpha_k}{(k!)^2} \sum_{i,j=0}^{k} S_k^{(i)} S_k^{(j)} \sum_{r=0}^{n} \binom{n}{r} (-1)^{j-i} \cdot$$
$$\cdot \frac{(k+j-i)!}{(k+j-i-r)!} \frac{(n-1+k+i-j)!}{(k+i-j+r-1)!} \lim_{s \to 1} (s-1)^{k+j-i-r} s^{k+i-j+r-1} \quad (23)$$

Now, after taking the limit $s \to 1$, in the summation over $r$ survives only the term for which the exponent in $(s-1)^{k+j-i-r}$ is zero, i.e. when $r = k+j-i$. Hence:

$$\lambda_n = n \sum_{k=1}^{\infty} \frac{\alpha_k}{(k!)^2} \sum_{i,j=0}^{k} S_k^{(i)} S_k^{(j)} (-1)^{j-i} \frac{1}{(n-k-j+i)!} \frac{(n-1+k+i-j)!}{(2k-1)!} = \quad (24)$$
$$= n \sum_{k=1}^{\infty} \frac{\alpha_k}{(k!)^2} \sum_{i,j=0}^{k} \left| S_k^{(i)} S_k^{(j)} \right| \binom{n+k+i-j-1}{2k-1}$$

Defining rational coefficients $\beta_{nk}$:

$$\beta_{nk} = \frac{n}{(k!)^2} \sum_{i,j=1}^{k} \left| S_k^{(i)} S_k^{(j)} \right| \binom{n+k+i-j-1}{2k-1} \quad (25)$$

we finally have:

$$\lambda_n = \sum_{k=1}^{\infty} \beta_{nk} \, \alpha_k \quad (26)$$

The rational coefficients $\beta_{nk}$ are positive for integer $n$ and increase with the increase of $k$. Nevertheless, because the coefficients $\alpha_k$ tend to zero much faster (see Fig. 1), the series representation (26) of $\lambda_n$ converges quite quickly. Therefore, it can be effectively used to calculate several thousand of initial values of $\lambda_n$, provided that we have previously tabulated numerical values of $\alpha_k$, according to the procedure outlined above. For instance, if we have precomputed values of $f(s)$ in the nodes with precision $20,000$ digits we can compute $\alpha_k$ up to $k = 3361$, since, as usual, precision of $\alpha_k$ decreases to zero with $k$. In order to increase this range, one should also (proportionally) increase the precision of the precomputed values of $f(s)$ in the nodes.

The numerical values of coefficients $\lambda_n$ calculated in this way agree perfectly with values obtained by other methods, see e.g. [7], [6], [12].

One thing should be stressed out here. In the formula (26), it is enough to sum up to some maximum $k$, more precisely, to that $k$ for which the precision of $\alpha_k$ drops to zero. Most importantly, all digits of $\lambda_n$ obtained from this summation are significant ([13], for details see comment on formula (23) and especially Fig. 9).

## 5 Analytic extension of $\lambda_n$

In what follows we are going to extend the range of index $n$ from positive integers to arbitrary complex numbers. Simply note that the values of the parameter $n$ in (26) no longer need to be positive integers as in the formula (4); on the contrary: they can be negative, rational and even complex. Thus (26) is in fact *an extension of the formula (4) to the whole plane of complex numbers*. For example:

$$\lambda(1/2) = 0.005774507219796948948...$$
$$\lambda(i) = -0.0231018963998549040 0...$$
$$\lambda(1+i) = 0.00001237486681209165... + i \, 0.04619760137033736709...$$

Our purpose is to demonstrate that $\lambda(n)$ expressed by (26) is an entire function, so it has complex zeros that can be investigated numerically. It is also, as we will show, an even function, which is not obvious when looking at (4) or (26).

As mentioned above, the variable $n$ is no longer limited to integers. This is because in the formula (25) it appears in a binomial symbol and not, for example, as an index in the Stirling numbers, where it should be integer by its very nature. Therefore we can use the following identity:

$$\binom{x}{k} = \frac{(-1)^k}{k!}(-x)_k = \frac{1}{k!}\sum_{p=0}^{k} S_k^{(p)} x^p \tag{27}$$

and put $x = n + k + i - j - 1$. After some tedious but elementary calculations we get:

$$\beta_{nk} = \frac{n}{(k!)^2} \frac{1}{\Gamma(2k)} \sum_{q=1,3,5,\ldots}^{2k-1} \chi_{kq} n^q \tag{28}$$

where

$$\chi_{kq} = \sum_{i,j=1}^{k} \sum_{p=q}^{2k-1} \left| S_k^{(i)} S_k^{(j)} \right| S_{2k-1}^{(p)} \binom{p}{q} (k+i-j-1)^{p-q} \tag{29}$$

(As usual, the power in (29) should be taken 1 when $p = q$.) It is easy to check that all the coefficients $\chi_{kq}$ are integer and form a triangular matrix. Moreover, we see from (28) that $\beta_{nk}$ are actually polynomials in variable $n$ with rational coefficients. Below there are several of them labelled by the integer parameter $k$:

| $k$ | $\beta_{nk}$ |
|---|---|
| 1 | $n^2$ |
| 2 | $\frac{1}{12}n^2\left(1 + 2n^2\right)$ |
| 3 | $\frac{1}{1080}n^2\left(11 + 40n^2 + 9n^4\right)$ |
| 4 | $\frac{1}{120960}n^2\left(151 + 756n^2 + 329n^4 + 24n^6\right)$ |
| 5 | $\frac{1}{36288000}n^2\left(5148 + 32163n^2 + 20167n^4 + 2902n^6 + 100n^8\right)$ |
| ... | ... |

In a naturally way, the integer parameter $k$ enumerates consecutive polynomials $\beta_{nk}$, whereas $n$ is an independent (in general complex) variable. Each $\beta_{nk}$ is of degree $2k$ and has, of course, $2k$ roots. The distribution of these complex roots, even if not particularly interesting from the theoretical point of view, is quite aesthetic, see Fig. 2.

Since $\chi_{kq} = 0$ for even values of $q$ we have only even exponents of $n^q$ in (28). Introducing real coefficients $\nu_q$:

$$\nu_q = \sum_{k=1}^{\infty} \frac{\chi_{k,q-1}}{\Gamma(2k)} \frac{\alpha_k}{(k!)^2} \tag{30}$$

we finally have $\lambda(n)$ in the form of a "polynomial of infinite order", i.e. the entire function $\lambda(n)$ that interpolates $\lambda_n$ at integer values:

$$\lambda(n) = \sum_{q=1}^{\infty} \nu_{2q} n^{2q} \tag{31}$$

One can easily prove that coefficients $\nu_q$ (30) with even $q = 2, 4, 6, \ldots$ have alternating signs and tend quickly to zero whereas those with odd $q = 1, 3, 5, \ldots$ are all exactly zero:

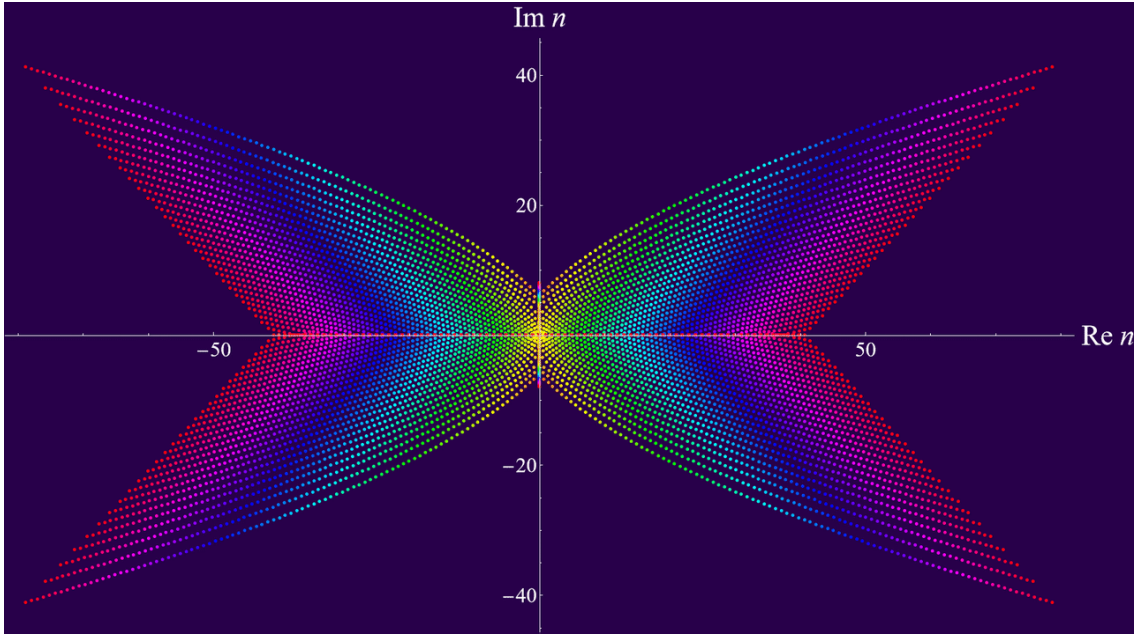

Figure 2: Distribution of complex roots of polynomials $\beta_{nk}$, labelled by the discrete parameter $k$, with $n$ as a variable, in the range of $k = 1, ..., 100$. A given color represents a specific value of $k$, e.g. the outermost red dots correspond to $k = 100$.

| $q$ | $\nu_q$ | $q$ | $\nu_q$ |
|---|---|---|---|
| 2 | $+0.0230980228342410477676...$ | 52 | $-3.7155576119842961977361... \cdot 10^{-128}$ |
| 4 | $-3.0937168341526547389855... \cdot 10^{-6}$ | 54 | $+6.4925834305243738440098... \cdot 10^{-134}$ |
| 6 | $+3.9956386236413045765970... \cdot 10^{-10}$ | 56 | $-1.0542169424270268808847... \cdot 10^{-139}$ |
| 8 | $-3.2782435703993271130506... \cdot 10^{-14}$ | 58 | $+1.5947415560818552164824... \cdot 10^{-145}$ |
| 10 | $+1.7639262984059719177765... \cdot 10^{-18}$ | 60 | $-2.2529430737360214162773... \cdot 10^{-151}$ |
| 12 | $-6.5962789103290383827162... \cdot 10^{-23}$ | 62 | $+2.9791465180379417592460... \cdot 10^{-157}$ |
| 14 | $+1.8026144481824669799319... \cdot 10^{-27}$ | 64 | $-3.6951707245044534871607... \cdot 10^{-163}$ |
| 16 | $-3.7473615439195016430971... \cdot 10^{-32}$ | 66 | $+4.3076497628345653415034... \cdot 10^{-169}$ |
| 18 | $+6.1181065540448221333456... \cdot 10^{-37}$ | 68 | $-4.7284613013367503741230... \cdot 10^{-175}$ |
| 20 | $-8.0481782096050416485754... \cdot 10^{-42}$ | 70 | $+4.8959376380361236711655... \cdot 10^{-181}$ |
| 22 | $+8.7102416117399043482412... \cdot 10^{-47}$ | 72 | $-4.7896987733528280311391... \cdot 10^{-187}$ |
| 24 | $-7.8906845737106387228314... \cdot 10^{-52}$ | 74 | $+4.4342154465175331585110... \cdot 10^{-193}$ |
| 26 | $+6.0708111424369601646918... \cdot 10^{-57}$ | 76 | $-3.8904971712441437847101... \cdot 10^{-199}$ |
| 28 | $-4.0158742275419841019323... \cdot 10^{-62}$ | 78 | $+3.2395370039463542056496... \cdot 10^{-205}$ |
| 30 | $+2.3084497121362229752239... \cdot 10^{-67}$ | 80 | $-2.5634747017952700780901... \cdot 10^{-211}$ |
| 32 | $-1.1637783963925315569620... \cdot 10^{-72}$ | 82 | $+1.9301600998675873062503... \cdot 10^{-217}$ |
| 34 | $+5.1872821372126991798104... \cdot 10^{-78}$ | 84 | $-1.3845196512412514444324... \cdot 10^{-223}$ |
| 36 | $-2.0588853472216592896166... \cdot 10^{-83}$ | 86 | $+9.4720703821567096405369... \cdot 10^{-230}$ |
| 38 | $+7.3233511971829375434125... \cdot 10^{-89}$ | 88 | $-6.1873710462545554764105... \cdot 10^{-236}$ |
| 40 | $-2.3477311629431825139375... \cdot 10^{-94}$ | 90 | $+3.8631076789948577973158... \cdot 10^{-242}$ |
| 42 | $+6.8183353226084464811575... \cdot 10^{-100}$ | 92 | $-2.3076543721182030913720... \cdot 10^{-248}$ |
| 44 | $-1.8022721811108732683536... \cdot 10^{-105}$ | 94 | $+1.3201495439226364084845... \cdot 10^{-254}$ |
| 46 | $+4.3542483564414335558191... \cdot 10^{-111}$ | 96 | $-7.2392147610413477676082... \cdot 10^{-261}$ |
| 48 | $-9.6524434850174081529763... \cdot 10^{-117}$ | 98 | $+3.8085250799294983302726... \cdot 10^{-267}$ |
| 50 | $+1.9703090445322328649401... \cdot 10^{-122}$ | 100 | $-1.9239101875299874269219... \cdot 10^{-273}$ |

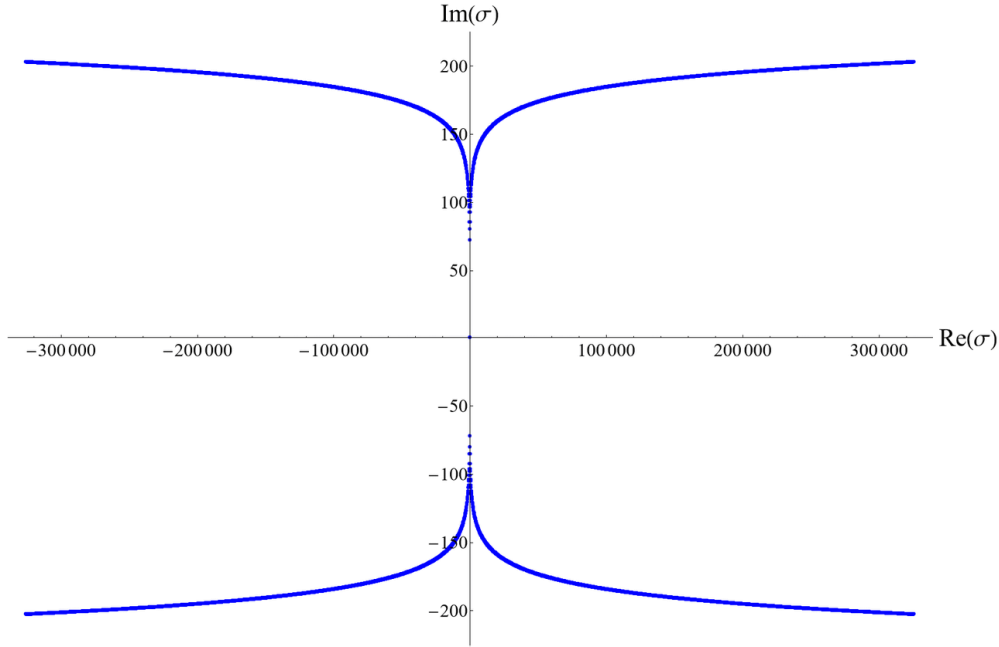

Figure 3: Distribution of 3665 zeros $\sigma_k$ of $\lambda(s)$ on the complex plane. The scale on both axes is linear but with different units.

In the summary of the considerations made so far, it can be stated that formula (31) works well for complex variable $n$. It is also clear that $\lambda_n=\lambda(n)$ is an entire even function, where $n$, as mentioned above, no longer needs to be limited to integers. Therefore, from now on, we will make change in the notation: $n \to s$.

We are quite aware that with extensions of this type there always exists a subtle problem of unambiguity: is the proposed extension unique? For instance, it is well-known that in the relatively simple case of extension of the ordinary factorial to an analytic function there are infinitely many possibilities. However, the integral proposed by Euler is the simplest, most natural.

So quite consciously, we will remain faithful to numerical experiments that effectively support intuition, in the spirit of the program presented by Jonathan Borwein. And, for the time being, we will leave the question of mathematical rigor.

## 6 Distribution of zeros

According to the general theory of entire functions [9] expression (31) has complex zeros which may be investigated numerically. To achieve this the function **FindRoot** implemented in Wolfram's Mathematica is particularly useful as well as rich graphical capabilities of this remarkable program. For example, choosing the appropriate starting point $s_0$ in

$$\textbf{FindRoot}[\lambda(s),\{s,s_0\},\textbf{WorkingPrecision} \to 100]$$

one can find several thousand complex zeros of $\lambda(s)$ (see Appendix below).

The first glance at the distribution of complex zeros of $\lambda(s)$ reveals four symmetrical branches plus one double zero in the center (Fig. 3). However, logarithmic rescaling of the real axis makes the four branches almost straight lines (Fig. 4). The mechanism of the formation of a few initial zeros is shown in Figs. 5 and 6.

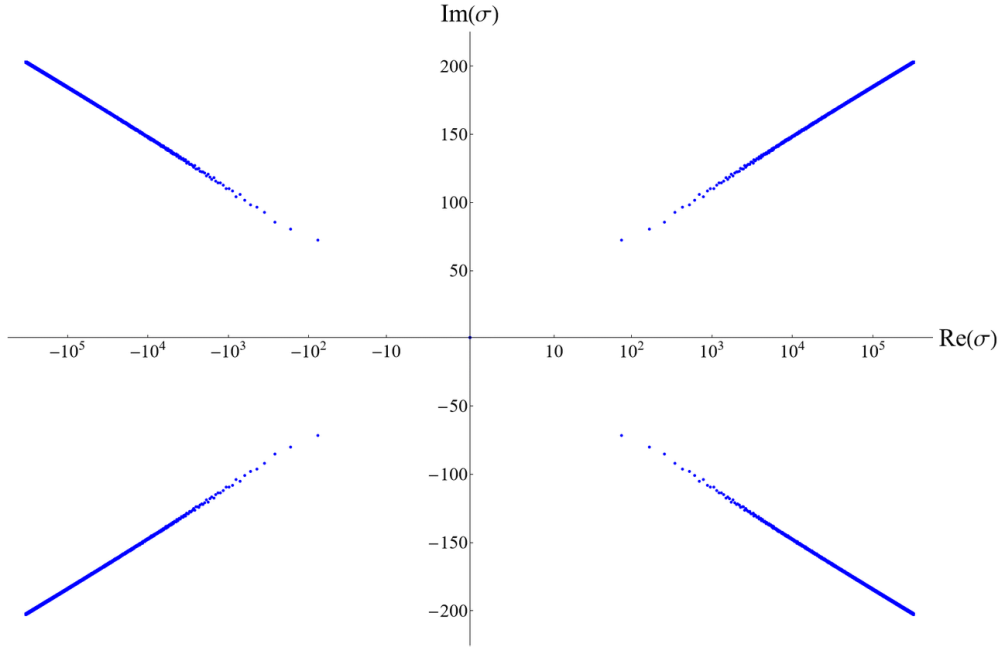

Figure 4: The same points as in the previous Figure 3, but now the horizontal axis has a logarithmic scale.

Numerical experiments clearly suggests that the growth of $\mathrm{Re}\,\sigma_k$ with $k$ is nearly linear while growth of $\mathrm{Im}\,\sigma_k$ is almost perfectly logarithmic. Combining these dependencies we obtain simple formula:

$$\mathrm{Im}\,\sigma_k \approx 16 \log(\mathrm{Re}\,\sigma_k) \tag{32}$$

(Factor 16 is from the fitting procedure.) More detailed inspection of these zeros reveals certain surprising structure in their distribution (Figs. 7 and 8). Another unexpected effect which can be grasped only by computer experiments is certain strong, long-range correlation between values of the real and imaginary part of the zeros. This can be demonstrated by simply taking high-order finite differences, e.g. of the order of 700 or more. We'll give just two examples here. These are described in details below, in Figs. 9 and 10. Many other examples of such correlations can also be given, which is currently being investigated and will be the subject of another publication.

Using Jensen's theorem, an independent test of the numerical correctness of several dozen initial zeros was carried out. By $n(r)$ we denote the number of zeros of a given integer function:

$$f(s) = \sum_{n=0}^{\infty} c_n s^n \tag{33}$$

in a disk of radius $r$, $|s| \leq r$, and by $N(r)$ the integral:

$$N(r) = \int_0^r \frac{n(t)}{t} dt \qquad \text{with } n(0) = 0 \tag{34}$$

Jensen's theorem (cf. [2], §1.2.1) relates the distribution of zeros of an entire function to its growth. It states that:

$$N(r) = \frac{1}{2\pi} \int_0^{2\pi} \ln\left|f\left(re^{i\varphi}\right)\right| d\varphi - \ln|f(0)| \tag{35}$$

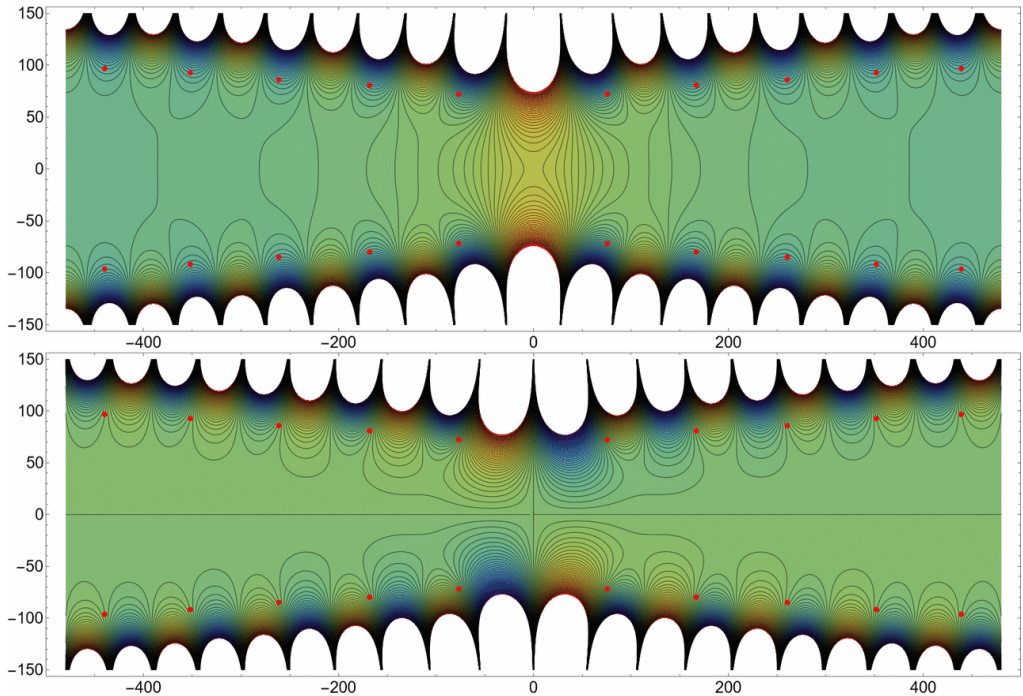

Figure 5: Contour plots of the real (upper) and imaginary (lower) part of $\lambda(x+iy)$. Complex zeros $\sigma_k$ of $\lambda(s)$ are marked by red dots.

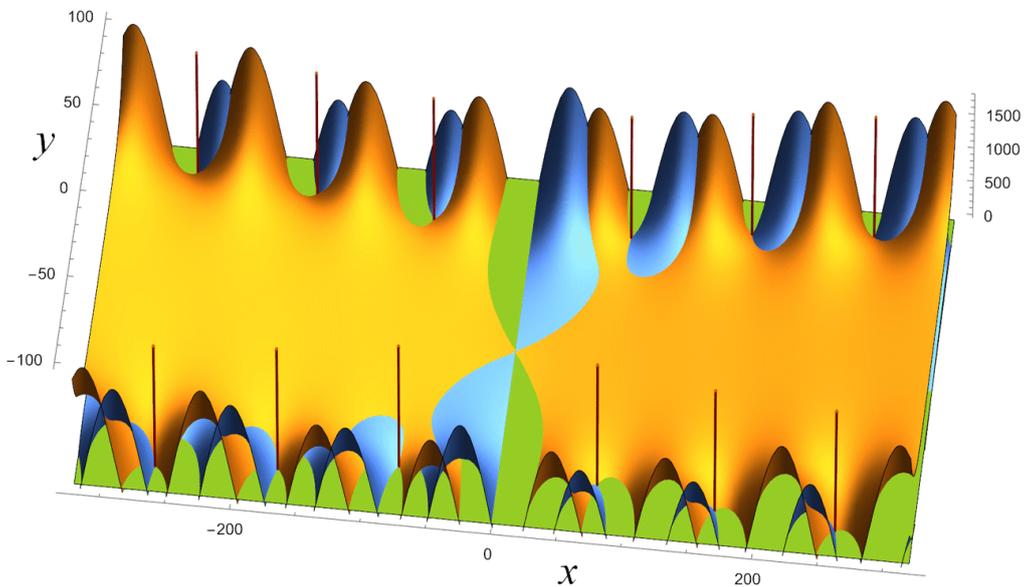

Figure 6: Three-dimensional view of the real (orange) and imaginary (blue) parts of $\lambda(s) = \lambda(x+iy)$. Both surfaces intersect the plane $(x,y)$ (green color) along certain complicated curves. These curves in turn intersect themselves in certain definite points marked by vertical dark red lines, that is, in the complex zeros of $\lambda(s)$.

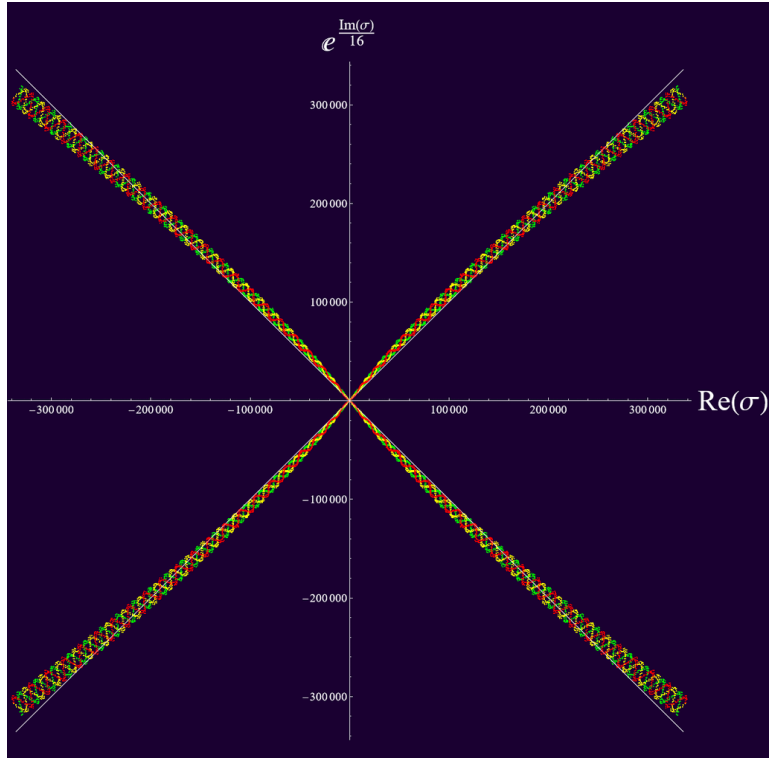

Figure 7: Distribution of 3665 complex zeros $\sigma_k$ of $\lambda(s)$. Vertical axis is rescaled as $\exp(\text{Im}(\sigma)/16)$. Amplitude of oscillations is multiplied by a factor of 5 for their better visualization.

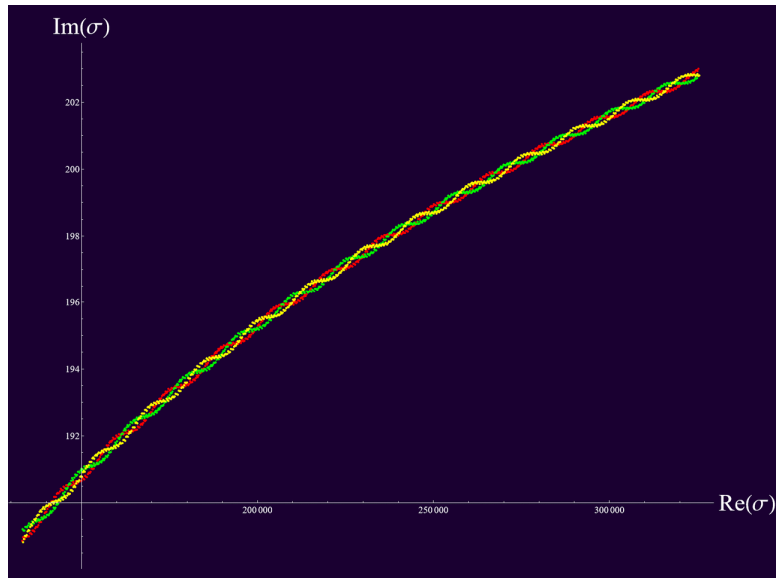

Figure 8: Unexpected distribution of complex zeros $\sigma_n$ in the range $n = 1500, ..., 3665$. It resembles something like "three-phase electric current" or – loosely speaking – a triple helix. In order to visually enhance the effect, the colors of the points were selected accordingly: zeros number $1500, 1503, 15006, ...$ are red, $1501, 1504, 15007, ...$ are green and $1502, 1505, 15008, ...$ are yellow.

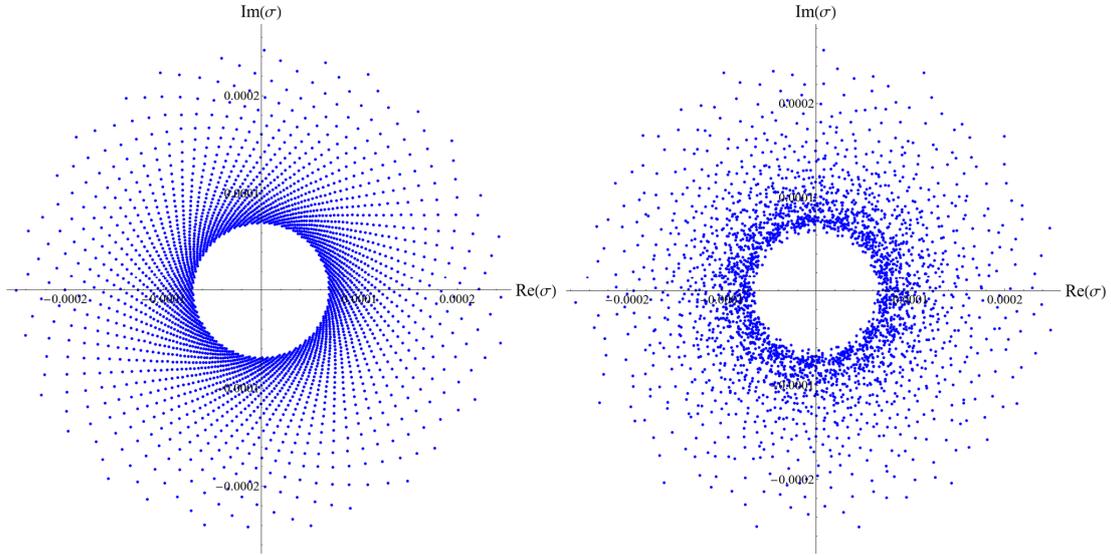

Figure 9: Another example of long-range correlations among zeros $\sigma_k$ of the complex function $\lambda(s)$. Finite differences of the order of 700 for the sequence of 3665 zeros $\sigma_k$ (left). On the right side we have the same, but a random perturbation – a complex number with an absolute value of no more than $4 \cdot 10^{-5}$ – has been added to each zero. The destruction of the regular pattern is clearly visible.

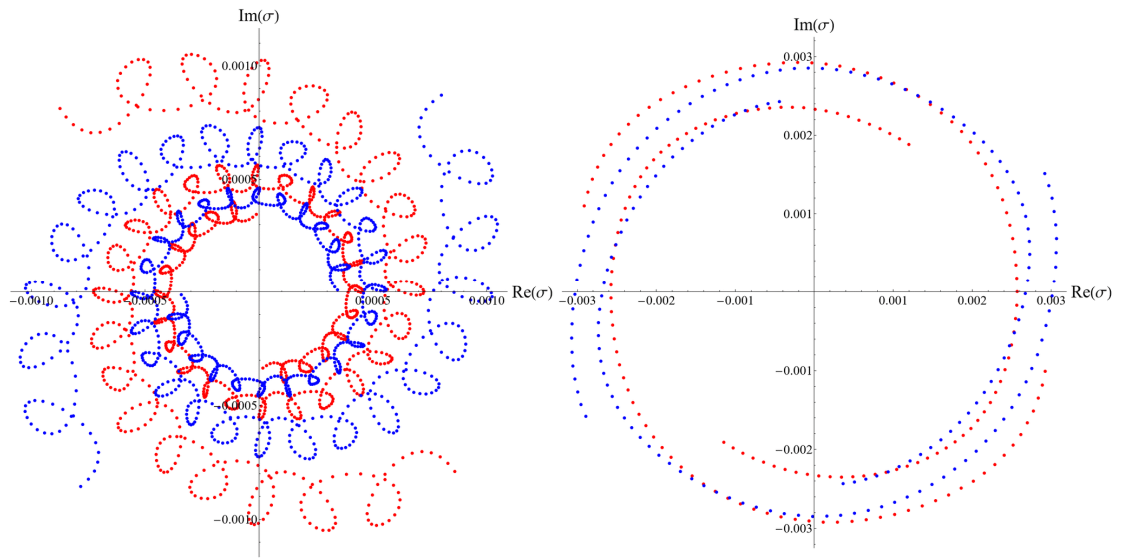

Figure 10: Finite differences of the order of 700 for the sequence of 3665 zeros $\sigma_k$ of the entire function $\lambda(s)$ but now every second zero was taken (left) or every fourth zero was taken (right).

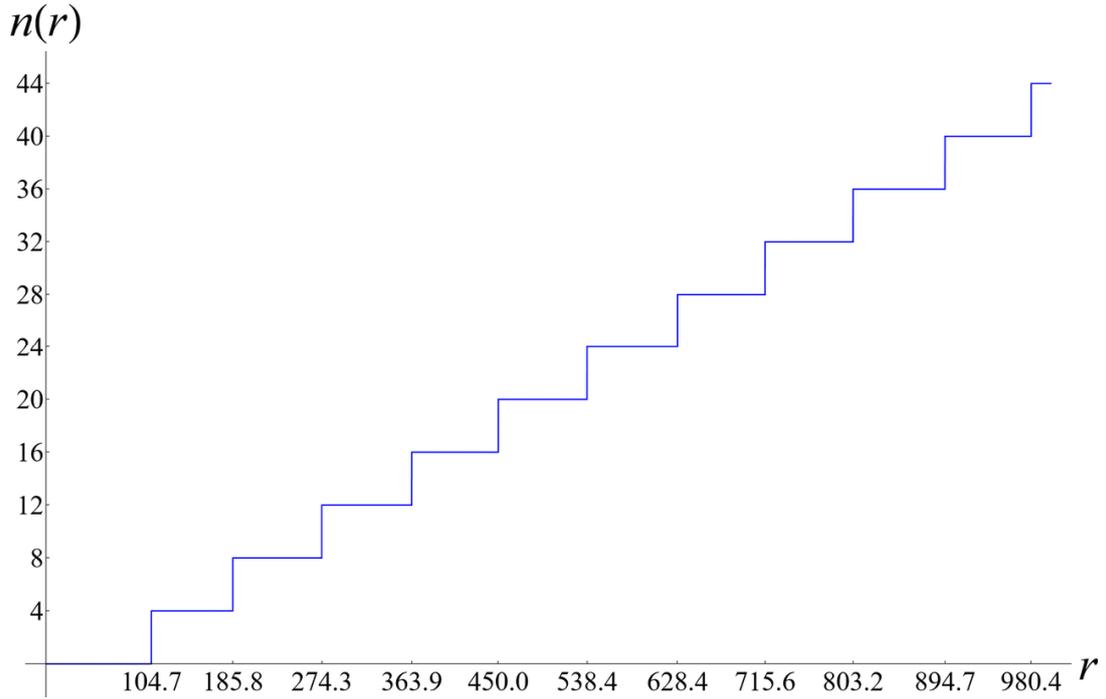

Figure 11: Whenever $r$ exceeds the absolute value of the successive zero $|\sigma_k|$, the value of $n(r)$ given by (36) jumps by 4 as expected.

By elementary calculations we get that for (33) the following relation holds:

$$n(r) = \frac{r}{2\pi}\int_0^{2\pi} h(r,\varphi)d\varphi \qquad (36)$$

where

$$h(r,\varphi) \equiv \frac{1}{2}\frac{d}{dr}\ln\left[f(r,\varphi)^2 + g(r,\varphi)^2\right] \qquad (37)$$

and

$$f(r,\varphi) \equiv \sum_{k=1}^{\infty} c_k r^{2k-2} \cos\left[2(k-1)\varphi\right] \qquad (38)$$

$$g(r,\varphi) \equiv \sum_{k=1}^{\infty} c_k r^{2k-2} \sin\left[2(k-1)\varphi\right] \qquad (39)$$

Figure 11 illustrates how the formula (36) works. (Of course, in (38) and (39) the coefficients $c_k$ have been replaced by (30). In order to satisfy condition $n(0) = 0$ the series (31) has been divided by $n^2$.)

## 7  $\lambda(s)$ as a product over its zeros

Note that the complex zeros $\sigma_k$ of the lambda function introduced in this work are in a sense "dual" to the complex zeros $\rho_n$ of the Riemann zeta function (or the xi function (1)) in the sense that both these families of zeros contain the same and complete information about this extremely important function.

Since $\lambda(s)$ is an entire function, we can represent it, by Weierstrass's theorem, as an infinite product over its zeros:

$$\lambda(s) = \text{const } s^2 \prod_{\text{all zeros } \sigma} \left(1 - \frac{s}{\sigma}\right) \quad (40)$$

(No additional prime factors are needed here because $\lambda(s)$ is of finite order [8].) Remembering that lambda zeros $\sigma$ come in fours (complex conjugate, different signs, i.e. $\sigma_k$, $-\sigma_k$, $\sigma_k^*$, $-\sigma_k^*$) we have:

$$\begin{aligned}\lambda(s) &= \text{const } s^2 \prod_{k=1}^{\infty} \left(1 - \frac{s}{\sigma_k}\right)\left(1 + \frac{s}{\sigma_k}\right)\left(1 - \frac{s}{\sigma_k^*}\right)\left(1 + \frac{s}{\sigma_k^*}\right) = \\ &= \text{const } s^2 \prod_{k=1}^{\infty} \left(1 + \frac{s^4 - 2s^2\left((\operatorname{Re}\sigma_k)^2 - (\operatorname{Im}\sigma_k)^2\right)}{|\sigma_k|^4}\right)\end{aligned} \quad (41)$$

The constant in (41) can be calculated with high accuracy. Note that it is fairly close to (but of course different from) the first Keiper-Li coefficient:

$$\begin{aligned}\text{const} &= 0.0230988022834241047767624316\ldots \\ \lambda_1 &= \frac{1}{2}(2 + \gamma - \ln 4\pi) = \\ &= 0.02309570896612103381431024790\ldots\end{aligned} \quad (42)$$

This is because the first zero is relatively large: $|\sigma_1| = 104.7027678\ldots$

Formula (41) is based on the assumption that zeros $\sigma_k$ occur in fours. At least this is the case for small values of these zeros. Of course, we do not know if such distribution of zeros will persist on the global scale. In any case, *this configuration of zeros guarantees that the Riemann hypothesis will be true*. Its violation, as we shall see, is associated with a radical deviation from the logarithmic dependence (32). Then, instead of complex fours of zeros, there will only be pairs of real zeros differing in sign.

## 8 Violation of the Riemann hypothesis?

Finally, by analyzing the formula (41) from the point of view of the arrangement of complex zeros $\sigma_k$, we can consider a hypothetical violation of the RH. As mentioned before, all attempts to use the Li's criterion (5) for this purpose are rather hopeless: in order to determine that the $k^{\text{th}}$ zero of Riemann zeta deviates from the critical line, one would have to find a negative coefficient $\lambda_n$. As estimated by Oesterle [14], such $n$ would be of the order $k^2$. Since it is known that there is no counterexample to the RH up to $k \approx 1,236 \cdot 10^{13}$ (i.e. height of the zero $\approx 3 \cdot 10^{12}$) [15], and most probably much higher, that means that all $\lambda_n$ up to $n \approx 1,528 \cdot 10^{26}$ are nonnegative. It is then clear that direct searching for the breaking of the RH using Keiper-Li criterion is a hopeless task[3].

Such a hypothetical negative coefficient $\lambda_n$ would clearly mean that some factor in the product (41) should also be negative. It is easy to check by elementary calculations that this is impossible as long as we have four symmetrical zeros: $\sigma_k$, $-\sigma_k$, $\sigma_k^*$, $-\sigma_k^*$ (cf. Figures 3 and 4). But if the logarithmic relationship $\operatorname{Im}\sigma_k \approx 16 \log(\operatorname{Re}\sigma_k)$ were broken, then – as the numerical experiments clearly suggest – zeros $\sigma_k$ would come closer to the real axis and ultimately there would be only two kinds of them: $+\sigma_k$ and $-\sigma_k$. This behavior would cause the coefficient $\lambda_n$ to become negative for some large index $n$ due to oscillations increasing with $n$. This phenomenon is illustrated in Figures 12 and 13.

Anyway, the resolution of Riemann's hypothesis is still as far away as it was when it was formulated over a century and a half ago.

---

[3]It should be stressed out that the Oesterle's estimation is only asymptotical. For example, if we consider quite absurd situation that already the first zero of zeta is off the critical line, then – depending on how much it deviates from this line – the first $\lambda_n$ that should "feel" this deviation and become negative would be for $n$ equal approx. 7500.

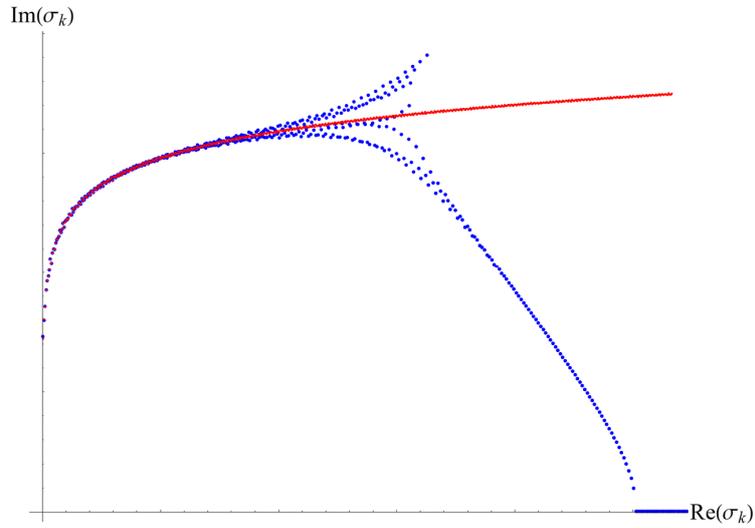

Figure 12: The distribution of zeros $\sigma_k$ in the case of the truth of the Riemann hypothesis is marked with red dots. (Only one quarter of the complex plane is visible here.) Suppose, however, that somewhere far away there is a complex zero of the Riemann zeta function $\rho_n$ that deviates from the critical line. (Of course, there would then be 4 such zeros: $\rho_n$, $1 - \rho_n$, $\rho_n^*$, $1 - \rho_n^*$.) The distribution of zeros of the function $\lambda(s)$ in this (extremely) hypothetical situation is illustrated by blue points. For the initial values of $n$ the distribution of blue points does not differ much from the distribution of red points. But as $n$ increases, there will appear oscillations, followed by a sudden drop of the zeros $\sigma_k$ to the real axis. Eventually, starting with some $k$, all zeros will line up exactly on the real axis. Then there will be only two types of them: $+\sigma_k$ and $-\sigma_k$.

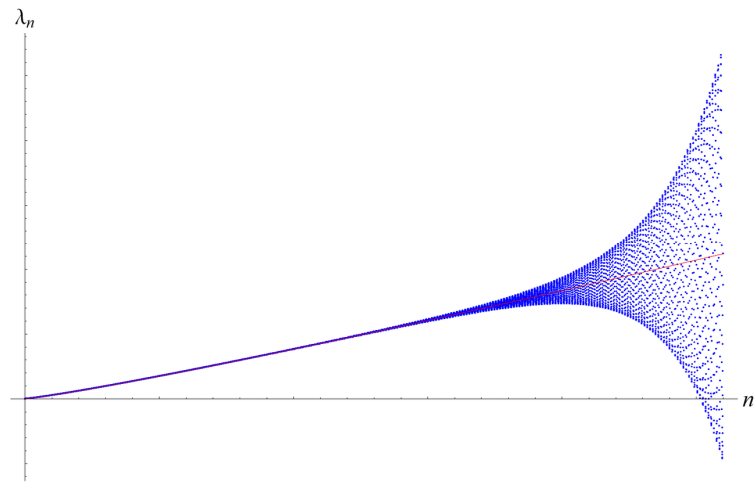

Figure 13: The figure shows the reaction of coefficients $\lambda_n$ to the behavior of zeros $\sigma_k$ caused by the deviation of even single complex zero $\rho_k$ of the Riemann zeta function from the critical line. The red dots here illustrate the distribution of $\lambda_n$ in case RH is true. (This distribution is described by an approximate law, see [12], eqn. 6). However, as $n$ increases, rising oscillations appear and eventually $\lambda_n$ becomes negative. This will happen for $n$ equal exactly to that $k$, from which the zeros $\sigma_k$ are already purely real (see previous Figure 11).

# 9 Appendix: Tables of complex zeros of $\lambda(s)$

Tables attached at the end of this article contain numerical values of the first 3520 zeros $\sigma_k$ with precision of 14 significant digits. Note that each number $x + iy$ in this table represents in fact four complex zeros $\pm x \pm iy$ located symmetrically with respect to both axes of the complex plane. There is also one double zero $\sigma_0 = 0$.

| $k$ | $\sigma_k$ | $k$ | $\sigma_k$ | $k$ | $\sigma_k$ |
|---|---|---|---|---|---|
| 1 | 76.010927161420 + 72.007003457304 $i$ | 41 | 3642.225226563 + 131.860719171 3 $i$ | 121 | 10 749.9644602275 + 149.471248288 7 $i$ |
| 2 | 167.539562791 69 + 80.244129087 91 $i$ | 42 | 3730.677895123 7 + 131.795263591 6 $i$ | 122 | 10 839.1946631452 + 148.616204350 5 $i$ |
| 3 | 260.728214212 09 + 85.357463900 80 $i$ | 43 | 3819.957050148 1 + 131.570684528 6 $i$ | 123 | 10 928.5391982744 + 149.463628115 8 $i$ |
| 4 | 352.002452273 35 + 92.402416325 40 $i$ | 44 | 3909.273772495 6 + 132.987633473 1 $i$ | 124 | 11 016.6036286267 + 149.692273657 9 $i$ |
| 5 | 439.571876677 63 + 96.441411893 07 $i$ | 45 | 3996.909139603 3 + 133.128493161 9 $i$ | 125 | 11 105.7310867881 + 149.179521871 1 $i$ |
| 6 | 529.411623473 50 + 98.056238508 06 $i$ | 46 | 4086.693481186 2 + 132.477640014 1 $i$ | 126 | 11 194.9373738856 + 149.734644118 2 $i$ |
| 7 | 620.226314479 12 + 101.152114153 58 $i$ | 47 | 4175.869073381 8 + 134.266617025 0 $i$ | 127 | 11 283.4314777224 + 150.123430359 4 $i$ |
| 8 | 707.770737608 77 + 105.535362822 04 $i$ | 48 | 4263.373927633 1 + 133.999037358 8 $i$ | 128 | 11 372.0104600234 + 149.638891719 2 $i$ |
| 9 | 796.428218702 51 + 104.084953378 42 $i$ | 49 | 4353.493568154 4 + 133.676745009 5 $i$ | 129 | 11 461.6887477782 + 149.903461828 $i$ |
| 10 | 888.114550894 48 + 108.307057371 71 $i$ | 50 | 4442.027435274 4 + 135.159512212 4 $i$ | 130 | 11 549.8573339925 + 150.731511342 8 $i$ |
| 11 | 974.251859896 19 + 109.913142910 15 $i$ | 51 | 4530.499266615 8 + 135.028903898 9 $i$ | 131 | 11 638.5700642497 + 149.806060032 1 $i$ |
| 12 | 1064.304773564 44 + 109.854011228 78 $i$ | 52 | 4619.484530425 8 + 134.810807745 3 $i$ | 132 | 11 728.2873266933 + 150.515432675 1 $i$ |
| 13 | 1154.141917026 17 + 112.360419270 19 $i$ | 53 | 4709.096661557 4 + 135.793352682 1 $i$ | 133 | 11 816.3459575779 + 150.834619907 4 $i$ |
| 14 | 1242.003593652 37 + 113.940580250 20 $i$ | 54 | 4796.729175177 9 + 136.321218153 5 $i$ | 134 | 11 905.2325217035 + 150.412203591 3 $i$ |
| 15 | 1330.670592228 6 + 114.104146617 22 $i$ | 55 | 4886.231175513 6 + 135.491361201 2 $i$ | 135 | 11 994.6588275406 + 150.702667599 3 $i$ |
| 16 | 1421.046977936 4 + 115.253818870 9 $i$ | 56 | 4975.643599519 7 + 136.745939832 6 $i$ | 136 | 12 083.0552480968 + 151.233510838 9 $i$ |
| 17 | 1508.999713386 7 + 117.783062873 4 $i$ | 57 | 5063.156938224 9 + 137.250839954 6 $i$ | 137 | 12 171.7090891669 + 150.914943996 0 $i$ |
| 18 | 1597.116048303 2 + 116.820537296 2 $i$ | 58 | 5153.020098126 0 + 136.289632242 5 $i$ | 138 | 12 261.1836991537 + 150.788578309 7 $i$ |
| 19 | 1688.146659952 9 + 118.276728821 0 $i$ | 59 | 5241.872980133 2 + 137.824207242 0 $i$ | 139 | 12 349.7231462976 + 151.850613181 5 $i$ |
| 20 | 1775.238723854 0 + 120.600671428 9 $i$ | 60 | 5330.252842206 7 + 137.310748266 6 $i$ | 140 | 12 438.0581010419 + 151.021670894 4 $i$ |
| 21 | 1864.490087819 2 + 119.130400966 0 $i$ | 61 | 5419.059475571 3 + 137.473503523 0 $i$ | 141 | 12 527.9292739115 + 151.389330873 $i$ |
| 22 | 1954.034649865 1 + 121.225888750 0 $i$ | 62 | 5508.698397266 2 + 138.236162754 8 $i$ | 142 | 12 616.1081088050 + 151.969334525 0 $i$ |
| 23 | 2042.637783450 4 + 122.281603074 9 $i$ | 63 | 5596.750951914 1 + 138.957749380 5 $i$ | 143 | 12 704.8063852375 + 151.515611863 $i$ |
| 24 | 2130.668883249 4 + 122.099680613 2 $i$ | 64 | 5685.602418009 5 + 138.007360463 4 $i$ | 144 | 12 794.2376142036 + 151.644093542 1 $i$ |
| 25 | 2220.902934155 6 + 122.770053871 4 $i$ | 65 | 5775.382096223 2 + 139.042883143 7 $i$ | 145 | 12 882.7958493931 + 152.297081828 0 $i$ |
| 26 | 2309.306514006 1 + 124.560320121 3 $i$ | 66 | 5863.122723748 7 + 139.765237501 3 $i$ | 146 | 12 971.3737565111 + 152.012623399 1 $i$ |
| 27 | 2397.125408355 0 + 124.205125195 7 $i$ | 67 | 5952.357725930 7 + 138.666330081 2 $i$ | 147 | 13 060.6636698092 + 151.757734834 8 $i$ |
| 28 | 2487.829478461 4 + 124.321692708 4 $i$ | 68 | 6041.751877328 9 + 139.964872225 6 $i$ | 148 | 13 149.5792413774 + 152.765254620 1 $i$ |
| 29 | 2575.598115341 8 + 126.804823590 3 $i$ | 69 | 6129.842431092 7 + 140.149424057 8 $i$ | 149 | 13 237.6149518046 + 152.231107774 7 $i$ |
| 30 | 2664.106663975 5 + 125.421912058 5 $i$ | 70 | 6218.823350193 8 + 139.860583690 4 $i$ | 150 | 13 327.4761517680 + 152.187615585 8 $i$ |
| 31 | 2754.087768044 1 + 126.707423810 7 $i$ | 71 | 6308.279480165 0 + 140.216275796 1 $i$ | 151 | 13 415.9072275012 + 153.030765632 $i$ |
| 32 | 2842.471441751 3 + 127.652885177 8 $i$ | 72 | 6396.613812158 0 + 141.140084345 7 $i$ | 152 | 13 504.3750167731 + 152.554296775 0 $i$ |
| 33 | 2930.656386253 2 + 127.761908799 4 $i$ | 73 | 6485.053541006 2 + 140.405097436 0 $i$ | 153 | 13 593.8291935832 + 152.590141042 0 $i$ |
| 34 | 3020.602665389 6 + 127.686269554 3 $i$ | 74 | 6575.151132354 9 + 140.919948522 5 $i$ | 154 | 13 682.5504494008 + 153.178761742 8 $i$ |
| 35 | 3109.208382503 1 + 129.302214838 0 $i$ | 75 | 6662.935524608 0 + 141.834685553 0 $i$ | 155 | 13 770.9706976325 + 153.098708701 3 $i$ |
| 36 | 3197.091561520 1 + 129.405444291 4 $i$ | 76 | 6751.781948520 6 + 140.992413101 9 $i$ | 156 | 13 860.2528172355 + 152.709929844 1 $i$ |
| 37 | 3287.341215791 5 + 128.732274941 6 $i$ | 77 | 6841.584202825 5 + 141.680378533 9 $i$ | 157 | 13 949.2804550844 + 153.527643447 6 $i$ |
| 38 | 3375.738054151 0 + 131.085562479 9 $i$ | 78 | 6929.492104908 5 + 142.317169754 2 $i$ | 158 | 14 037.3162065361 + 153.458724347 2 $i$ |
| 39 | 3463.729928244 2 + 130.203792498 2 $i$ | 79 | 7018.542545116 3 + 141.836434573 3 $i$ | 159 | 14 126.9455936565 + 152.905917492 $i$ |
| 40 | 3553.875051872 3 + 130.582785351 7 $i$ | 80 | 7107.747447008 2 + 142.094356726 $i$ | 160 | 14 215.6902614 92 + 154.002829802 $i$ |
| | | 81 | 7196.561722540 4 + 143.029017889 1 $i$ | | |
| | | 82 | 7284.596219142 8 + 142.439359466 1 $i$ | | |
| | | 83 | 7374.663462604 3 + 142.628927041 9 $i$ | | |
| | | 84 | 7462.876433166 3 + 143.718799127 7 $i$ | | |
| | | 85 | 7551.290199029 0 + 142.969980823 7 $i$ | | |
| | | 86 | 7641.204236141 5 + 143.266432765 1 $i$ | | |
| | | 87 | 7729.331300671 2 + 144.195554315 4 $i$ | | |
| | | 88 | 7818.109456761 + 143.578537186 6 $i$ | | |
| | | 89 | 7907.322920354 2 + 143.847172899 8 $i$ | | |
| | | 90 | 7996.333282719 6 + 144.610272217 2 $i$ | | |
| | | 91 | 8084.287062408 1 + 144.404164245 2 $i$ | | |
| | | 92 | 8174.161785882 7 + 144.098122436 3 $i$ | | |
| | | 93 | 8262.736232635 8 + 145.328225769 0 $i$ | | |
| | | 94 | 8350.866984486 3 + 144.825197376 1 $i$ | | |
| | | 95 | 8440.790913843 9 + 144.690702969 9 $i$ | | |
| | | 96 | 8529.113100048 2 + 145.767232200 5 $i$ | | |
| | | 97 | 8617.666740932 + 145.355247637 5 $i$ | | |
| | | 98 | 8707.050392220 6 + 145.312601439 3 $i$ | | |
| | | 99 | 8795.906287185 7 + 146.024327180 5 $i$ | | |
| | | 100 | 8884.135674555 3 + 146.176891166 2 $i$ | | |
| | | 101 | 8973.536042549 7 + 145.463255052 3 $i$ | | |
| | | 102 | 9062.614435120 1 + 146.775756452 9 $i$ | | |
| | | 103 | 9150.501722294 2 + 146.448800852 8 $i$ | | |
| | | 104 | 9240.254577523 4 + 146.107224894 3 $i$ | | |
| | | 105 | 9329.007189337 7 + 147.140894042 0 $i$ | | |
| | | 106 | 9417.226319888 4 + 146.905541351 8 $i$ | | |
| | | 107 | 9506.658737698 0 + 146.700312451 0 $i$ | | |
| | | 108 | 9595.576110659 3 + 147.370185795 7 $i$ | | |
| | | 109 | 9683.943613498 + 147.611895207 5 $i$ | | |
| | | 110 | 9772.877162944 3 + 146.939390575 8 $i$ | | |
| | | 111 | 9862.518528225 2 + 147.986122382 2 $i$ | | |
| | | 112 | 9950.183267954 8 + 147.925109381 9 $i$ | | |
| | | 113 | 10 039.6941472721 + 147.417092921 9 $i$ | | |
| | | 114 | 10 128.7963790775 + 148.357001663 4 $i$ | | |
| | | 115 | 10 216.8756407438 + 148.403140360 4 $i$ | | |
| | | 116 | 10 306.2598431987 + 147.950401636 0 $i$ | | |
| | | 117 | 10 395.2074136730 + 148.591299130 0 $i$ | | |
| | | 118 | 10 483.7297715573 + 148.969331022 5 $i$ | | |
| | | 119 | 10 572.4306093637 + 148.282485590 5 $i$ | | |
| | | 120 | 10 662.1108038245 + 148.952159471 9 $i$ | | |

| $k$ | $\sigma_k$ | $k$ | $\sigma_k$ | $k$ | $\sigma_k$ |
|---|---|---|---|---|---|
| 161 | 14 304.007501980 + 153.570176037 $i$ | 201 | 17 858.048855259 + 157.169141351 $i$ | 241 | 21 412.004448118 + 159.780774518 $i$ | 281 | 24 966.378243578 + 162.188913598 $i$
| 162 | 14 393.398344027 + 153.455711149 $i$ | 202 | 17 947.441487325 + 156.959486642 $i$ | 242 | 21 501.581508407 + 160.126582338 $i$ | 282 | 25 055.252775706 + 162.811854297 $i$
| 163 | 14 482.267001097 + 154.008317631 $i$ | 203 | 18 036.124408878 + 157.707112934 $i$ | 243 | 21 589.734604810 + 160.587439964 $i$ | 283 | 25 143.640856290 + 162.785148181 $i$
| 164 | 14 570.622289110 + 154.101308126 $i$ | 204 | 18 124.482175157 + 157.370800501 $i$ | 244 | 21 678.680662251 + 159.937188234 $i$ | 284 | 25 232.831907712 + 162.420081249 $i$
| 165 | 14 659.790469820 + 153.628464906 $i$ | 205 | 18 214.025584992 + 157.138876634 $i$ | 245 | 21 768.018060887 + 160.372798766 $i$ | 285 | 25 321.804886031 + 162.871931165 $i$
| 166 | 14 748.955798831 + 154.301139532 $i$ | 206 | 18 302.661629088 + 158.088460793 $i$ | 246 | 21 856.345076022 + 160.682089531 $i$ | 286 | 25 410.305198896 + 163.023057106 $i$
| 167 | 14 837.115146721 + 154.514372269 $i$ | 207 | 18 391.060447702 + 157.434764321 $i$ | 247 | 21 945.246557259 + 160.257056493 $i$ | 287 | 25 499.148201321 + 162.567164125 $i$
| 168 | 14 926.331513930 + 153.713887383 $i$ | 208 | 18 480.566656569 + 157.563588672 $i$ | 248 | 22 034.408200133 + 160.454502520 $i$ | 288 | 25 588.593938281 + 163.015354840 $i$
| 169 | 15 015.536050686 + 154.816376320 $i$ | 209 | 18 569.159383373 + 158.079686207 $i$ | 249 | 22 123.130189471 + 160.938290660 $i$ | 289 | 25 676.657490686 + 163.243017968 $i$
| 170 | 15 103.603053665 + 154.526887820 $i$ | 210 | 18 657.729508571 + 157.836608450 $i$ | 250 | 22 211.535306992 + 160.490932275 $i$ | 290 | 25 765.834955544 + 162.655616905 $i$
| 171 | 15 192.940512682 + 154.332311865 $i$ | 211 | 18 746.911485146 + 157.660490874 $i$ | 251 | 22 301.205451784 + 160.552195613 $i$ | 291 | 25 855.006560522 + 163.275289379 $i$
| 172 | 15 281.999779665 + 154.797639298 $i$ | 212 | 18 835.922451348 + 158.317745311 $i$ | 252 | 22 389.486898750 + 161.191218847 $i$ | 292 | 25 943.312156538 + 163.306496483 $i$
| 173 | 15 370.295925334 + 154.993742840 $i$ | 213 | 18 924.095389925 + 158.202892562 $i$ | 253 | 22 478.186680437 + 160.626866710 $i$ | 293 | 26 032.344159968 + 162.948883172 $i$
| 174 | 15 459.320529515 + 154.569343535 $i$ | 214 | 19 013.572710993 + 157.771833686 $i$ | 254 | 22 567.672237817 + 160.819816000 $i$ | 294 | 26 121.500718474 + 163.313115574 $i$
| 175 | 15 548.652064597 + 154.990788290 $i$ | 215 | 19 102.425785396 + 158.635778327 $i$ | 255 | 22 656.058191275 + 161.311606168 $i$ | 295 | 26 209.998019820 + 163.524231779 $i$
| 176 | 15 636.853837730 + 155.427239513 $i$ | 216 | 19 190.630133401 + 158.341574674 $i$ | 256 | 22 744.821323197 + 160.879539379 $i$ | 296 | 26 298.707655380 + 163.140481701 $i$
| 177 | 15 725.767186569 + 154.643857836 $i$ | 217 | 19 280.174234472 + 158.093228391 $i$ | 257 | 22 834.042061067 + 160.988499398 $i$ | 297 | 26 388.181194058 + 163.371421916 $i$
| 178 | 15 815.298743896 + 155.465018325 $i$ | 218 | 19 368.850276947 + 158.767226481 $i$ | 258 | 22 922.835251069 + 161.405893723 $i$ | 298 | 26 476.465445740 + 163.846248766 $i$
| 179 | 15 903.292295441 + 155.548971541 $i$ | 219 | 19 457.386112863 + 158.597071483 $i$ | 259 | 23 011.139695622 + 161.229791127 $i$ | 299 | 26 565.310323141 + 163.154485316 $i$
| 180 | 15 992.512277937 + 155.052331681 $i$ | 220 | 19 546.454806343 + 158.371262303 $i$ | 260 | 23 100.740685044 + 160.981240994 $i$ | 300 | 26 654.723697245 + 163.643612762 $i$
| 181 | 16 081.600652846 + 155.577162380 $i$ | 221 | 19 635.641869498 + 158.834160374 $i$ | 261 | 23 189.283325812 + 161.767918579 $i$ | 301 | 26 742.957534918 + 163.847082752 $i$
| 182 | 16 170.057933216 + 155.864449362 $i$ | 222 | 19 723.766858210 + 159.013093226 $i$ | 262 | 23 277.748238527 + 161.293084989 $i$ | 302 | 26 831.933541061 + 163.496700057 $i$
| 183 | 16 258.868175248 + 155.424143803 $i$ | 223 | 19 813.084109254 + 158.406799705 $i$ | 263 | 23 367.268089338 + 161.250571116 $i$ | 303 | 26 921.135444407 + 163.664676071 $i$
| 184 | 16 348.253935392 + 155.665809628 $i$ | 224 | 19 902.155571136 + 159.166589423 $i$ | 264 | 23 455.788805463 + 161.880700991 $i$ | 304 | 27 009.662908158 + 164.069398415 $i$
| 185 | 16 436.615109208 + 156.302304403 $i$ | 225 | 19 990.312255086 + 159.168712614 $i$ | 265 | 23 544.413538974 + 161.504264981 $i$ | 305 | 27 098.366699412 + 163.693664646 $i$
| 186 | 16 525.321433855 + 155.568624008 $i$ | 226 | 20 079.663362631 + 158.637217260 $i$ | 266 | 23 633.634662698 + 161.470910691 $i$ | 306 | 27 187.709679313 + 163.706837856 $i$
| 187 | 16 614.935878457 + 155.995938497 $i$ | 227 | 20 168.600523258 + 159.431730056 $i$ | 267 | 23 722.500751228 + 161.934750133 $i$ | 307 | 27 276.269600537 + 164.358645511 $i$
| 188 | 16 703.036687525 + 156.517008615 $i$ | 228 | 20 257.055568200 + 159.245917855 $i$ | 268 | 23 810.870211347 + 161.871359247 $i$ | 308 | 27 364.791838736 + 163.719138050 $i$
| 189 | 16 792.014738726 + 155.802158939 $i$ | 229 | 20 345.989788259 + 159.045880800 $i$ | 269 | 23 900.164148599 + 161.456680325 $i$ | 309 | 27 454.390005037 + 163.997693041 $i$
| 190 | 16 881.284791060 + 156.342608960 $i$ | 230 | 20 435.293102923 + 159.378361890 $i$ | 270 | 23 989.133765199 + 162.281976791 $i$ | 310 | 27 542.676703872 + 164.371016889 $i$
| 191 | 16 969.784496907 + 156.575121573 $i$ | 231 | 20 523.530876105 + 159.720824627 $i$ | 271 | 24 077.315135171 + 161.916082401 $i$ | 311 | 27 631.498268153 + 163.998516904 $i$
| 192 | 17 058.391543613 + 156.336874241 $i$ | 232 | 20 612.523295883 + 159.082018505 $i$ | 272 | 24 166.813608406 + 161.732992184 $i$ | 312 | 27 720.745630475 + 164.078523153 $i$
| 193 | 17 147.912991598 + 156.330126587 $i$ | 233 | 20 701.900599952 + 159.685710072 $i$ | 273 | 24 255.547287756 + 162.379849359 $i$ | 313 | 27 809.390949482 + 164.523947875 $i$
| 194 | 17 236.351019693 + 157.005045855 $i$ | 234 | 20 790.024896074 + 159.897190233 $i$ | 274 | 24 344.013721135 + 162.138315228 $i$ | 314 | 27 897.969282109 + 164.225587135 $i$
| 195 | 17 324.866623368 + 156.522706219 $i$ | 235 | 20 879.160631554 + 159.262900019 $i$ | 275 | 24 433.238571929 + 161.941293861 $i$ | 315 | 27 987.258650310 + 164.113561259 $i$
| 196 | 17 414.536694789 + 156.549925808 $i$ | 236 | 20 968.339626803 + 159.921003969 $i$ | 276 | 24 522.149720706 + 162.416400013 $i$ | 316 | 28 076.049521288 + 164.757639310 $i$
| 197 | 17 502.835286755 + 157.344203085 $i$ | 237 | 21 056.646548942 + 159.961690387 $i$ | 277 | 24 610.589974732 + 162.470170943 $i$ | 317 | 28 164.335726530 + 164.336400460 $i$
| 198 | 17 591.505221736 + 156.605011675 $i$ | 238 | 21 145.629968783 + 159.712208858 $i$ | 278 | 24 699.634765138 + 162.000770173 $i$ | 318 | 28 253.987630353 + 164.310589733 $i$
| 199 | 17 680.957378544 + 156.999452303 $i$ | 239 | 21 234.884512180 + 159.885868982 $i$ | 279 | 24 788.908593176 + 162.663174406 $i$ | 319 | 28 342.424806446 + 164.874098195 $i$
| 200 | 17 769.462277886 + 157.319558965 $i$ | 240 | 21 323.318072123 + 160.362947201 $i$ | 280 | 24 876.928351556 + 162.605492328 $i$ | 320 | 28 431.067374870 + 164.498310900 $i$

| $k$ | $\sigma_k$ | $k$ | $\sigma_k$ | $k$ | $\sigma_k$ |
|---|---|---|---|---|---|
| 321 | 28 520.361497837 + 164.493079320 $i$ | 361 | 32 074.382915360 + 166.368849273 $i$ | 401 | 35 628.345308819 + 168.339638834 $i$ |
| 322 | 28 609.089019771 + 164.912452520 $i$ | 362 | 32 162.870096823 + 166.918780714 $i$ | 402 | 35 716.587642054 + 168.410863410 $i$ |
| 323 | 28 697.588146811 + 164.784275375 $i$ | 363 | 32 251.536416417 + 166.397527453 $i$ | 403 | 35 805.722031582 + 167.957877756 $i$ |
| 324 | 28 786.812302874 + 164.510991197 $i$ | 364 | 32 340.974207914 + 166.526172276 $i$ | 404 | 35 894.816523307 + 168.454105330 $i$ |
| 325 | 28 875.759548811 + 165.152739184 $i$ | 365 | 32 429.381523665 + 167.060450367 $i$ | 405 | 35 983.208213863 + 168.459418566 $i$ |
| 326 | 28 964.029120844 + 164.937972796 $i$ | 366 | 32 518.127146671 + 166.473981766 $i$ | 406 | 36 072.199705676 + 168.229433635 $i$ |
| 327 | 29 053.467048292 + 164.592910446 $i$ | 367 | 32 607.449215688 + 166.773528034 $i$ | 407 | 36 161.368119557 + 168.393034534 $i$ |
| 328 | 29 142.204791646 + 165.366829585 $i$ | 368 | 32 695.983550986 + 167.015535579 $i$ | 408 | 36 249.843005513 + 168.740194264 $i$ |
| 329 | 29 230.654757240 + 164.975833486 $i$ | 369 | 32 784.681058482 + 166.837227628 $i$ | 409 | 36 338.610753992 + 168.219591527 $i$ |
| 330 | 29 319.939581944 + 164.916402985 $i$ | 370 | 32 873.934083927 + 166.667345720 $i$ | 410 | 36 428.028036745 + 168.569399519 $i$ |
| 331 | 29 408.789834821 + 165.299394578 $i$ | 371 | 32 962.631804197 + 167.291167437 $i$ | 411 | 36 516.287295265 + 168.855278527 $i$ |
| 332 | 29 497.249741937 + 165.304445319 $i$ | 372 | 33 051.109629758 + 166.893862151 $i$ | 412 | 36 605.258932332 + 168.319677838 $i$ |
| 333 | 29 586.358638668 + 164.926362969 $i$ | 373 | 33 140.540745370 + 166.814260693 $i$ | 413 | 36 694.471263422 + 168.721355053 $i$ |
| 334 | 29 675.452956443 + 165.492393653 $i$ | 374 | 33 229.157920199 + 167.454855478 $i$ | 414 | 36 782.916003192 + 168.891137091 $i$ |
| 335 | 29 763.723500729 + 165.485305039 $i$ | 375 | 33 317.663527351 + 166.947214910 $i$ | 415 | 36 871.774132339 + 168.549467652 $i$ |
| 336 | 29 852.922900081 + 164.988717916 $i$ | 376 | 33 407.083572537 + 167.070815783 $i$ | 416 | 36 960.955782162 + 168.730086185 $i$ |
| 337 | 29 942.011898376 + 165.741189838 $i$ | 377 | 33 495.682689586 + 167.409678615 $i$ | 417 | 37 049.602966118 + 169.059850449 $i$ |
| 338 | 30 030.228946051 + 165.488771934 $i$ | 378 | 33 584.299657957 + 167.245245140 $i$ | 418 | 37 138.145817128 + 168.689202885 $i$ |
| 339 | 30 119.536783325 + 165.334038977 $i$ | 379 | 33 673.474296354 + 167.062337419 $i$ | 419 | 37 227.676658755 + 168.758167415 $i$ |
| 340 | 30 208.464611617 + 165.656580379 $i$ | 380 | 33 762.392711663 + 167.575842205 $i$ | 420 | 37 315.995146008 + 169.253939877 $i$ |
| 341 | 30 296.920026703 + 165.798769252 $i$ | 381 | 33 850.698415561 + 167.416815890 $i$ | 421 | 37 404.807190240 + 168.737260565 $i$ |
| 342 | 30 385.917938975 + 165.380918113 $i$ | 382 | 33 940.114656753 + 167.111780095 $i$ | 422 | 37 494.121003494 + 168.947943489 $i$ |
| 343 | 30 475.117356846 + 165.790468210 $i$ | 383 | 34 028.886902276 + 167.781533277 $i$ | 423 | 37 582.605025330 + 169.282845284 $i$ |
| 344 | 30 563.436139898 + 166.013629237 $i$ | 384 | 34 117.270730543 + 167.464247509 $i$ | 424 | 37 671.361800005 + 168.925071728 $i$ |
| 345 | 30 652.427363318 + 165.418572321 $i$ | 385 | 34 206.656013169 + 167.325612274 $i$ | 425 | 37 760.566978067 + 169.023688628 $i$ |
| 346 | 30 741.726369435 + 166.030308819 $i$ | 386 | 34 295.389765766 + 167.826239038 $i$ | 426 | 37 849.305075444 + 169.362094261 $i$ |
| 347 | 30 829.895394608 + 166.047074078 $i$ | 387 | 34 383.953522634 + 167.624024863 $i$ | 427 | 37 937.766782475 + 169.146186853 $i$ |
| 348 | 30 919.090307511 + 165.676753206 $i$ | 388 | 34 473.006469032 + 167.454868596 $i$ | 428 | 38 027.212426346 + 168.963281896 $i$ |
| 349 | 31 008.107457456 + 166.063668189 $i$ | 389 | 34 562.097997346 + 167.856330991 $i$ | 429 | 38 115.764144625 + 169.639362717 $i$ |
| 350 | 31 096.634243966 + 166.229967366 $i$ | 390 | 34 650.372521711 + 167.913753325 $i$ | 430 | 38 204.369573097 + 169.139398841 $i$ |
| 351 | 31 185.460030464 + 165.863920011 $i$ | 391 | 34 739.623988018 + 167.419592323 $i$ | 431 | 38 293.728893934 + 169.188052260 $i$ |
| 352 | 31 274.775434008 + 166.084901867 $i$ | 392 | 34 828.613697452 + 168.104598044 $i$ | 432 | 38 382.320987003 + 169.644863206 $i$ |
| 353 | 31 363.154285515 + 166.473109197 $i$ | 393 | 34 916.946488327 + 167.940629990 $i$ | 433 | 38 470.969091663 + 169.302577689 $i$ |
| 354 | 31 451.945300979 + 165.903629015 $i$ | 394 | 35 006.170260646 + 167.599122557 $i$ | 434 | 38 560.149282399 + 169.298611504 $i$ |
| 355 | 31 541.391746110 + 166.296125968 $i$ | 395 | 35 095.109555569 + 168.190894200 $i$ | 435 | 38 648.990962685 + 169.675425184 $i$ |
| 356 | 31 629.624408325 + 166.563142918 $i$ | 396 | 35 183.579750107 + 168.020605353 $i$ | 436 | 38 737.428334260 + 169.559702787 $i$ |
| 357 | 31 718.590524487 + 166.063101891 $i$ | 397 | 35 272.589401387 + 167.863237446 $i$ | 437 | 38 826.713950952 + 169.257379663 $i$ |
| 358 | 31 807.810801978 + 166.457067334 $i$ | 398 | 35 361.759445562 + 168.105472595 $i$ | 438 | 38 915.584280079 + 169.908154497 $i$ |
| 359 | 31 896.319686434 + 166.583518723 $i$ | 399 | 35 450.077890290 + 168.365804414 $i$ | 439 | 39 003.899684137 + 169.561970983 $i$ |
| 360 | 31 985.022254132 + 166.369870117 $i$ | 400 | 35 539.129648587 + 167.808849618 $i$ | 440 | 39 093.333383578 + 169.455072647 $i$ |
| 441 | 39 182.025478480 + 169.944692973 $i$ | | | | |
| 442 | 39 270.582475189 + 169.720817010 $i$ | | | | |
| 443 | 39 359.745032766 + 169.571430746 $i$ | | | | |
| 444 | 39 448.661205109 + 169.968335472 $i$ | | | | |
| 445 | 39 537.117907101 + 169.942194609 $i$ | | | | |
| 446 | 39 626.213215505 + 169.579327920 $i$ | | | | |
| 447 | 39 715.334946647 + 170.125659241 $i$ | | | | |
| 448 | 39 803.513911542 + 170.029633877 $i$ | | | | |
| 449 | 39 892.899736403 + 169.693828959 $i$ | | | | |
| 450 | 39 981.721113758 + 170.244239773 $i$ | | | | |
| 451 | 40 070.237031048 + 170.114403143 $i$ | | | | |
| 452 | 40 159.302522692 + 169.854532642 $i$ | | | | |
| 453 | 40 248.332427823 + 170.262394399 $i$ | | | | |
| 454 | 40 336.822781407 + 170.290247186 $i$ | | | | |
| 455 | 40 425.723120555 + 169.930079872 $i$ | | | | |
| 456 | 40 515.031541131 + 170.332178123 $i$ | | | | |
| 457 | 40 603.214574676 + 170.459315812 $i$ | | | | |
| 458 | 40 692.397880735 + 169.967064646 $i$ | | | | |
| 459 | 40 781.455171744 + 170.510300143 $i$ | | | | |
| 460 | 40 869.868149317 + 170.472306615 $i$ | | | | |
| 461 | 40 958.881806986 + 170.190001480 $i$ | | | | |
| 462 | 41 047.986051173 + 170.484242878 $i$ | | | | |
| 463 | 41 136.492064619 + 170.662183620 $i$ | | | | |
| 464 | 41 225.310619827 + 170.297104542 $i$ | | | | |
| 465 | 41 314.646508611 + 170.509980630 $i$ | | | | |
| 466 | 41 402.956100616 + 170.868312242 $i$ | | | | |
| 467 | 41 491.897084551 + 170.299979613 $i$ | | | | |
| 468 | 41 581.180347396 + 170.717664631 $i$ | | | | |
| 469 | 41 669.501390583 + 170.837609146 $i$ | | | | |
| 470 | 41 758.483516680 + 170.528673341 $i$ | | | | |
| 471 | 41 847.603322216 + 170.702305766 $i$ | | | | |
| 472 | 41 936.196930451 + 171.012279780 $i$ | | | | |
| 473 | 42 024.894922173 + 170.640862212 $i$ | | | | |
| 474 | 42 114.223655061 + 170.747328626 $i$ | | | | |
| 475 | 42 202.759979375 + 171.183742284 $i$ | | | | |
| 476 | 42 291.371638062 + 170.656136823 $i$ | | | | |
| 477 | 42 380.850507898 + 170.920268984 $i$ | | | | |
| 478 | 42 469.196015221 + 171.195489249 $i$ | | | | |
| 479 | 42 558.062029716 + 170.848356876 $i$ | | | | |
| 480 | 42 647.209776449 + 170.956753705 $i$ | | | | |

| $k$ | $\sigma_k$ | $k$ | $\sigma_k$ | $k$ | $\sigma_k$ |
|---|---|---|---|---|---|
| 481 | 42 735.911458916 + 171.314760974 $i$ | 521 | 46 289.650276259 + 172.580864889 $i$ | 561 | 49 843.476696929 + 173.511595503 $i$ |
| 482 | 42 824.502973863 + 171.006575390 $i$ | 522 | 46 378.518722988 + 172.087088910 $i$ | 562 | 49 932.682148275 + 173.310569053 $i$ |
| 483 | 42 913.788866224 + 170.962919617 $i$ | 523 | 46 467.842609405 + 172.463976177 $i$ | 563 | 50 021.603730408 + 173.752244122 $i$ |
| 484 | 43 002.496398389 + 171.472707525 $i$ | 524 | 46 556.144895406 + 172.626842212 $i$ | 564 | 50 110.078828813 + 173.613773285 $i$ |
| 485 | 43 090.949951563 + 171.083349133 $i$ | 525 | 46 645.142232064 + 172.217916027 $i$ | 565 | 50 199.140488021 + 173.466858654 $i$ |
| 486 | 43 180.459039006 + 171.072135467 $i$ | 526 | 46 734.278494419 + 172.541245008 $i$ | 566 | 50 288.221075841 + 173.690105893 $i$ |
| 487 | 43 268.905811152 + 171.551844688 $i$ | 527 | 46 822.810733846 + 172.667108263 $i$ | 567 | 50 376.603827256 + 173.874894018 $i$ |
| 488 | 43 357.650850191 + 171.182265016 $i$ | 528 | 46 911.600345382 + 172.434131434 $i$ | 568 | 50 465.655767901 + 173.377742788 $i$ |
| 489 | 43 446.828485884 + 171.201240896 $i$ | 529 | 47 000.855522164 + 172.448349487 $i$ | 569 | 50 554.796835565 + 173.914337735 $i$ |
| 490 | 43 535.588673807 + 171.572082154 $i$ | 530 | 47 089.374006524 + 172.922451717 $i$ | 570 | 50 643.124996893 + 173.867975533 $i$ |
| 491 | 43 624.140838147 + 171.396305475 $i$ | 531 | 47 178.090972848 + 172.401462198 $i$ | 571 | 50 732.238820083 + 173.524842073 $i$ |
| 492 | 43 713.340659992 + 171.184477272 $i$ | 532 | 47 267.439956829 + 172.617775208 $i$ | 572 | 50 821.279498795 + 173.962371210 $i$ |
| 493 | 43 802.221268828 + 171.739204467 $i$ | 533 | 47 355.888946344 + 172.978242380 $i$ | 573 | 50 909.751848014 + 173.933647769 $i$ |
| 494 | 43 890.588592973 + 171.473001707 $i$ | 534 | 47 444.675994707 + 172.497203188 $i$ | 574 | 50 998.725549494 + 173.692762121 $i$ |
| 495 | 43 979.968054983 + 171.266548536 $i$ | 535 | 47 533.936817593 + 172.754772426 $i$ | 575 | 51 087.830282219 + 173.895492806 $i$ |
| 496 | 44 068.691733031 + 171.869007851 $i$ | 536 | 47 622.479379463 + 172.945639112 $i$ | 576 | 51 176.348256687 + 174.160581798 $i$ |
| 497 | 44 157.207048935 + 171.503740732 $i$ | 537 | 47 711.215722346 + 172.759587425 $i$ | 577 | 51 265.163823960 + 173.671145051 $i$ |
| 498 | 44 246.445421885 + 171.484150021 $i$ | 538 | 47 800.419492267 + 172.662096595 $i$ | 578 | 51 354.482554674 + 174.041959377 $i$ |
| 499 | 44 335.270569256 + 171.795840254 $i$ | 539 | 47 889.138498647 + 173.179767369 $i$ | 579 | 51 442.799446685 + 174.220403718 $i$ |
| 500 | 44 423.784059198 + 171.769749653 $i$ | 540 | 47 977.648476698 + 172.752299203 $i$ | 580 | 51 531.803554582 + 173.763510385 $i$ |
| 501 | 44 512.891523034 + 171.454247194 $i$ | 541 | 48 067.039318130 + 172.791021599 $i$ | 581 | 51 620.935936885 + 174.135078687 $i$ |
| 502 | 44 601.938519634 + 171.955099146 $i$ | 542 | 48 155.635664652 + 173.258008456 $i$ | 582 | 51 709.427766708 + 174.238732000 $i$ |
| 503 | 44 690.246114497 + 171.848711906 $i$ | 543 | 48 244.225720350 + 172.823373529 $i$ | 583 | 51 798.299390409 + 173.953877675 $i$ |
| 504 | 44 779.466186891 + 171.515212318 $i$ | 544 | 48 333.565171465 + 172.940164919 $i$ | 584 | 51 887.467891496 + 174.100179972 $i$ |
| 505 | 44 868.464818494 + 172.107750948 $i$ | 545 | 48 422.167397323 + 173.229822636 $i$ | 585 | 51 976.052161722 + 174.383441179 $i$ |
| 506 | 44 956.791646532 + 171.873937067 $i$ | 546 | 48 510.844864437 + 173.050975167 $i$ | 586 | 52 064.719064561 + 174.040727294 $i$ |
| 507 | 45 046.057837159 + 171.726075333 $i$ | 547 | 48 599.972970210 + 172.900057947 $i$ | 587 | 52 154.138253963 + 174.126056662 $i$ |
| 508 | 45 134.916071350 + 172.036143775 $i$ | 548 | 48 688.856541297 + 173.385491055 $i$ | 588 | 52 242.491955406 + 174.540959899 $i$ |
| 509 | 45 223.476824439 + 172.117044558 $i$ | 549 | 48 777.269329179 + 173.157541883 $i$ | 589 | 52 331.363755943 + 174.034865277 $i$ |
| 510 | 45 312.434205979 + 171.726750692 $i$ | 550 | 48 866.618430765 + 172.927609232 $i$ | 590 | 52 420.573002322 + 174.293680613 $i$ |
| 511 | 45 401.596223705 + 172.155077742 $i$ | 551 | 48 955.340216219 + 173.520135167 $i$ | 591 | 52 509.119452542 + 174.529905134 $i$ |
| 512 | 45 489.952624151 + 172.237828095 $i$ | 552 | 49 043.843306217 + 173.175878988 $i$ | 592 | 52 597.883717152 + 174.217772919 $i$ |
| 513 | 45 578.985039916 + 171.775394493 $i$ | 553 | 49 133.130366284 + 173.106259345 $i$ | 593 | 52 687.072469327 + 174.307501362 $i$ |
| 514 | 45 668.176110669 + 172.303257430 $i$ | 554 | 49 221.877941622 + 173.526340281 $i$ | 594 | 52 775.769916668 + 174.612238855 $i$ |
| 515 | 45 756.450156989 + 172.256414858 $i$ | 555 | 49 310.476884578 + 173.325840275 $i$ | 595 | 52 864.326341861 + 174.361015374 $i$ |
| 516 | 45 845.611750634 + 171.950287188 $i$ | 556 | 49 399.543956229 + 173.176252046 $i$ | 596 | 52 953.683780842 + 174.259439321 $i$ |
| 517 | 45 934.583853706 + 172.301433666 $i$ | 557 | 49 488.559402499 + 173.540430136 $i$ | 597 | 53 042.261594097 + 174.836308023 $i$ |
| 518 | 46 023.153646492 + 172.400413123 $i$ | 558 | 49 576.912584801 + 173.531986554 $i$ | 598 | 53 130.916009556 + 174.316287222 $i$ |
| 519 | 46 111.989411657 + 172.080759035 $i$ | 559 | 49 666.142451461 + 173.130509015 $i$ | 599 | 53 220.195191118 + 174.458135553 $i$ |
| 520 | 46 201.263963634 + 172.310772290 $i$ | 560 | 49 755.080335784 + 173.748295878 $i$ | 600 | 53 308.821219003 + 174.792472649 $i$ |
| $k$ | $\sigma_k$ | | | | |
| 601 | 53 397.488610976 + 174.493815305 $i$ | | | | |
| 602 | 53 486.661250440 + 174.498595024 $i$ | | | | |
| 603 | 53 575.466366433 + 174.828308178 $i$ | | | | |
| 604 | 53 663.958452580 + 174.677139627 $i$ | | | | |
| 605 | 53 753.226946207 + 174.448060046 $i$ | | | | |
| 606 | 53 842.032222920 + 175.021894991 $i$ | | | | |
| 607 | 53 930.462065717 + 174.663904572 $i$ | | | | |
| 608 | 54 019.823951294 + 174.603107178 $i$ | | | | |
| 609 | 54 108.492357278 + 175.030512132 $i$ | | | | |
| 610 | 54 197.129186012 + 174.797561450 $i$ | | | | |
| 611 | 54 286.235187460 + 174.674424917 $i$ | | | | |
| 612 | 54 375.142162502 + 175.047902145 $i$ | | | | |
| 613 | 54 463.627078714 + 174.973755657 $i$ | | | | |
| 614 | 54 552.753644106 + 174.671684197 $i$ | | | | |
| 615 | 54 641.774714066 + 175.172794021 $i$ | | | | |
| 616 | 54 730.076358479 + 175.023103993 $i$ | | | | |
| 617 | 54 819.395032597 + 174.744944996 $i$ | | | | |
| 618 | 54 908.186784437 + 175.261821195 $i$ | | | | |
| 619 | 54 996.768543350 + 175.075450097 $i$ | | | | |
| 620 | 55 085.798149751 + 174.899249023 $i$ | | | | |
| 621 | 55 174.840596066 + 175.237086553 $i$ | | | | |
| 622 | 55 263.296394985 + 175.230735702 $i$ | | | | |
| 623 | 55 352.280133126 + 174.927905858 $i$ | | | | |
| 624 | 55 441.472344740 + 175.298649120 $i$ | | | | |
| 625 | 55 529.739558574 + 175.369684511 $i$ | | | | |
| 626 | 55 618.931276246 + 174.938487554 $i$ | | | | |
| 627 | 55 707.911232705 + 175.440701534 $i$ | | | | |
| 628 | 55 796.388466049 + 175.356825839 $i$ | | | | |
| 629 | 55 885.399074340 + 175.132558417 $i$ | | | | |
| 630 | 55 974.472340564 + 175.385286166 $i$ | | | | |
| 631 | 56 062.981905621 + 175.530294394 $i$ | | | | |
| 632 | 56 151.851307521 + 175.178006069 $i$ | | | | |
| 633 | 56 241.106078911 + 175.429103482 $i$ | | | | |
| 634 | 56 329.468981358 + 175.678891456 $i$ | | | | |
| 635 | 56 418.435762964 + 175.166679739 $i$ | | | | |
| 636 | 56 507.620632416 + 175.591387449 $i$ | | | | |
| 637 | 56 596.031107886 + 175.648609925 $i$ | | | | |
| 638 | 56 685.003683824 + 175.357300689 $i$ | | | | |
| 639 | 56 774.077986871 + 175.553803621 $i$ | | | | |
| 640 | 56 862.708859112 + 175.791434582 $i$ | | | | |

| $k$ | $\sigma_k$ | $k$ | $\sigma_k$ | $k$ | $\sigma_k$ |
|---|---|---|---|---|---|
| 641 | 56 951.413727674 + 175.424987795 $i$ | 681 | 60 505.528932610 + 176.297984297 $i$ | 721 | 64 059.617093774 + 177.330492519 $i$ | 761 | 67 613.558797476 + 178.239479353 $i$
| 642 | 57 040.713269022 + 175.582199027 $i$ | 682 | 60 594.626242106 + 176.756427577 $i$ | 722 | 64 148.367912814 + 177.674524631 $i$ | 762 | 67 702.228480785 + 178.546273257 $i$
| 643 | 57 129.215480828 + 175.924422775 $i$ | 683 | 60 682.967635364 + 176.682676049 $i$ | 723 | 64 236.975260055 + 177.479920430 $i$ | 763 | 67 790.873262565 + 178.256028622 $i$
| 644 | 57 217.947194606 + 175.456463706 $i$ | 684 | 60 772.120059289 + 176.429932191 $i$ | 724 | 64 326.053177692 + 177.360899045 $i$ | 764 | 67 880.162584107 + 178.213867470 $i$
| 645 | 57 307.304729209 + 175.699840804 $i$ | 685 | 60 861.069826665 + 176.752971614 $i$ | 725 | 64 415.023267597 + 177.713275830 $i$ | 765 | 67 968.737967216 + 178.711111355 $i$
| 646 | 57 395.682596841 + 175.939434922 $i$ | 686 | 60 949.658459250 + 176.797676836 $i$ | 726 | 64 503.437064688 + 177.640272955 $i$ | 766 | 68 057.439044145 + 178.212741390 $i$
| 647 | 57 484.609751280 + 175.605527411 $i$ | 687 | 61 038.506750708 + 176.511699277 $i$ | 727 | 64 592.650603942 + 177.311027428 $i$ | 767 | 68 146.688342095 + 178.395132719 $i$
| 648 | 57 573.695354785 + 175.715128503 $i$ | 688 | 61 127.738635731 + 176.753637042 $i$ | 728 | 64 681.542063342 + 177.883614169 $i$ | 768 | 68 235.296315626 + 178.637765671 $i$
| 649 | 57 662.401021236 + 176.014304866 $i$ | 689 | 61 216.133003175 + 176.957525489 $i$ | 729 | 64 769.996457184 + 177.606539806 $i$ | 769 | 68 324.003124536 + 178.395364587 $i$
| 650 | 57 751.015626271 + 175.721953914 $i$ | 690 | 61 305.071584364 + 176.504365992 $i$ | 730 | 64 859.194102110 + 177.479816787 $i$ | 770 | 68 413.167048249 + 178.381250579 $i$
| 651 | 57 840.304154756 + 175.712608007 $i$ | 691 | 61 394.288540930 + 176.880889827 $i$ | 731 | 64 948.064453465 + 177.841844908 $i$ | 771 | 68 501.926273290 + 178.702710935 $i$
| 652 | 57 928.939990362 + 176.150015886 $i$ | 692 | 61 482.656105731 + 176.981686696 $i$ | 732 | 65 036.592266597 + 177.736162812 $i$ | 772 | 68 590.484910886 + 178.520604268 $i$
| 653 | 58 017.518910447 + 175.771304457 $i$ | 693 | 61 571.674985793 + 176.616216025 $i$ | 733 | 65 125.664894606 + 177.555447244 $i$ | 773 | 68 679.732805302 + 178.340013084 $i$
| 654 | 58 106.918206795 + 175.795642709 $i$ | 694 | 61 660.737602038 + 176.915019452 $i$ | 734 | 65 214.664873390 + 177.820893805 $i$ | 774 | 68 768.483683468 + 178.862223179 $i$
| 655 | 58 195.397169950 + 176.214798871 $i$ | 695 | 61 749.310258258 + 177.034477183 $i$ | 735 | 65 303.129061903 + 177.936168135 $i$ | 775 | 68 856.998590837 + 178.493440783 $i$
| 656 | 58 284.185207874 + 175.837357042 $i$ | 696 | 61 838.133359881 + 176.771032097 $i$ | 736 | 65 392.176486609 + 177.475797542 $i$ | 776 | 68 946.306529900 + 178.483259359 $i$
| 657 | 58 373.306667456 + 175.909011105 $i$ | 697 | 61 927.320845579 + 176.841565422 $i$ | 737 | 65 481.255747079 + 178.012377050 $i$ | 777 | 69 034.970986625 + 178.857409247 $i$
| 658 | 58 462.098009995 + 176.209106949 $i$ | 698 | 62 015.868779182 + 177.240200401 $i$ | 738 | 65 569.636546046 + 177.893553718 $i$ | 778 | 69 123.650741389 + 178.610821243 $i$
| 659 | 58 550.645060857 + 176.004564322 $i$ | 699 | 62 104.620797340 + 176.728043515 $i$ | 739 | 65 658.761361218 + 177.628468868 $i$ | 779 | 69 212.723811201 + 178.526083790 $i$
| 660 | 58 639.844548403 + 175.870717352 $i$ | 700 | 62 193.919223624 + 177.002029564 $i$ | 740 | 65 747.739587803 + 178.005271400 $i$ | 780 | 69 301.629823996 + 178.872841582 $i$
| 661 | 58 728.687893225 + 176.363376524 $i$ | 701 | 62 282.375878159 + 177.241467312 $i$ | 741 | 65 836.253101621 + 177.978979433 $i$ | 781 | 69 390.128273579 + 178.754364464 $i$
| 662 | 58 817.127128501 + 176.054132649 $i$ | 702 | 62 371.211171149 + 176.840651512 $i$ | 742 | 65 925.239332127 + 177.746648170 $i$ | 782 | 69 479.271988079 + 178.505826481 $i$
| 663 | 58 906.461835074 + 175.938785459 $i$ | 703 | 62 460.412888596 + 177.064395884 $i$ | 743 | 66 014.311076742 + 177.962898127 $i$ | 783 | 69 568.216513525 + 178.984841158 $i$
| 664 | 58 995.159819997 + 176.443164974 $i$ | 704 | 62 548.968605706 + 177.241013753 $i$ | 744 | 66 102.833129708 + 178.160768172 $i$ | 784 | 69 656.611196807 + 178.791044286 $i$
| 665 | 59 083.744417798 + 176.097783799 $i$ | 705 | 62 637.744070946 + 177.028755754 $i$ | 745 | 66 191.701509027 + 177.719362028 $i$ | 785 | 69 745.887782925 + 178.575908915 $i$
| 666 | 59 172.942317190 + 176.091543665 $i$ | 706 | 62 726.900830603 + 176.983204570 $i$ | 746 | 66 280.956721058 + 178.083435276 $i$ | 786 | 69 834.668431678 + 179.046071915 $i$
| 667 | 59 261.745616126 + 176.378524107 $i$ | 707 | 62 815.612787723 + 177.440094421 $i$ | 747 | 66 369.287421940 + 178.183605437 $i$ | 787 | 69 923.277211403 + 178.823464347 $i$
| 668 | 59 350.311411915 + 176.313742278 $i$ | 708 | 62 904.177086416 + 177.012037712 $i$ | 748 | 66 458.329323475 + 177.806361875 $i$ | 788 | 70 012.299233990 + 178.708065611 $i$
| 669 | 59 439.399437524 + 176.033525196 $i$ | 709 | 62 993.540773126 + 177.101441724 $i$ | 749 | 66 547.408492829 + 178.145159285 $i$ | 789 | 70 101.317685007 + 178.998660795 $i$
| 670 | 59 528.392718834 + 176.524982328 $i$ | 710 | 63 082.087462193 + 177.474044705 $i$ | 750 | 66 635.924848416 + 178.204716503 $i$ | 790 | 70 189.780871216 + 178.995910767 $i$
| 671 | 59 616.766833826 + 176.364072694 $i$ | 711 | 63 170.780592917 + 177.094580224 $i$ | 751 | 66 724.809587934 + 177.957067093 $i$ | 791 | 70 278.835846090 + 178.689538824 $i$
| 672 | 59 706.004582651 + 176.099919071 $i$ | 712 | 63 260.040812920 + 177.187156983 $i$ | 752 | 66 813.954655354 + 178.101201397 $i$ | 792 | 70 367.902474825 + 179.061232004 $i$
| 673 | 59 794.897614499 + 176.611745263 $i$ | 713 | 63 348.644294443 + 177.467700877 $i$ | 753 | 66 902.528144131 + 178.355344954 $i$ | 793 | 70 456.256719085 + 179.092336425 $i$
| 674 | 59 883.333043051 + 176.399344175 $i$ | 714 | 63 437.371657412 + 177.150379435 $i$ | 754 | 66 991.271174351 + 177.987722897 $i$ | 794 | 70 545.444256710 + 178.707186942 $i$
| 675 | 59 972.556336962 + 176.249633687 $i$ | 715 | 63 526.469222501 + 177.150379435 $i$ | 755 | 67 080.578979520 + 178.132532399 $i$ | 795 | 70 634.372269238 + 179.189584389 $i$
| 676 | 60 061.382556324 + 176.570088569 $i$ | 716 | 63 615.328297368 + 177.602117568 $i$ | 756 | 67 168.994423953 + 178.484336052 $i$ | 796 | 70 722.903472474 + 179.065845331 $i$
| 677 | 60 150.002237399 + 176.586613774 $i$ | 717 | 63 703.793476518 + 177.328491472 $i$ | 757 | 67 257.897764575 + 177.979285053 $i$ | 797 | 70 811.911835770 + 178.875440719 $i$
| 678 | 60 238.943352312 + 176.243185061 $i$ | 718 | 63 793.112438740 + 177.186002425 $i$ | 758 | 67 347.035622610 + 178.279691587 $i$ | 798 | 70 900.944826619 + 179.107474227 $i$
| 679 | 60 328.085359601 + 176.658706619 $i$ | 719 | 63 881.810214403 + 177.702364822 $i$ | 759 | 67 435.622968917 + 178.440816916 $i$ | 799 | 70 989.467756184 + 179.238819873 $i$
| 680 | 60 416.436358822 + 176.658091909 $i$ | 720 | 63 970.382582962 + 177.340092762 $i$ | 760 | 67 524.401927745 + 178.157442893 $i$ | 800 | 71 078.384758482 + 178.885556538 $i$

| $k$ | $\sigma_k$ | $k$ | $\sigma_k$ | $k$ | $\sigma_k$ |
|---|---|---|---|---|---|
| 801 | 71 167.569056125 + 179.167540812 $i$ | 841 | 74 721.348525927 + 180.105313773 $i$ | 881 | 78 275.130859169 + 180.801993059 $i$ | 921 | 81 829.000368674 + 181.518892147 $i$ |
| 802 | 71 255.963713138 + 179.334674753 $i$ | 842 | 74 809.865373470 + 179.883375229 $i$ | 882 | 78 363.882598349 + 180.567587733 $i$ | 922 | 81 917.808423671 + 181.133975621 $i$ |
| 803 | 71 344.958647730 + 178.883701508 $i$ | 843 | 74 899.036764095 + 179.757401705 $i$ | 883 | 78 452.965851872 + 180.499475872 $i$ | 923 | 82 007.040493383 + 181.327146194 $i$ |
| 804 | 71 434.080341923 + 179.300551664 $i$ | 844 | 74 987.869084632 + 180.080697988 $i$ | 884 | 78 541.801541869 + 180.908472419 $i$ | 924 | 82 095.487051237 + 181.603590675 $i$ |
| 805 | 71 522.536821627 + 179.307036735 $i$ | 845 | 75 076.500680985 + 180.031867751 $i$ | 885 | 78 630.305303188 + 180.611306867 $i$ | 925 | 82 184.411537498 + 181.136844675 $i$ |
| 806 | 71 611.513311656 + 179.050909644 $i$ | 846 | 75 165.444180736 + 179.755131548 $i$ | 886 | 78 719.614510375 + 180.533140353 $i$ | 926 | 82 273.517658776 + 181.447079989 $i$ |
| 807 | 71 700.563180779 + 179.244794015 $i$ | 847 | 75 254.565391483 + 180.148408807 $i$ | 887 | 78 808.272333918 + 180.979506219 $i$ | 927 | 82 362.098713881 + 181.549979570 $i$ |
| 808 | 71 789.188933760 + 179.439842541 $i$ | 848 | 75 342.919675255 + 180.100900268 $i$ | 888 | 78 896.903895802 + 180.625787844 $i$ | 928 | 82 450.911263404 + 181.320096121 $i$ |
| 809 | 71 877.937112560 + 179.096361836 $i$ | 849 | 75 432.069841603 + 179.796799118 $i$ | 889 | 78 986.106407954 + 180.650336538 $i$ | 929 | 82 540.060620844 + 181.378425833 $i$ |
| 810 | 71 967.208554606 + 179.261176721 $i$ | 850 | 75 521.064216896 + 180.217353676 $i$ | 890 | 79 074.837357585 + 180.944717968 $i$ | 930 | 82 628.692449447 + 181.670024561 $i$ |
| 811 | 72 055.667302400 + 179.546868482 $i$ | 851 | 75 609.489348422 + 180.131781306 $i$ | 891 | 79 163.474423468 + 180.772322428 $i$ | 931 | 82 717.391354650 + 181.362854414 $i$ |
| 812 | 72 144.502004657 + 179.120560403 $i$ | 852 | 75 698.626177612 + 179.894115647 $i$ | 892 | 79 252.571701356 + 180.653554395 $i$ | 932 | 82 806.648236305 + 181.387563600 $i$ |
| 813 | 72 233.758195853 + 179.356233812 $i$ | 853 | 75 787.536203838 + 180.212020937 $i$ | 893 | 79 341.481865064 + 180.997238999 $i$ | 933 | 82 895.218821578 + 181.792887786 $i$ |
| 814 | 72 322.184075658 + 179.562944081 $i$ | 854 | 75 876.151918767 + 180.237127504 $i$ | 894 | 79 429.947266066 + 180.885395569 $i$ | 934 | 82 983.950410201 + 181.327650375 $i$ |
| 815 | 72 411.126251279 + 179.228117006 $i$ | 855 | 75 965.036931053 + 179.958444584 $i$ | 895 | 79 519.163043540 + 180.614218387 $i$ | 935 | 83 073.176629779 + 181.536564001 $i$ |
| 816 | 72 500.168903638 + 179.371071982 $i$ | 856 | 76 054.207443109 + 180.200741501 $i$ | 896 | 79 608.002495764 + 181.136839060 $i$ | 936 | 83 161.768479926 + 181.722226135 $i$ |
| 817 | 72 588.887617992 + 179.623096367 $i$ | 857 | 76 142.611596567 + 180.371160890 $i$ | 897 | 79 696.504984208 + 180.846796963 $i$ | 937 | 83 250.524749836 + 181.506233369 $i$ |
| 818 | 72 677.521175629 + 179.333775367 $i$ | 858 | 76 231.617553242 + 179.939743288 $i$ | 898 | 79 785.703302014 + 180.767168653 $i$ | 938 | 83 339.652951079 + 181.480443580 $i$ |
| 819 | 72 766.803440988 + 179.356075107 $i$ | 859 | 76 320.735463084 + 180.316291115 $i$ | 899 | 79 874.522692548 + 181.077842371 $i$ | 939 | 83 428.393371299 + 181.815524473 $i$ |
| 820 | 72 855.399064832 + 179.748820300 $i$ | 860 | 76 409.155702786 + 180.375444561 $i$ | 900 | 79 963.104381532 + 180.990404121 $i$ | 940 | 83 517.005103002 + 181.580959545 $i$ |
| 821 | 72 944.072536065 + 179.350955630 $i$ | 861 | 76 498.194534808 + 180.049505926 $i$ | 901 | 80 052.172169623 + 180.790609585 $i$ | 941 | 83 606.218802268 + 181.464809141 $i$ |
| 822 | 73 033.369324520 + 179.429243801 $i$ | 862 | 76 587.209978465 + 180.326634106 $i$ | 902 | 80 141.132641772 + 181.091922387 $i$ | 942 | 83 694.946410520 + 181.948018837 $i$ |
| 823 | 73 121.895381857 + 179.796334768 $i$ | 863 | 76 675.798348929 + 180.428129467 $i$ | 903 | 80 229.632716662 + 181.124322414 $i$ | 943 | 83 783.533326587 + 181.557444472 $i$ |
| 824 | 73 210.705123130 + 179.407857426 $i$ | 864 | 76 764.642103892 + 180.167507160 $i$ | 904 | 80 318.679437765 + 180.738551977 $i$ | 944 | 83 872.789871758 + 181.586263048 $i$ |
| 825 | 73 299.788125761 + 179.516374690 $i$ | 865 | 76 853.809775617 + 180.278933543 $i$ | 905 | 80 407.727359953 + 181.251832226 $i$ | 945 | 83 961.441567948 + 181.918334697 $i$ |
| 826 | 73 388.571918769 + 179.774547861 $i$ | 866 | 76 942.350966992 + 180.591933235 $i$ | 906 | 80 496.130308715 + 181.082105940 $i$ | 946 | 84 050.157944640 + 181.679034910 $i$ |
| 827 | 73 477.150033164 + 179.579098441 $i$ | 867 | 77 031.137256020 + 180.121941743 $i$ | 907 | 80 585.279040085 + 180.888393574 $i$ | 947 | 84 139.226534672 + 181.614066744 $i$ |
| 828 | 73 566.367942806 + 179.463541130 $i$ | 868 | 77 120.407174035 + 180.418324349 $i$ | 908 | 80 674.206881310 + 181.211673071 $i$ | 948 | 84 228.104571693 + 181.929362449 $i$ |
| 829 | 73 655.136203573 + 179.912868982 $i$ | 869 | 77 208.847784096 + 180.572824005 $i$ | 909 | 80 762.757615727 + 181.183489703 $i$ | 949 | 84 316.618189777 + 181.799090510 $i$ |
| 830 | 73 743.655525245 + 179.596541560 $i$ | 870 | 77 297.737314534 + 180.239057195 $i$ | 910 | 80 851.748395510 + 180.950685392 $i$ | 950 | 84 405.794691409 + 181.587576996 $i$ |
| 831 | 73 832.961897323 + 179.528030244 $i$ | 871 | 77 386.879985249 + 180.437621309 $i$ | 911 | 80 940.786967048 + 181.185006207 $i$ | 951 | 84 494.663036507 + 182.037641978 $i$ |
| 832 | 73 921.615585889 + 179.964748903 $i$ | 872 | 77 475.457683767 + 180.612834380 $i$ | 912 | 81 029.309377623 + 181.332681740 $i$ | 952 | 84 583.130820321 + 181.820215605 $i$ |
| 833 | 74 010.268716664 + 179.648673606 $i$ | 873 | 77 564.261711046 + 180.369246770 $i$ | 913 | 81 118.229999415 + 180.932668512 $i$ | 953 | 84 672.389668081 + 181.657203166 $i$ |
| 834 | 74 099.438260680 + 179.642036587 $i$ | 874 | 77 653.384806308 + 180.374217359 $i$ | 914 | 81 207.420253695 + 181.293136524 $i$ | 954 | 84 761.143709326 + 182.068061439 $i$ |
| 835 | 74 188.214030066 + 179.914109287 $i$ | 875 | 77 742.084609781 + 180.777773720 $i$ | 915 | 81 295.778700076 + 181.344754196 $i$ | 955 | 84 849.765914295 + 181.851893874 $i$ |
| 836 | 74 276.815689243 + 179.829105692 $i$ | 876 | 77 830.704139611 + 180.352868869 $i$ | 916 | 81 384.858500367 + 181.007100305 $i$ | 956 | 84 938.816503295 + 181.770750405 $i$ |
| 837 | 74 365.911872612 + 179.587996967 $i$ | 877 | 77 920.032511915 + 180.472769887 $i$ | 917 | 81 473.863874930 + 181.325842935 $i$ | 957 | 85 027.787738814 + 182.014180927 $i$ |
| 838 | 74 454.862756472 + 180.049410042 $i$ | 878 | 78 008.542809239 + 180.790641239 $i$ | 918 | 81 562.417585839 + 181.377928057 $i$ | 958 | 85 116.262736958 + 182.021551056 $i$ |
| 839 | 74 543.274365025 + 179.844385262 $i$ | 879 | 78 097.327147212 + 180.434444683 $i$ | 919 | 81 651.330127399 + 181.133158237 $i$ | 959 | 85 205.362435047 + 181.719101431 $i$ |
| 840 | 74 632.521638747 + 179.649013579 $i$ | 880 | 78 186.502341432 + 180.523808493 $i$ | 920 | 81 740.431446215 + 181.277836792 $i$ | 960 | 85 294.349991776 + 182.112530196 $i$ |

| $k$ | $\sigma_k$ | $k$ | $\sigma_k$ | $k$ | $\sigma_k$ |
|---|---|---|---|---|---|
| 961 | 85 382.780671823 + 182.072382871 $i$ | 1001 | 88 936.787398498 + 182.524841140 $i$ | 1041 | 92 490.771119378 + 183.119758781 $i$ |
| 962 | 85 471.942364099 + 181.744325225 $i$ | 1002 | 89 025.920690750 + 182.535184242 $i$ | 1042 | 92 579.881221161 + 183.161867142 $i$ |
| 963 | 85 560.843849240 + 182.203771885 $i$ | 1003 | 89 114.693393400 + 182.799076126 $i$ | 1043 | 92 668.549316488 + 183.508508993 $i$ |
| 964 | 85 649.397964330 + 182.047641269 $i$ | 1004 | 89 203.318692937 + 182.678664177 $i$ | 1044 | 92 757.222668025 + 183.106067377 $i$ |
| 965 | 85 738.412613746 + 181.903891084 $i$ | 1005 | 89 292.414156322 + 182.479220127 $i$ | 1045 | 92 846.522601730 + 183.244174688 $i$ |
| 966 | 85 827.425457750 + 182.118848573 $i$ | 1006 | 89 381.323017681 + 182.921195670 $i$ | 1046 | 92 935.001489332 + 183.515061868 $i$ |
| 967 | 85 915.959919163 + 182.222780368 $i$ | 1007 | 89 469.778074968 + 182.687790043 $i$ | 1047 | 93 023.851898501 + 183.184347086 $i$ |
| 968 | 86 004.904610335 + 181.871573124 $i$ | 1008 | 89 559.033074186 + 182.543429600 $i$ | 1048 | 93 112.976303414 + 183.288032801 $i$ |
| 969 | 86 094.031833272 + 182.188391162 $i$ | 1009 | 89 647.801129223 + 182.955955448 $i$ | 1049 | 93 201.630339025 + 183.526893956 $i$ |
| 970 | 86 182.454477565 + 182.285382145 $i$ | 1010 | 89 736.386182903 + 182.730369640 $i$ | 1050 | 93 290.366299600 + 183.283017519 $i$ |
| 971 | 86 271.481763261 + 181.889029126 $i$ | 1011 | 89 825.526743452 + 182.627443157 $i$ | 1051 | 93 379.470814176 + 183.271497928 $i$ |
| 972 | 86 360.543900663 + 182.282981449 $i$ | 1012 | 89 914.356272103 + 182.930731381 $i$ | 1052 | 93 468.272652485 + 183.618566158 $i$ |
| 973 | 86 449.025649266 + 182.265177932 $i$ | 1013 | 90 002.978848949 + 182.855804618 $i$ | 1053 | 93 556.805889044 + 183.321635590 $i$ |
| 974 | 86 538.024962675 + 182.040859803 $i$ | 1014 | 90 091.967915948 + 182.627768266 $i$ | 1054 | 93 646.113627135 + 183.294985087 $i$ |
| 975 | 86 627.045584686 + 182.219159654 $i$ | 1015 | 90 181.033272613 + 182.978661373 $i$ | 1055 | 93 734.731754387 + 183.683639247 $i$ |
| 976 | 86 715.658802839 + 182.397591165 $i$ | 1016 | 90 269.399484692 + 182.925111024 $i$ | 1056 | 93 823.422973156 + 183.342873686 $i$ |
| 977 | 86 804.464241653 + 182.065029319 $i$ | 1017 | 90 358.608749669 + 182.650135519 $i$ | 1057 | 93 912.596579135 + 183.388817911 $i$ |
| 978 | 86 893.689529343 + 182.236381185 $i$ | 1018 | 90 447.505419698 + 183.042993565 $i$ | 1058 | 94 001.314001374 + 183.647382481 $i$ |
| 979 | 86 982.131202095 + 182.488459714 $i$ | 1019 | 90 536.002930192 + 182.945090450 $i$ | 1059 | 94 089.980977742 + 183.477581592 $i$ |
| 980 | 87 071.049802708 + 182.075598762 $i$ | 1020 | 90 625.122443429 + 182.733548650 $i$ | 1060 | 94 179.077010208 + 183.361584634 $i$ |
| 981 | 87 160.211598663 + 182.320016324 $i$ | 1021 | 90 714.018346558 + 183.038172487 $i$ | 1061 | 94 267.940807875 + 183.710523637 $i$ |
| 982 | 87 248.682995624 + 182.488389042 $i$ | 1022 | 90 802.632783209 + 183.035369019 $i$ | 1062 | 94 356.453800313 + 183.556642336 $i$ |
| 983 | 87 337.624274649 + 182.168997438 $i$ | 1023 | 90 891.547900302 + 182.786533582 $i$ | 1063 | 94 445.664832413 + 183.342889468 $i$ |
| 984 | 87 426.657387409 + 182.337360748 $i$ | 1024 | 90 980.681390275 + 183.034061872 $i$ | 1064 | 94 534.460640698 + 183.825333126 $i$ |
| 985 | 87 515.365979267 + 182.538173321 $i$ | 1025 | 91 069.096564204 + 183.157299770 $i$ | 1065 | 94 623.013838594 + 183.528610634 $i$ |
| 986 | 87 604.033412815 + 182.265365091 $i$ | 1026 | 91 158.142221839 + 182.756284235 $i$ | 1066 | 94 712.207500144 + 183.480027669 $i$ |
| 987 | 87 693.299410188 + 182.309947852 $i$ | 1027 | 91 247.196001554 + 183.141653174 $i$ | 1067 | 94 800.986342211 + 183.756120486 $i$ |
| 988 | 87 781.858781047 + 182.664920649 $i$ | 1028 | 91 335.658873507 + 183.141188319 $i$ | 1068 | 94 889.610435145 + 183.670362250 $i$ |
| 989 | 87 870.608082926 + 182.259585113 $i$ | 1029 | 91 424.692863696 + 182.863894520 $i$ | 1069 | 94 978.666853683 + 183.474404783 $i$ |
| 990 | 87 959.837003527 + 182.385062735 $i$ | 1030 | 91 513.688858961 + 183.134144559 $i$ | 1070 | 95 067.610356584 + 183.797369566 $i$ |
| 991 | 88 048.385645017 + 182.683711790 $i$ | 1031 | 91 602.290014389 + 183.204912477 $i$ | 1071 | 95 156.120655344 + 183.761141022 $i$ |
| 992 | 88 137.201787314 + 182.324670990 $i$ | 1032 | 91 691.146806080 + 182.948944650 $i$ | 1072 | 95 245.196883990 + 183.452081469 $i$ |
| 993 | 88 226.290474530 + 182.450547821 $i$ | 1033 | 91 780.293344333 + 183.094078418 $i$ | 1073 | 95 334.190067235 + 183.910189241 $i$ |
| 994 | 88 315.041440863 + 182.660606181 $i$ | 1034 | 91 868.822040244 + 183.347982218 $i$ | 1074 | 95 422.612688172 + 183.729732325 $i$ |
| 995 | 88 403.652471155 + 182.483265846 $i$ | 1035 | 91 957.665222641 + 182.912242506 $i$ | 1075 | 95 511.800188070 + 183.584095273 $i$ |
| 996 | 88 492.878390862 + 182.379764211 $i$ | 1036 | 92 046.883628599 + 183.206960794 $i$ | 1076 | 95 600.667793773 + 183.860075638 $i$ |
| 997 | 88 581.585373676 + 182.805682780 $i$ | 1037 | 92 135.316278041 + 183.324468360 $i$ | 1077 | 95 689.250586623 + 183.840810518 $i$ |
| 998 | 88 670.181020611 + 182.474125633 $i$ | 1038 | 92 224.272445087 + 183.023237867 $i$ | 1078 | 95 778.250318194 + 183.615079144 $i$ |
| 999 | 88 759.451490321 + 182.452493053 $i$ | 1039 | 92 313.340199995 + 183.199485943 $i$ | 1079 | 95 867.271191855 + 183.862674080 $i$ |
| 1000 | 88 848.085816443 + 182.834334323 $i$ | 1040 | 92 401.940640520 + 183.385477526 $i$ | 1080 | 95 955.790367844 + 183.957430862 $i$ |
| 1081 | 96 044.757203375 + 183.593863749 $i$ | | | | |
| 1082 | 96 133.881369935 + 183.948026848 $i$ | | | | |
| 1083 | 96 222.265889907 + 183.964208152 $i$ | | | | |
| 1084 | 96 311.374251515 + 183.665901788 $i$ | | | | |
| 1085 | 96 400.330675483 + 183.972320628 $i$ | | | | |
| 1086 | 96 488.916151558 + 183.997712827 $i$ | | | | |
| 1087 | 96 577.829623472 + 183.763010243 $i$ | | | | |
| 1088 | 96 666.909717793 + 183.929792469 $i$ | | | | |
| 1089 | 96 755.479533168 + 184.135605040 $i$ | | | | |
| 1090 | 96 844.328566781 + 183.750710046 $i$ | | | | |
| 1091 | 96 933.515461589 + 183.985317734 $i$ | | | | |
| 1092 | 97 021.970653528 + 184.179791628 $i$ | | | | |
| 1093 | 97 110.914150514 + 183.774803769 $i$ | | | | |
| 1094 | 97 200.004847240 + 184.078426838 $i$ | | | | |
| 1095 | 97 288.575509515 + 184.139922981 $i$ | | | | |
| 1096 | 97 377.425219521 + 183.939869868 $i$ | | | | |
| 1097 | 97 466.545102291 + 183.987586391 $i$ | | | | |
| 1098 | 97 555.161568015 + 184.278134889 $i$ | | | | |
| 1099 | 97 643.911491749 + 183.942986780 $i$ | | | | |
| 1100 | 97 733.131826526 + 184.023767389 $i$ | | | | |
| 1101 | 97 821.688834573 + 184.356000060 $i$ | | | | |
| 1102 | 97 910.464663280 + 183.930507305 $i$ | | | | |
| 1103 | 97 999.661888473 + 184.142040133 $i$ | | | | |
| 1104 | 98 088.236450754 + 184.299135997 $i$ | | | | |
| 1105 | 98 177.044071272 + 184.094752097 $i$ | | | | |
| 1106 | 98 266.134816320 + 184.068422802 $i$ | | | | |
| 1107 | 98 354.867079294 + 184.406709378 $i$ | | | | |
| 1108 | 98 443.510609511 + 184.133280879 $i$ | | | | |
| 1109 | 98 532.719469932 + 184.081695987 $i$ | | | | |
| 1110 | 98 621.412408667 + 184.495503047 $i$ | | | | |
| 1111 | 98 710.038658428 + 184.120885394 $i$ | | | | |
| 1112 | 98 799.285565383 + 184.191341123 $i$ | | | | |
| 1113 | 98 887.923732930 + 184.456329716 $i$ | | | | |
| 1114 | 98 976.649264820 + 184.234545358 $i$ | | | | |
| 1115 | 99 065.726242597 + 184.196773820 $i$ | | | | |
| 1116 | 99 154.579426294 + 184.481699483 $i$ | | | | |
| 1117 | 99 243.109616005 + 184.338287649 $i$ | | | | |
| 1118 | 99 332.309196181 + 184.163310547 $i$ | | | | |
| 1119 | 99 421.118817893 + 184.586417599 $i$ | | | | |
| 1120 | 99 509.642739768 + 184.340017937 $i$ | | | | |

| $k$ | $\sigma_k$ | $k$ | $\sigma_k$ | $k$ | $\sigma_k$ |
|---|---|---|---|---|---|
| 1121 | 99 598.875603009 + 184.237399826 $i$ | 1161 | 103 152.795515809 + 184.907820661 $i$ | 1201 | 106 706.780291224 + 185.473094973 $i$ |
| 1122 | 99 687.619998789 + 184.605313656 $i$ | 1162 | 103 241.500260133 + 185.094415064 $i$ | 1202 | 106 795.286527444 + 185.681173759 $i$ |
| 1123 | 99 776.268500618 + 184.385955534 $i$ | 1163 | 103 330.160011979 + 184.928551827 $i$ | 1203 | 106 884.188524006 + 185.282017848 $i$ |
| 1124 | 99 865.325847888 + 184.313995665 $i$ | 1164 | 103 419.378680348 + 184.837455150 $i$ | 1204 | 106 973.353373926 + 185.575118739 $i$ |
| 1125 | 99 954.246654214 + 184.549634790 $i$ | 1165 | 103 508.045573843 + 185.244190750 $i$ | 1205 | 107 061.795302531 + 185.658524199 $i$ |
| 1126 | 100 042.760074680 + 184.548629616 $i$ | 1166 | 103 596.698794200 + 184.888746939 $i$ | 1206 | 107 150.789872451 + 185.380669940 $i$ |
| 1127 | 100 131.873400786 + 184.248619431 $i$ | 1167 | 103 685.933990004 + 184.925256621 $i$ | 1207 | 107 239.809452292 + 185.562023200 $i$ |
| 1128 | 100 220.804947997 + 184.669562947 $i$ | 1168 | 103 774.564123504 + 185.247892783 $i$ | 1208 | 107 328.441147602 + 185.720615647 $i$ |
| 1129 | 100 309.290293158 + 184.558918806 $i$ | 1169 | 103 863.290925637 + 184.947822005 $i$ | 1209 | 107 417.260853544 + 185.444290965 $i$ |
| 1130 | 100 398.438665860 + 184.303590781 $i$ | 1170 | 103 952.406156520 + 184.985399723 $i$ | 1210 | 107 506.376942925 + 185.538764841 $i$ |
| 1131 | 100 487.323216346 + 184.718392912 $i$ | 1171 | 104 041.178249977 + 185.225244677 $i$ | 1211 | 107 595.016621470 + 185.818789463 $i$ |
| 1132 | 100 575.884325678 + 184.551585615 $i$ | 1172 | 104 129.809411327 + 185.082165342 $i$ | 1212 | 107 683.731551463 + 185.442309655 $i$ |
| 1133 | 100 664.925524036 + 184.438300271 $i$ | 1173 | 104 218.922910753 + 184.927760209 $i$ | 1213 | 107 773.002961709 + 185.609494665 $i$ |
| 1134 | 100 753.899851873 + 184.627971894 $i$ | 1174 | 104 307.788238695 + 185.333262362 $i$ | 1214 | 107 861.479855668 + 185.832463918 $i$ |
| 1135 | 100 842.436990566 + 184.725517337 $i$ | 1175 | 104 396.274694782 + 185.084124553 $i$ | 1215 | 107 950.368633448 + 185.506976006 $i$ |
| 1136 | 100 931.424405520 + 184.387277312 $i$ | 1176 | 104 485.538355031 + 184.991251494 $i$ | 1216 | 108 039.444918902 + 185.641602319 $i$ |
| 1137 | 101 020.506109495 + 184.713714399 $i$ | 1177 | 104 574.262860821 + 185.358329879 $i$ | 1217 | 108 128.126219657 + 185.843123670 $i$ |
| 1138 | 101 108.933774269 + 184.751172166 $i$ | 1178 | 104 662.895280632 + 185.128459714 $i$ | 1218 | 108 216.855121040 + 185.597216412 $i$ |
| 1139 | 101 197.998800307 + 184.421335671 $i$ | 1179 | 104 752.014478524 + 185.057968998 $i$ | 1219 | 108 305.979979035 + 185.615691556 $i$ |
| 1140 | 101 287.015774991 + 184.780637810 $i$ | 1180 | 104 840.843773894 + 185.334254313 $i$ | 1220 | 108 394.726356873 + 185.923262753 $i$ |
| 1141 | 101 375.512782407 + 184.743419109 $i$ | 1181 | 104 929.454797686 + 185.244206958 $i$ | 1221 | 108 483.318272069 + 185.636338853 $i$ |
| 1142 | 101 464.530885302 + 184.548178002 $i$ | 1182 | 105 018.489305466 + 185.051041375 $i$ | 1222 | 108 572.610883659 + 185.632268688 $i$ |
| 1143 | 101 553.530106243 + 184.718339029 $i$ | 1183 | 105 107.492421011 + 185.378334880 $i$ | 1223 | 108 661.189358055 + 185.984858252 $i$ |
| 1144 | 101 642.136397020 + 184.871531530 $i$ | 1184 | 105 195.896253690 + 185.309456489 $i$ | 1224 | 108 749.944306630 + 185.651902396 $i$ |
| 1145 | 101 730.976852697 + 184.546491786 $i$ | 1185 | 105 285.123254579 + 185.054002832 $i$ | 1225 | 108 839.081606027 + 185.710551869 $i$ |
| 1146 | 101 820.163531111 + 184.747463926 $i$ | 1186 | 105 373.952664227 + 185.450921425 $i$ | 1226 | 108 927.784941086 + 185.945861748 $i$ |
| 1147 | 101 908.610698568 + 184.948367054 $i$ | 1187 | 105 462.522478140 + 185.315547148 $i$ | 1227 | 109 016.476167856 + 185.782246538 $i$ |
| 1148 | 101 997.567280365 + 184.551194676 $i$ | 1188 | 105 551.608907395 + 185.138187631 $i$ | 1228 | 109 105.581828323 + 185.680318053 $i$ |
| 1149 | 102 086.667881689 + 184.832999916 $i$ | 1189 | 105 640.506828265 + 185.431515302 $i$ | 1229 | 109 194.410734325 + 186.015428819 $i$ |
| 1150 | 102 175.190853288 + 184.940770598 $i$ | 1190 | 105 729.114436566 + 185.400177663 $i$ | 1230 | 109 282.947524823 + 185.827100801 $i$ |
| 1151 | 102 264.115505079 + 184.643510284 $i$ | 1191 | 105 818.061963322 + 185.175669831 $i$ | 1231 | 109 372.167031798 + 185.683454652 $i$ |
| 1152 | 102 353.150361205 + 184.829218376 $i$ | 1192 | 105 907.152218991 + 185.425648633 $i$ | 1232 | 109 460.933120185 + 186.106678006 $i$ |
| 1153 | 102 441.837413074 + 184.991241596 $i$ | 1193 | 105 995.578417336 + 185.510044153 $i$ | 1233 | 109 549.513867098 + 185.804516049 $i$ |
| 1154 | 102 530.547662228 + 184.731875006 $i$ | 1194 | 106 084.658043277 + 185.143843984 $i$ | 1234 | 109 638.699915906 + 185.798586379 $i$ |
| 1155 | 102 619.791475054 + 184.788277975 $i$ | 1195 | 106 173.659611461 + 185.530887439 $i$ | 1235 | 109 727.462967293 + 186.043022955 $i$ |
| 1156 | 102 708.320319990 + 185.113290088 $i$ | 1196 | 106 262.151763900 + 185.482998414 $i$ | 1236 | 109 816.110147574 + 185.942637726 $i$ |
| 1157 | 102 797.132101079 + 184.706215490 $i$ | 1197 | 106 351.193310279 + 185.257804433 $i$ | 1237 | 109 905.158581471 + 185.772388620 $i$ |
| 1158 | 102 886.308070720 + 184.874220098 $i$ | 1198 | 106 440.173921058 + 185.503555874 $i$ | 1238 | 109 994.097080620 + 186.095495211 $i$ |
| 1159 | 102 974.868788996 + 185.109554981 $i$ | 1199 | 106 528.768332501 + 185.557444454 $i$ | 1239 | 110 082.598198214 + 186.003220487 $i$ |
| 1160 | 103 063.701157794 + 184.790373726 $i$ | 1200 | 106 617.655584035 + 185.314794478 $i$ | 1240 | 110 171.710534512 + 185.769656141 $i$ |

| $k$ | $\sigma_k$ |
|---|---|
| 1241 | 110 260.661098345 + 186.180408303 $i$ |
| 1242 | 110 349.110959584 + 185.996422156 $i$ |
| 1243 | 110 438.311745737 + 185.867945790 $i$ |
| 1244 | 110 527.123479619 + 186.133330396 $i$ |
| 1245 | 110 615.756317687 + 186.106003133 $i$ |
| 1246 | 110 704.747401981 + 185.880879805 $i$ |
| 1247 | 110 793.752974270 + 186.148164034 $i$ |
| 1248 | 110 882.269076741 + 186.192703007 $i$ |
| 1249 | 110 971.274540642 + 185.867245639 $i$ |
| 1250 | 111 060.338840657 + 186.221881222 $i$ |
| 1251 | 111 148.760634677 + 186.203460568 $i$ |
| 1252 | 111 237.883373467 + 185.935572350 $i$ |
| 1253 | 111 326.799784075 + 186.234970503 $i$ |
| 1254 | 111 415.408474643 + 186.235478830 $i$ |
| 1255 | 111 504.329364889 + 186.020575187 $i$ |
| 1256 | 111 593.402605096 + 186.192745922 $i$ |
| 1257 | 111 681.950231877 + 186.355563348 $i$ |
| 1258 | 111 770.836417793 + 185.998308459 $i$ |
| 1259 | 111 859.996135072 + 186.259413185 $i$ |
| 1260 | 111 948.450121367 + 186.381358747 $i$ |
| 1261 | 112 037.423144800 + 186.035981765 $i$ |
| 1262 | 112 126.488108165 + 186.317894425 $i$ |
| 1263 | 112 215.049590687 + 186.360245037 $i$ |
| 1264 | 112 303.933019845 + 186.175734682 $i$ |
| 1265 | 112 393.023798752 + 186.233905134 $i$ |
| 1266 | 112 481.644384515 + 186.506270434 $i$ |
| 1267 | 112 570.418508950 + 186.145570159 $i$ |
| 1268 | 112 659.613195209 + 186.291843693 $i$ |
| 1269 | 112 748.164688005 + 186.548800745 $i$ |
| 1270 | 112 836.977809091 + 186.158696057 $i$ |
| 1271 | 112 926.141020705 + 186.373905400 $i$ |
| 1272 | 113 014.714079138 + 186.508860561 $i$ |
| 1273 | 113 103.551861021 + 186.302368794 $i$ |
| 1274 | 113 192.616080121 + 186.301472600 $i$ |
| 1275 | 113 281.353444923 + 186.616635762 $i$ |
| 1276 | 113 370.002734255 + 186.315182301 $i$ |
| 1277 | 113 459.219284764 + 186.332373549 $i$ |
| 1278 | 113 547.877657949 + 186.674487918 $i$ |
| 1279 | 113 636.546363070 + 186.322984897 $i$ |
| 1280 | 113 725.776234134 + 186.416646870 $i$ |

| $k$ | $\sigma_k$ | $k$ | $\sigma_k$ | $k$ | $\sigma_k$ |
|---|---|---|---|---|---|
| 1281 | 113 814.396046532 + 186.639700637 $i$ | 1321 | 117 368.300149048 + 187.102321943 $i$ | 1361 | 120 922.067179640 + 187.546219638 $i$ |
| 1282 | 113 903.147658226 + 186.433522619 $i$ | 1322 | 117 457.059543112 + 186.858949894 $i$ | 1362 | 121 011.167237876 + 187.214623828 $i$ |
| 1283 | 113 992.225073486 + 186.406451557 $i$ | 1323 | 117 546.277499523 + 186.931248915 $i$ | 1363 | 121 100.126463252 + 187.603707274 $i$ |
| 1284 | 114 081.049657838 + 186.677095213 $i$ | 1324 | 117 634.789666643 + 187.221139270 $i$ | 1364 | 121 188.647721856 + 187.511454336 $i$ |
| 1285 | 114 169.608092055 + 186.514463138 $i$ | 1325 | 117 723.644396587 + 186.818282262 $i$ | 1365 | 121 277.692464527 + 187.328363951 $i$ |
| 1286 | 114 258.813103837 + 186.372649833 $i$ | 1326 | 117 812.785667417 + 187.029588404 $i$ | 1366 | 121 366.651434209 + 187.557062728 $i$ |
| 1287 | 114 347.575852633 + 186.779566110 $i$ | 1327 | 117 901.356035369 + 187.193328981 $i$ | 1367 | 121 455.248907858 + 187.602232693 $i$ |
| 1288 | 114 436.157992635 + 186.504940339 $i$ | 1328 | 117 990.191273646 + 186.914282648 $i$ | 1368 | 121 544.162721259 + 187.360428455 $i$ |
| 1289 | 114 525.366046960 + 186.448767510 $i$ | 1329 | 118 079.285724223 + 187.037131011 $i$ | 1369 | 121 633.256733001 + 187.544795953 $i$ |
| 1290 | 114 614.091079864 + 186.773068944 $i$ | 1330 | 118 167.973611561 + 187.210899299 $i$ | 1370 | 121 721.765213875 + 187.709521329 $i$ |
| 1291 | 114 702.756926565 + 186.566341300 $i$ | 1331 | 118 256.672109670 + 187.021412419 $i$ | 1371 | 121 810.710036625 + 187.326545166 $i$ |
| 1292 | 114 791.827973835 + 186.510247724 $i$ | 1332 | 118 345.860766297 + 186.969988479 $i$ | 1372 | 121 899.816340389 + 187.634084067 $i$ |
| 1293 | 114 880.718723107 + 186.738195835 $i$ | 1333 | 118 434.521705130 + 187.353530243 $i$ | 1373 | 121 988.284496671 + 187.684582405 $i$ |
| 1294 | 114 969.257860022 + 186.707145850 $i$ | 1334 | 118 523.209940776 + 186.964955193 $i$ | 1374 | 122 077.295736310 + 187.422322768 $i$ |
| 1295 | 115 058.376833174 + 186.436924431 $i$ | 1335 | 118 612.412363877 + 187.068930668 $i$ | 1375 | 122 166.283083682 + 187.619639705 $i$ |
| 1296 | 115 147.280894618 + 186.859437049 $i$ | 1336 | 118 701.046722260 + 187.330431598 $i$ | 1376 | 122 254.935912701 + 187.739735359 $i$ |
| 1297 | 115 235.783924500 + 186.682812813 $i$ | 1337 | 118 789.792329670 + 187.042309131 $i$ | 1377 | 122 343.744459163 + 187.472813964 $i$ |
| 1298 | 115 324.934463241 + 186.516034947 $i$ | 1338 | 118 878.893188606 + 187.098101555 $i$ | 1378 | 122 432.877033432 + 187.604882322 $i$ |
| 1299 | 115 413.805113148 + 186.870173883 $i$ | 1339 | 118 967.653042659 + 187.321383476 $i$ | 1379 | 122 521.481538633 + 187.827293158 $i$ |
| 1300 | 115 502.367054615 + 186.705920365 $i$ | 1340 | 119 056.300608937 + 187.167682464 $i$ | 1380 | 122 610.250907592 + 187.473985168 $i$ |
| 1301 | 115 591.431478189 + 186.614008990 $i$ | 1341 | 119 145.433830536 + 187.043103676 $i$ | 1381 | 122 699.481117398 + 187.654903389 $i$ |
| 1302 | 115 680.372587458 + 186.797532070 $i$ | 1342 | 119 234.246914220 + 187.419206040 $i$ | 1382 | 122 787.955379287 + 187.841054457 $i$ |
| 1303 | 115 768.926901874 + 186.870423109 $i$ | 1343 | 119 322.778170776 + 187.161854420 $i$ | 1383 | 122 876.874986009 + 187.529810127 $i$ |
| 1304 | 115 857.929444436 + 186.544121023 $i$ | 1344 | 119 412.041348923 + 187.104081779 $i$ | 1384 | 122 965.922558761 + 187.689097693 $i$ |
| 1305 | 115 946.972710073 + 186.897663391 $i$ | 1345 | 119 500.724537261 + 187.434680037 $i$ | 1385 | 123 054.613481588 + 187.848307324 $i$ |
| 1306 | 116 035.423081502 + 186.876523794 $i$ | 1346 | 119 589.401204274 + 187.205156325 $i$ | 1386 | 123 143.340845057 + 187.617522717 $i$ |
| 1307 | 116 124.512129629 + 186.595273584 $i$ | 1347 | 119 678.505993598 + 187.160743261 $i$ | 1387 | 123 232.487717944 + 187.654530034 $i$ |
| 1308 | 116 213.478812215 + 186.935117755 $i$ | 1348 | 119 767.328714146 + 187.409411281 $i$ | 1388 | 123 321.182143194 + 187.933079704 $i$ |
| 1309 | 116 302.009529338 + 186.884043228 $i$ | 1349 | 119 855.932716289 + 187.310007618 $i$ | 1389 | 123 409.837042456 + 187.643361419 $i$ |
| 1310 | 116 391.033227844 + 186.696968513 $i$ | 1350 | 119 945.002633312 + 187.143280854 $i$ | 1390 | 123 499.095671254 + 187.667762639 $i$ |
| 1311 | 116 480.006113653 + 186.878343499 $i$ | 1351 | 120 033.948704210 + 187.464759110 $i$ | 1391 | 123 587.660312222 + 187.987239079 $i$ |
| 1312 | 116 568.620361044 + 187.004826511 $i$ | 1352 | 120 122.400795095 + 187.365539251 $i$ | 1392 | 123 676.453173434 + 187.649593413 $i$ |
| 1313 | 116 657.486216377 + 186.681118803 $i$ | 1353 | 120 211.620707533 + 187.141361499 $i$ | 1393 | 123 765.556885241 + 187.746975013 $i$ |
| 1314 | 116 746.635609754 + 186.916329400 $i$ | 1354 | 120 300.418074836 + 187.540002719 $i$ | 1394 | 123 854.275446232 + 187.949968443 $i$ |
| 1315 | 116 835.092802770 + 187.064470801 $i$ | 1355 | 120 389.030644905 + 187.354561079 $i$ | 1395 | 123 942.970115447 + 187.771350518 $i$ |
| 1316 | 116 924.080991665 + 186.692466962 $i$ | 1356 | 120 478.089569622 + 187.234282645 $i$ | 1396 | 124 032.075054589 + 187.702759132 $i$ |
| 1317 | 117 013.138482725 + 186.991329437 $i$ | 1357 | 120 567.000913448 + 187.501791043 $i$ | 1397 | 124 120.884876999 + 188.024936401 $i$ |
| 1318 | 117 101.681262849 + 187.044445115 $i$ | 1358 | 120 655.589930994 + 187.449202860 $i$ | 1398 | 124 209.445441213 + 187.798196022 $i$ |
| 1319 | 117 190.606988930 + 186.793558369 $i$ | 1359 | 120 744.574287307 + 187.249805964 $i$ | 1399 | 124 298.662949183 + 187.717383565 $i$ |
| 1320 | 117 279.653428861 + 186.967657551 $i$ | 1360 | 120 833.622320260 + 187.503747812 $i$ | 1400 | 124 387.398562732 + 188.094613175 $i$ |
| | | | | 1401 | 124 476.019651738 + 187.792286236 $i$ |
| | | | | 1402 | 124 565.190257606 + 187.810648483 $i$ |
| | | | | 1403 | 124 653.934616051 + 188.037718118 $i$ |
| | | | | 1404 | 124 742.609186133 + 187.925865482 $i$ |
| | | | | 1405 | 124 831.655790617 + 187.773411530 $i$ |
| | | | | 1406 | 124 920.572675169 + 188.093900940 $i$ |
| | | | | 1407 | 125 009.081493911 + 187.970601770 $i$ |
| | | | | 1408 | 125 098.235009971 + 187.778695603 $i$ |
| | | | | 1409 | 125 187.109534379 + 188.151573677 $i$ |
| | | | | 1410 | 125 275.610717111 + 187.977046429 $i$ |
| | | | | 1411 | 125 364.811579815 + 187.860970805 $i$ |
| | | | | 1412 | 125 453.593193025 + 188.124137630 $i$ |
| | | | | 1413 | 125 542.258889098 + 188.065442418 $i$ |
| | | | | 1414 | 125 631.233274633 + 187.865814214 $i$ |
| | | | | 1415 | 125 720.241220982 + 188.144425171 $i$ |
| | | | | 1416 | 125 808.744862693 + 188.137858775 $i$ |
| | | | | 1417 | 125 897.788863076 + 187.859058568 $i$ |
| | | | | 1418 | 125 986.805703510 + 188.206065315 $i$ |
| | | | | 1419 | 126 075.257880993 + 188.144497574 $i$ |
| | | | | 1420 | 126 164.375737970 + 187.917918114 $i$ |
| | | | | 1421 | 126 253.275055042 + 188.213362577 $i$ |
| | | | | 1422 | 126 341.904046123 + 188.182651835 $i$ |
| | | | | 1423 | 126 430.823827278 + 187.987991208 $i$ |
| | | | | 1424 | 126 519.888322804 + 188.172775352 $i$ |
| | | | | 1425 | 126 608.425468663 + 188.298321975 $i$ |
| | | | | 1426 | 126 697.349513786 + 187.957098366 $i$ |
| | | | | 1427 | 126 786.465526208 + 188.241504069 $i$ |
| | | | | 1428 | 126 874.931945039 + 188.313816821 $i$ |
| | | | | 1429 | 126 963.936251964 + 188.005223130 $i$ |
| | | | | 1430 | 127 052.957056727 + 188.276395157 $i$ |
| | | | | 1431 | 127 141.538008074 + 188.312589093 $i$ |
| | | | | 1432 | 127 230.443209886 + 188.112915093 $i$ |
| | | | | 1433 | 127 319.494606004 + 188.199862241 $i$ |
| | | | | 1434 | 127 408.126521827 + 188.450109633 $i$ |
| | | | | 1435 | 127 496.922331753 + 188.076966378 $i$ |
| | | | | 1436 | 127 586.101287724 + 188.271950831 $i$ |
| | | | | 1437 | 127 674.634690535 + 188.458728839 $i$ |
| | | | | 1438 | 127 763.488537826 + 188.113075487 $i$ |
| | | | | 1439 | 127 852.620795045 + 188.324380586 $i$ |
| | | | | 1440 | 127 941.194925567 + 188.437517532 $i$ |

| $k$ | $\sigma_k$ | $k$ | $\sigma_k$ | $k$ | $\sigma_k$ | $k$ | $\sigma_k$ |
| --- | --- | --- | --- | --- | --- | --- | --- |
| 1441 | 128 030.047 266 88 + 188.231 469 990 $i$ | 1481 | 131 583.987 953 048 + 188.565 101 333 $i$ | 1521 | 135 138.109 108 121 + 188.982 600 502 $i$ | 1561 | 138 692.036 085 810 + 189.548 897 408 $i$ |
| 1442 | 128 119.101 941 505 + 188.264 891 356 $i$ | 1482 | 131 673.117 575 420 + 188.827 627 958 $i$ | 1522 | 135 226.886 998 773 + 189.377 046 171 $i$ | 1562 | 138 780.758 708 423 + 189.717 039 946 $i$ |
| 1443 | 128 207.836 587 947 + 188.541 836 972 $i$ | 1483 | 131 761.579 744 634 + 188.911 439 282 $i$ | 1523 | 135 315.529 894 847 + 189.156 498 185 $i$ | 1563 | 138 869.457 986 712 + 189.540 448 933 $i$ |
| 1444 | 128 296.492 615 529 + 188.234 630 353 $i$ | 1484 | 131 850.587 122 668 + 188.574 563 500 $i$ | 1524 | 135 404.582 526 136 + 189.082 686 673 $i$ | 1564 | 138 958.574 926 909 + 189.496 070 528 $i$ |
| 1445 | 128 385.725 913 594 + 188.300 543 675 $i$ | 1485 | 131 939.609 564 481 + 188.885 131 244 $i$ | 1525 | 135 493.482 533 421 + 189.315 858 534 $i$ | 1565 | 139 047.353 137 654 + 189.799 639 125 $i$ |
| 1446 | 128 474.339 582 216 + 188.577 306 260 $i$ | 1486 | 132 028.169 442 088 + 188.894 829 394 $i$ | 1526 | 135 582.064 027 844 + 189.269 865 569 $i$ | 1566 | 139 135.943 922 630 + 189.551 846 586 $i$ |
| 1447 | 128 563.056 071 571 + 188.250 750 676 $i$ | 1487 | 132 117.103 577 417 + 188.678 644 823 $i$ | 1527 | 135 671.094 108 964 + 189.073 164 860 $i$ | 1567 | 139 225.162 847 129 + 189.518 315 241 $i$ |
| 1448 | 128 652.255 486 852 + 188.413 607 872 $i$ | 1488 | 132 206.141 295 853 + 188.843 673 214 $i$ | 1528 | 135 760.082 590 863 + 189.335 301 153 $i$ | 1568 | 139 313.862 240 029 + 189.852 212 810 $i$ |
| 1449 | 128 740.873 437 366 + 188.562 750 793 $i$ | 1489 | 132 294.768 729 042 + 188.971 213 668 $i$ | 1529 | 135 848.556 832 132 + 189.345 333 828 $i$ | 1569 | 139 402.528 619 575 + 189.555 984 564 $i$ |
| 1450 | 128 829.652 818 963 + 188.350 180 690 $i$ | 1490 | 132 383.575 510 121 + 188.715 983 318 $i$ | 1530 | 135 937.673 794 652 + 189.046 098 584 $i$ | 1570 | 139 491.680 177 086 + 189.588 161 564 $i$ |
| 1451 | 128 918.712 487 621 + 188.339 518 955 $i$ | 1491 | 132 472.747 454 830 + 188.823 287 517 $i$ | 1531 | 136 026.594 058 535 + 189.431 340 672 $i$ | 1571 | 139 580.407 159 140 + 189.802 381 809 $i$ |
| 1452 | 129 007.522 096 483 + 188.607 736 787 $i$ | 1492 | 132 561.270 654 170 + 189.079 484 791 $i$ | 1532 | 136 115.141 192 091 + 189.304 754 080 $i$ | 1572 | 139 669.103 499 852 + 189.679 528 652 $i$ |
| 1453 | 129 096.109 939 394 + 188.413 607 872 $i$ | 1493 | 132 650.155 810 741 + 188.671 635 056 $i$ | 1533 | 136 204.191 082 976 + 189.158 288 673 $i$ | 1573 | 139 758.151 513 984 + 189.549 149 503 $i$ |
| 1454 | 129 185.308 679 831 + 188.312 783 257 $i$ | 1494 | 132 739.258 514 471 + 188.924 327 361 $i$ | 1534 | 136 293.128 378 764 + 189.373 399 644 $i$ | 1574 | 139 847.046 239 352 + 189.864 398 223 $i$ |
| 1455 | 129 274.040 490 326 + 188.699 902 410 $i$ | 1495 | 132 827.841 954 906 + 189.033 232 700 $i$ | 1535 | 136 381.735 605 129 + 189.408 161 786 $i$ | 1575 | 139 935.569 192 060 + 189.711 991 486 $i$ |
| 1456 | 129 362.667 584 515 + 188.397 037 048 $i$ | 1496 | 132 916.693 684 428 + 188.781 856 504 $i$ | 1536 | 136 470.668 254 985 + 189.162 434 056 $i$ | 1576 | 140 024.745 790 724 + 189.553 91 007 $i$ |
| 1457 | 129 451.850 384 866 + 188.393 028 322 $i$ | 1497 | 133 005.771 227 743 + 188.899 723 434 $i$ | 1537 | 136 559.733 746 551 + 189.377 707 855 $i$ | 1577 | 140 113.563 189 604 + 189.914 150 555 $i$ |
| 1458 | 129 540.569 238 452 + 188.677 367 350 $i$ | 1498 | 133 094.445 590 349 + 189.075 230 044 $i$ | 1538 | 136 648.248 243 827 + 189.490 002 296 $i$ | 1578 | 140 202.127 930 430 + 189.720 904 066 $i$ |
| 1459 | 129 629.254 983 517 + 188.472 494 329 $i$ | 1499 | 133 183.181 579 840 + 188.857 500 541 $i$ | 1539 | 136 737.216 817 991 + 189.132 925 597 $i$ | 1579 | 140 291.293 708 672 + 189.620 601 754 $i$ |
| 1460 | 129 718.326 303 297 + 188.424 973 561 $i$ | 1500 | 133 272.342 841 234 + 188.850 789 229 $i$ | 1540 | 136 826.280 708 343 + 189.462 910 017 $i$ | 1580 | 140 380.067 971 156 + 189.896 215 634 $i$ |
| 1461 | 129 807.188 754 179 + 188.655 355 900 $i$ | 1501 | 133 360.994 625 353 + 189.198 386 172 $i$ | 1541 | 136 914.781 468 092 + 189.466 114 423 $i$ | 1581 | 140 468.757 313 053 + 189.797 205 556 $i$ |
| 1462 | 129 895.750 714 507 + 188.591 801 986 $i$ | 1502 | 133 449.710 296 570 + 188.804 970 992 $i$ | 1542 | 137 003.790 336 767 + 189.226 249 328 $i$ | 1582 | 140 557.722 516 934 + 189.630 506 934 $i$ |
| 1463 | 129 984.874 596 442 + 188.362 023 152 $i$ | 1503 | 133 538.903 906 046 + 188.957 732 855 $i$ | 1543 | 137 092.761 785 314 + 189.440 923 965 $i$ | 1583 | 140 646.726 687 116 + 189.908 687 397 $i$ |
| 1464 | 130 073.754 355 367 + 188.774 042 221 $i$ | 1504 | 133 627.521 383 909 + 189.152 758 564 $i$ | 1544 | 137 181.425 299 420 + 189.521 866 862 $i$ | 1584 | 140 735.221 921 419 + 189.865 304 454 $i$ |
| 1465 | 130 162.273 959 432 + 188.554 272 344 $i$ | 1505 | 133 716.288 366 358 + 188.899 067 658 $i$ | 1545 | 137 270.236 223 452 + 189.269 597 530 $i$ | 1585 | 140 824.303 728 863 + 189.619 561 281 $i$ |
| 1466 | 130 251.437 793 389 + 188.453 153 846 $i$ | 1506 | 133 805.388 553 019 + 188.958 570 108 $i$ | 1546 | 137 359.372 957 270 + 189.421 357 696 $i$ | 1586 | 140 913.259 756 527 + 189.966 402 650 $i$ |
| 1467 | 130 340.271 817 695 + 188.755 925 633 $i$ | 1507 | 133 894.127 820 214 + 189.163 897 398 $i$ | 1547 | 137 447.942 369 524 + 189.606 427 990 $i$ | 1587 | 141 001.760 554 238 + 189.877 925 673 $i$ |
| 1468 | 130 428.860 563 379 + 188.610 489 961 $i$ | 1508 | 133 982.791 976 259 + 189.000 819 582 $i$ | 1548 | 137 536.772 318 960 + 189.266 608 571 $i$ | 1588 | 141 090.878 183 643 + 189.672 717 060 $i$ |
| 1469 | 130 517.938 384 241 + 188.509 467 519 $i$ | 1509 | 134 071.937 188 117 + 188.908 234 867 $i$ | 1549 | 137 625.946 198 434 + 189.466 830 983 $i$ | 1589 | 141 179.747 009 253 + 189.964 742 196 $i$ |
| 1470 | 130 606.837 840 971 + 188.707 943 388 $i$ | 1510 | 134 160.709 591 143 + 189.260 182 888 $i$ | 1550 | 137 714.440 507 447 + 189.625 484 034 $i$ | 1590 | 141 268.387 557 952 + 189.915 500 085 $i$ |
| 1471 | 130 695.420 881 986 + 188.751 281 111 $i$ | 1511 | 134 249.285 375 916 + 188.984 437 468 $i$ | 1551 | 137 803.384 724 861 + 189.314 724 254 $i$ | 1591 | 141 357.325 390 838 + 189.737 863 948 $i$ |
| 1472 | 130 784.432 436 603 + 188.440 329 715 $i$ | 1512 | 134 338.529 887 399 + 188.968 511 216 $i$ | 1552 | 137 892.399 245 073 + 189.493 823 315 $i$ | 1592 | 141 446.369 19 233 + 189.925 048 08 $i$ |
| 1473 | 130 873.443 407 541 + 188.812 413 068 $i$ | 1513 | 134 427.196 353 607 + 189.270 853 634 $i$ | 1553 | 137 981.093 599 274 + 189.625 293 027 $i$ | 1593 | 141 534.898 978 92 + 190.027 269 16 $i$ |
| 1474 | 130 961.911 379 213 + 188.737 539 962 $i$ | 1514 | 134 515.908 273 291 + 189.024 932 892 $i$ | 1554 | 138 069.837 566 257 + 189.405 700 565 $i$ | 1594 | 141 623.866 918 39 + 189.697 439 89 $i$ |
| 1475 | 131 051.020 487 563 + 188.505 980 717 $i$ | 1515 | 134 604.984 118 517 + 189.013 145 982 $i$ | 1555 | 138 158.988 464 434 + 189.449 562 150 $i$ | 1595 | 141 712.932 517 46 + 189.998 104 63 $i$ |
| 1476 | 131 139.944 731 828 + 188.826 587 981 $i$ | 1516 | 134 693.810 352 979 + 189.252 850 423 $i$ | 1556 | 138 247.639 562 312 + 189.714 864 304 $i$ | 1596 | 141 801.416 941 89 + 190.029 016 42 $i$ |
| 1477 | 131 228.507 226 982 + 188.758 478 877 $i$ | 1517 | 134 782.426 092 463 + 189.131 520 007 $i$ | 1557 | 138 336.356 080 889 + 189.409 761 035 $i$ | 1597 | 141 890.444 627 30 + 189.751 854 81 $i$ |
| 1478 | 131 317.525 787 559 + 188.583 750 922 $i$ | 1518 | 134 871.510 594 699 + 188.979 617 183 $i$ | 1558 | 138 425.568 939 179 + 189.474 880 076 $i$ | 1598 | 141 979.426 800 41 + 190.018 710 67 $i$ |
| 1479 | 131 406.486 413 835 + 188.781 549 515 $i$ | 1519 | 134 960.406 299 687 + 189.310 595 341 $i$ | 1559 | 138 514.141 040 155 + 189.758 872 640 $i$ | 1599 | 142 068.030 080 51 + 190.039 205 28 $i$ |
| 1480 | 131 495.102 202 135 + 188.870 515 764 $i$ | 1520 | 135 048.911 803 242 + 189.168 033 094 $i$ | 1560 | 138 602.957 098 318 + 189.415 151 556 $i$ | 1600 | 142 156.938 184 13 + 189.834 751 50 $i$ |

| $k$ | $\sigma_k$ | $k$ | $\sigma_k$ | $k$ | $\sigma_k$ |
|---|---|---|---|---|---|
| 1601 | 142 245.97721038 + 189.95931465 $i$ | 1641 | 145 799.91595471 + 190.51438382 $i$ | 1681 | 149 353.66535097 + 190.91032410 $i$ |
| 1602 | 142 334.61232785 + 190.15979081 $i$ | 1642 | 145 888.40015777 + 190.39414047 $i$ | 1682 | 149 442.41228876 + 190.65965954 $i$ |
| 1603 | 142 423.41678754 + 189.79572913 $i$ | 1643 | 145 977.52714238 + 190.20584882 $i$ | 1683 | 149 531.47326019 + 190.66823125 $i$ |
| 1604 | 142 512.58763661 + 190.03271676 $i$ | 1644 | 146 066.40607812 + 190.50967451 $i$ | 1684 | 149 620.28791429 + 190.89144403 $i$ |
| 1605 | 142 601.10598176 + 190.15770754 $i$ | 1645 | 146 155.00360912 + 190.43285476 $i$ | 1685 | 149 708.91383112 + 190.76216024 $i$ |
| 1606 | 142 690.00184987 + 189.84783152 $i$ | 1646 | 146 244.02047721 + 190.26707468 $i$ | 1686 | 149 798.01871915 + 190.62577369 $i$ |
| 1607 | 142 779.09221359 + 190.05726128 $i$ | 1647 | 146 332.96999026 + 190.47510425 $i$ | 1687 | 149 886.86584393 + 190.95640773 $i$ |
| 1608 | 142 867.68138082 + 190.15767894 $i$ | 1648 | 146 421.58148787 + 190.53022596 $i$ | 1688 | 149 975.41337694 + 190.78046579 $i$ |
| 1609 | 142 956.54527127 + 189.94180870 $i$ | 1649 | 146 510.49432657 + 190.24401856 $i$ | 1689 | 150 064.60292385 + 190.64320697 $i$ |
| 1610 | 143 045.58834527 + 190.00870730 $i$ | 1650 | 146 599.59095393 + 190.52190589 $i$ | 1690 | 150 153.36338558 + 191.00548734 $i$ |
| 1611 | 143 134.31516909 + 190.25215455 $i$ | 1651 | 146 688.05722074 + 190.56498396 $i$ | 1691 | 150 242.01680457 + 190.77253130 $i$ |
| 1612 | 143 222.99339759 + 189.93987425 $i$ | 1652 | 146 777.09881110 + 190.26178011 $i$ | 1692 | 150 331.08136704 + 190.73914689 $i$ |
| 1613 | 143 312.21805003 + 190.03675112 $i$ | 1653 | 146 866.07954062 + 190.56747560 $i$ | 1693 | 150 419.95991650 + 190.93881233 $i$ |
| 1614 | 143 400.79331168 + 190.28481192 $i$ | 1654 | 146 954.65509434 + 190.55336699 $i$ | 1694 | 150 508.54804775 + 190.89852177 $i$ |
| 1615 | 143 489.58326728 + 189.96706532 $i$ | 1655 | 147 043.60690441 + 190.35820253 $i$ | 1695 | 150 597.60381091 + 190.69652273 $i$ |
| 1616 | 143 578.72940512 + 190.08096270 $i$ | 1656 | 147 132.62212082 + 190.51673415 $i$ | 1696 | 150 686.54222415 + 190.98701415 $i$ |
| 1617 | 143 667.34962802 + 190.27513253 $i$ | 1657 | 147 221.24371547 + 190.64199715 $i$ | 1697 | 150 775.05950123 + 190.94616234 $i$ |
| 1618 | 143 756.15510210 + 190.05334398 $i$ | 1658 | 147 310.08846577 + 190.36761634 $i$ | 1698 | 150 864.17036583 + 190.68079507 $i$ |
| 1619 | 143 845.20178248 + 190.06153719 $i$ | 1659 | 147 399.21946084 + 190.51831151 $i$ | 1699 | 150 953.06011846 + 191.06898878 $i$ |
| 1620 | 143 933.99525385 + 190.31924617 $i$ | 1660 | 147 487.75537401 + 190.72474504 $i$ | 1700 | 151 041.63937586 + 190.91069173 $i$ |
| 1621 | 144 022.60697680 + 190.10293889 $i$ | 1661 | 147 576.65613033 + 190.33065612 $i$ | 1701 | 151 130.68856448 + 190.78895194 $i$ |
| 1622 | 144 111.80684926 + 190.04312238 $i$ | 1662 | 147 665.73605266 + 190.61603850 $i$ | 1702 | 151 219.59956836 + 191.00048571 $i$ |
| 1623 | 144 200.50449943 + 190.39806516 $i$ | 1663 | 147 754.32466465 + 190.67310995 $i$ | 1703 | 151 308.22531792 + 191.02465482 $i$ |
| 1624 | 144 289.16875317 + 190.08757536 $i$ | 1664 | 147 843.19462663 + 190.45148642 $i$ | 1704 | 151 397.16965891 + 190.77193395 $i$ |
| 1625 | 144 378.34187752 + 190.12692107 $i$ | 1665 | 147 932.25575495 + 190.56520067 $i$ | 1705 | 151 486.20469015 + 191.02309536 $i$ |
| 1626 | 144 467.04755953 + 190.36346145 $i$ | 1666 | 148 020.92298600 + 190.73811186 $i$ | 1706 | 151 574.73241971 + 191.08807227 $i$ |
| 1627 | 144 555.74588867 + 190.16951594 $i$ | 1667 | 148 109.68491012 + 190.49364807 $i$ | 1707 | 151 663.72997996 + 190.75568027 $i$ |
| 1628 | 144 644.82051329 + 190.13453998 $i$ | 1668 | 148 198.82804533 + 190.53785416 $i$ | 1708 | 151 752.74755636 + 191.08969486 $i$ |
| 1629 | 144 733.66246126 + 190.36664026 $i$ | 1669 | 148 287.46949917 + 190.86372896 $i$ | 1709 | 151 841.26734324 + 191.06229845 $i$ |
| 1630 | 144 822.24528041 + 190.26689324 $i$ | 1670 | 148 376.21164437 + 190.45176455 $i$ | 1710 | 151 930.28976515 + 190.85304055 $i$ |
| 1631 | 144 911.37381724 + 190.07862917 $i$ | 1671 | 148 465.39274752 + 190.63480666 $i$ | 1711 | 152 019.24657439 + 191.06085379 $i$ |
| 1632 | 145 000.22345160 + 190.47415217 $i$ | 1672 | 148 553.98743985 + 190.78621799 $i$ | 1712 | 152 107.90133808 + 191.12029641 $i$ |
| 1633 | 145 088.76841555 + 190.22782494 $i$ | 1673 | 148 642.79506611 + 190.56052135 $i$ | 1713 | 152 196.73830899 + 190.88493496 $i$ |
| 1634 | 145 177.94069221 + 190.16920821 $i$ | 1674 | 148 731.87709953 + 190.61154411 $i$ | 1714 | 152 285.86218837 + 191.04120881 $i$ |
| 1635 | 145 266.73059418 + 190.44145455 $i$ | 1675 | 148 820.60160642 + 190.81923456 $i$ | 1715 | 152 374.40356083 + 191.20682390 $i$ |
| 1636 | 145 355.36656058 + 190.30660418 $i$ | 1676 | 148 909.29396174 + 190.63093437 $i$ | 1716 | 152 463.29747174 + 190.87215775 $i$ |
| 1637 | 145 444.43249881 + 190.18945622 $i$ | 1677 | 148 998.42942014 + 190.56870846 $i$ | 1717 | 152 552.41301484 + 191.09028306 $i$ |
| 1638 | 145 533.30597868 + 190.42251489 $i$ | 1678 | 149 087.16895829 + 190.91083949 $i$ | 1718 | 152 640.92961937 + 191.21364313 $i$ |
| 1639 | 145 621.91956854 + 190.41800742 $i$ | 1679 | 149 175.79838904 + 190.61372815 $i$ | 1719 | 152 729.88152716 + 190.91563284 $i$ |
| 1640 | 145 710.93216152 + 190.13303642 $i$ | 1680 | 149 265.01497601 + 190.63169849 $i$ | 1720 | 152 818.88116154 + 191.11522056 $i$ |
| | | | | 1721 | 152 907.57427509 + 191.21339833 $i$ |
| | | | | 1722 | 152 996.33252735 + 191.00414195 $i$ |
| | | | | 1723 | 153 085.47866579 + 191.06321153 $i$ |
| | | | | 1724 | 153 174.10437245 + 191.31559714 $i$ |
| | | | | 1725 | 153 262.86956634 + 190.98929453 $i$ |
| | | | | 1726 | 153 352.03993262 + 191.10613869 $i$ |
| | | | | 1727 | 153 440.62695253 + 191.34256141 $i$ |
| | | | | 1728 | 153 529.45526057 + 191.00558935 $i$ |
| | | | | 1729 | 153 618.52319665 + 191.16495878 $i$ |
| | | | | 1730 | 153 707.23799090 + 191.30092685 $i$ |
| | | | | 1731 | 153 795.95475081 + 191.12949608 $i$ |
| | | | | 1732 | 153 885.07345656 + 191.09341667 $i$ |
| | | | | 1733 | 153 973.81167259 + 191.39408230 $i$ |
| | | | | 1734 | 154 062.45304418 + 191.12867941 $i$ |
| | | | | 1735 | 154 151.65990617 + 191.12288218 $i$ |
| | | | | 1736 | 154 240.31809395 + 191.42794093 $i$ |
| | | | | 1737 | 154 329.03688059 + 191.14313199 $i$ |
| | | | | 1738 | 154 418.16693692 + 191.18452829 $i$ |
| | | | | 1739 | 154 506.88357343 + 191.38929318 $i$ |
| | | | | 1740 | 154 595.60004518 + 191.25118491 $i$ |
| | | | | 1741 | 154 684.64842095 + 191.14024909 $i$ |
| | | | | 1742 | 154 773.51345249 + 191.45103681 $i$ |
| | | | | 1743 | 154 862.05966815 + 191.27786552 $i$ |
| | | | | 1744 | 154 951.25281862 + 191.15646953 $i$ |
| | | | | 1745 | 155 040.02057968 + 191.49582971 $i$ |
| | | | | 1746 | 155 128.63911836 + 191.28260614 $i$ |
| | | | | 1747 | 155 217.77356527 + 191.21204429 $i$ |
| | | | | 1748 | 155 306.55328918 + 191.48080126 $i$ |
| | | | | 1749 | 155 395.24421069 + 191.35246800 $i$ |
| | | | | 1750 | 155 484.22274476 + 191.21968888 $i$ |
| | | | | 1751 | 155 573.20300724 + 191.48302133 $i$ |
| | | | | 1752 | 155 661.69839217 + 191.42661615 $i$ |
| | | | | 1753 | 155 750.82375101 + 191.20056305 $i$ |
| | | | | 1754 | 155 839.71435882 + 191.54324406 $i$ |
| | | | | 1755 | 155 928.26136488 + 191.43016065 $i$ |
| | | | | 1756 | 156 017.37019617 + 191.25338638 $i$ |
| | | | | 1757 | 156 106.22185814 + 191.54009154 $i$ |
| | | | | 1758 | 156 194.87614845 + 191.47371782 $i$ |
| | | | | 1759 | 156 283.82808672 + 191.30522718 $i$ |
| | | | | 1760 | 156 372.84377086 + 191.50216104 $i$ |

| $k$ | $\sigma_k$ | $k$ | $\sigma_k$ | $k$ | $\sigma_k$ |
|---|---|---|---|---|---|
| 1761 | 156 461.37863996 + 191.58069646 $i$ | 1801 | 160 015.26723037 + 191.74740266 $i$ | 1841 | 163 569.30448785 + 192.05498342 $i$ |
| 1762 | 156 550.37923553 + 191.25944706 $i$ | 1802 | 160 104.44663642 + 191.70962258 $i$ | 1842 | 163 658.35669301 + 192.10480638 $i$ |
| 1763 | 156 639.39779648 + 191.58264125 $i$ | 1803 | 160 193.18956017 + 191.95914125 $i$ | 1843 | 163 747.07816077 + 192.31766782 $i$ |
| 1764 | 156 727.91107066 + 191.56772212 $i$ | 1804 | 160 281.87109093 + 191.83038580 $i$ | 1844 | 163 835.79146914 + 192.10159077 $i$ |
| 1765 | 156 816.94705337 + 191.31736350 $i$ | 1805 | 160 370.92205278 + 191.70905798 $i$ | 1845 | 163 924.92245715 + 192.07955075 $i$ |
| 1766 | 156 905.89559623 + 191.58665344 $i$ | 1806 | 160 459.78432381 + 191.96212713 $i$ | 1846 | 164 013.63673981 + 192.39513102 $i$ |
| 1767 | 156 994.52171595 + 191.58971747 $i$ | 1807 | 160 548.40113696 + 191.91212089 $i$ | 1847 | 164 102.30482381 + 192.08379223 $i$ |
| 1768 | 157 083.43225009 + 191.38508211 $i$ | 1808 | 160 637.43265911 + 191.67761130 $i$ | 1848 | 164 191.49983720 + 192.14094235 $i$ |
| 1769 | 157 172.46153170 + 191.53692354 $i$ | 1809 | 160 726.39282462 + 192.04099602 $i$ | 1849 | 164 280.13980602 + 192.38586322 $i$ |
| 1770 | 157 261.08478866 + 191.69839127 $i$ | 1810 | 160 814.88335099 + 191.89145442 $i$ | 1850 | 164 368.91147304 + 192.13756636 $i$ |
| 1771 | 157 349.92355176 + 191.35515654 $i$ | 1811 | 160 904.03460115 + 191.74586489 $i$ | 1851 | 164 457.96461127 + 192.16316715 $i$ |
| 1772 | 157 439.07450201 + 191.60226485 $i$ | 1812 | 160 992.87337026 + 192.02685691 $i$ | 1852 | 164 546.76156165 + 192.37019724 $i$ |
| 1773 | 157 527.57078160 + 191.69039369 $i$ | 1813 | 161 081.50082362 + 191.93890597 $i$ | 1853 | 164 635.40187155 + 192.23710590 $i$ |
| 1774 | 157 616.52107789 + 191.41070100 $i$ | 1814 | 161 170.51103380 + 191.78545026 $i$ | 1854 | 164 724.52397080 + 192.11576789 $i$ |
| 1775 | 157 705.56199931 + 191.61069227 $i$ | 1815 | 161 259.45257927 + 192.00428124 $i$ | 1855 | 164 813.32593918 + 192.44411161 $i$ |
| 1776 | 157 794.16626853 + 191.70485568 $i$ | 1816 | 161 348.06036048 + 192.02861545 $i$ | 1856 | 164 901.91488481 + 192.23933554 $i$ |
| 1777 | 157 883.04238489 + 191.48338284 $i$ | 1817 | 161 437.00483579 + 191.76006820 $i$ | 1857 | 164 991.09775621 + 192.14673197 $i$ |
| 1778 | 157 972.08018512 + 191.57295278 $i$ | 1818 | 161 526.06130319 + 192.04815034 $i$ | 1858 | 165 079.83585285 + 192.47428168 $i$ |
| 1779 | 158 060.78213453 + 191.78565089 $i$ | 1819 | 161 614.54110204 + 192.05680150 $i$ | 1859 | 165 168.50961667 + 192.24259920 $i$ |
| 1780 | 158 149.49535377 + 191.48388366 $i$ | 1820 | 161 703.60874723 + 191.77762330 $i$ | 1860 | 165 257.58259416 + 192.22369957 $i$ |
| 1781 | 158 238.70980376 + 191.60036060 $i$ | 1821 | 161 792.54207222 + 192.08491731 $i$ | 1861 | 165 346.42371451 + 192.41196856 $i$ |
| 1782 | 158 327.25561823 + 191.81941543 $i$ | 1822 | 161 881.14602661 + 192.05447374 $i$ | 1862 | 165 435.04212945 + 192.37533911 $i$ |
| 1783 | 158 416.10220779 + 191.50483444 $i$ | 1823 | 161 970.10925935 + 191.86624315 $i$ | 1863 | 165 524.10926797 + 192.16133756 $i$ |
| 1784 | 158 505.19813722 + 191.63873295 $i$ | 1824 | 162 059.09740657 + 192.03035682 $i$ | 1864 | 165 613.00420882 + 192.48263220 $i$ |
| 1785 | 158 593.83721492 + 191.81256732 $i$ | 1825 | 162 147.72314025 + 192.14778761 $i$ | 1865 | 165 701.55773893 + 192.38942009 $i$ |
| 1786 | 158 682.64574803 + 191.58327099 $i$ | 1826 | 162 236.59557303 + 191.85865931 $i$ | 1866 | 165 790.66720388 + 192.17037695 $i$ |
| 1787 | 158 771.69380472 + 191.62206897 $i$ | 1827 | 162 325.69522072 + 192.05164363 $i$ | 1867 | 165 879.53158366 + 192.54076910 $i$ |
| 1788 | 158 860.47080934 + 191.85489737 $i$ | 1828 | 162 414.23548134 + 192.20297121 $i$ | 1868 | 165 968.12829358 + 192.36605266 $i$ |
| 1789 | 158 949.10132543 + 191.62536699 $i$ | 1829 | 162 503.15812995 + 191.84101797 $i$ | 1869 | 166 057.18922447 + 192.27082858 $i$ |
| 1790 | 159 038.30437026 + 191.60544292 $i$ | 1830 | 162 592.22000807 + 192.13269578 $i$ | 1870 | 166 146.07511779 + 192.46813342 $i$ |
| 1791 | 159 126.97046997 + 191.92590456 $i$ | 1831 | 162 680.79879015 + 192.15398699 $i$ | 1871 | 166 234.70761213 + 192.48418083 $i$ |
| 1792 | 159 215.67521117 + 191.61092050 $i$ | 1832 | 162 769.69951147 + 191.96280362 $i$ | 1872 | 166 323.67107771 + 192.23939494 $i$ |
| 1793 | 159 304.82872711 + 191.68001218 $i$ | 1833 | 162 858.74013717 + 192.06495546 $i$ | 1873 | 166 412.68535080 + 192.50771466 $i$ |
| 1794 | 159 393.51834397 + 191.88841732 $i$ | 1834 | 162 947.39691469 + 192.23743343 $i$ | 1874 | 166 501.20855839 + 192.52469219 $i$ |
| 1795 | 159 482.24723856 + 191.70255386 $i$ | 1835 | 163 036.18595617 + 191.97610948 $i$ | 1875 | 166 590.23673682 + 192.23718388 $i$ |
| 1796 | 159 571.31498222 + 191.66426525 $i$ | 1836 | 163 125.31847816 + 192.06091497 $i$ | 1876 | 166 679.22178157 + 192.56333754 $i$ |
| 1797 | 159 660.12830453 + 191.90909536 $i$ | 1837 | 163 213.93531252 + 192.31029885 $i$ | 1877 | 166 767.75785217 + 192.50422279 $i$ |
| 1798 | 159 748.74317462 + 191.77961200 $i$ | 1838 | 163 302.71553964 + 191.95243270 $i$ | 1878 | 166 856.78934412 + 192.32115983 $i$ |
| 1799 | 159 837.87651367 + 191.62108647 $i$ | 1839 | 163 391.88364449 + 192.14464297 $i$ | 1879 | 166 945.72556288 + 192.52330803 $i$ |
| 1800 | 159 926.68175825 + 192.00098932 $i$ | 1840 | 163 480.45585045 + 192.26450927 $i$ | 1880 | 167 034.37667422 + 192.57173096 $i$ |
| 1881 | 167 123.24255323 + 192.34294536 $i$ | | | | |
| 1882 | 167 212.34223775 + 192.51338059 $i$ | | | | |
| 1883 | 167 300.87806727 + 192.65554465 $i$ | | | | |
| 1884 | 167 389.81322957 + 192.31654120 $i$ | | | | |
| 1885 | 167 478.87214302 + 192.57014163 $i$ | | | | |
| 1886 | 167 567.42293428 + 192.63563450 $i$ | | | | |
| 1887 | 167 656.37814313 + 192.37140377 $i$ | | | | |
| 1888 | 167 745.36621201 + 192.57949481 $i$ | | | | |
| 1889 | 167 834.04717377 + 192.65626968 $i$ | | | | |
| 1890 | 167 922.83588930 + 192.45804233 $i$ | | | | |
| 1891 | 168 011.97051331 + 192.52220420 $i$ | | | | |
| 1892 | 168 100.57017868 + 192.76070110 $i$ | | | | |
| 1893 | 168 189.38055400 + 192.42370705 $i$ | | | | |
| 1894 | 168 278.51903130 + 192.58119416 $i$ | | | | |
| 1895 | 168 367.10390522 + 192.76543697 $i$ | | | | |
| 1896 | 168 455.95133071 + 192.46076733 $i$ | | | | |
| 1897 | 168 545.01972208 + 192.61881835 $i$ | | | | |
| 1898 | 168 633.70400956 + 192.73243319 $i$ | | | | |
| 1899 | 168 722.45141858 + 192.57704532 $i$ | | | | |
| 1900 | 168 811.57206122 + 192.54361750 $i$ | | | | |
| 1901 | 168 900.27347104 + 192.83743252 $i$ | | | | |
| 1902 | 168 988.95774121 + 192.55532464 $i$ | | | | |
| 1903 | 169 078.14990968 + 192.58821205 $i$ | | | | |
| 1904 | 169 166.78714012 + 192.85787362 $i$ | | | | |
| 1905 | 169 255.54368952 + 192.57501769 $i$ | | | | |
| 1906 | 169 344.64855235 + 192.63554681 $i$ | | | | |
| 1907 | 169 433.36471165 + 192.82645558 $i$ | | | | |
| 1908 | 169 522.08896846 + 192.67372523 $i$ | | | | |
| 1909 | 169 611.14390516 + 192.59284236 $i$ | | | | |
| 1910 | 169 699.98648600 + 192.89168891 $i$ | | | | |
| 1911 | 169 788.55681195 + 192.69210884 $i$ | | | | |
| 1912 | 169 877.75134563 + 192.60385514 $i$ | | | | |
| 1913 | 169 966.47924388 + 192.93183886 $i$ | | | | |
| 1914 | 170 055.14950490 + 192.70045826 $i$ | | | | |
| 1915 | 170 144.25949230 + 192.65622688 $i$ | | | | |
| 1916 | 170 233.03405422 + 192.90583960 $i$ | | | | |
| 1917 | 170 321.72180149 + 192.77326974 $i$ | | | | |
| 1918 | 170 410.73191140 + 192.66093190 $i$ | | | | |
| 1919 | 170 499.66743479 + 192.91085970 $i$ | | | | |
| 1920 | 170 588.18419474 + 192.85111300 $i$ | | | | |

| $k$ | $\sigma_k$ | $k$ | $\sigma_k$ | $k$ | $\sigma_k$ |
|---|---|---|---|---|---|
| 1921 | 170 677.34080121 + 192.63107455 $i$ | 1961 | 174 231.31686530 + 193.09113755 $i$ | 2001 | 177 785.21857303 + 193.43677718 $i$ |
| 1922 | 170 766.16847461 + 192.97698943 $i$ | 1962 | 174 319.98980879 + 193.27387024 $i$ | 2002 | 177 873.87659584 + 193.60391366 $i$ |
| 1923 | 170 854.76340733 + 192.83783654 $i$ | 1963 | 174 408.74531557 + 193.09643946 $i$ | 2003 | 177 962.68467904 + 193.32516904 $i$ |
| 1924 | 170 943.86127692 + 192.69319065 $i$ | 1964 | 174 497.80853279 + 193.06114490 $i$ | 2004 | 178 051.80556604 + 193.45264193 $i$ |
| 1925 | 171 032.70266222 + 192.96502069 $i$ | 1965 | 174 586.59969467 + 193.30913876 $i$ | 2005 | 178 140.40342084 + 193.65594126 $i$ |
| 1926 | 171 121.35894030 + 192.88520252 $i$ | 1966 | 174 675.23466477 + 193.15144367 $i$ | 2006 | 178 229.22829315 + 193.31910525 $i$ |
| 1927 | 171 210.32932811 + 192.73216842 $i$ | 1967 | 174 764.37610057 + 193.03489475 $i$ | 2007 | 178 318.36439268 + 193.51295199 $i$ |
| 1928 | 171 299.32115993 + 192.93593195 $i$ | 1968 | 174 853.14743821 + 193.38842785 $i$ | 2008 | 178 406.92574702 + 193.61908616 $i$ |
| 1929 | 171 387.86113017 + 192.98406373 $i$ | 1969 | 174 941.76712099 + 193.12484773 $i$ | 2009 | 178 495.81256832 + 193.41416710 $i$ |
| 1930 | 171 476.88637707 + 192.68073417 $i$ | 1970 | 175 030.94129694 + 193.11432457 $i$ | 2010 | 178 584.83537247 + 193.47242724 $i$ |
| 1931 | 171 565.86718800 + 193.01910891 $i$ | 1971 | 175 119.65834383 + 193.34564931 $i$ | 2011 | 178 673.55957138 + 193.67903324 $i$ |
| 1932 | 171 654.40005841 + 192.95870105 $i$ | 1972 | 175 208.37169066 + 193.20710151 $i$ | 2012 | 178 762.28312048 + 193.44294589 $i$ |
| 1933 | 171 743.44414071 + 192.74875956 $i$ | 1973 | 175 297.40493561 + 193.10364734 $i$ | 2013 | 178 851.41990141 + 193.46146668 $i$ |
| 1934 | 171 832.37307440 + 193.00846322 $i$ | 1974 | 175 386.27277969 + 193.36338731 $i$ | 2014 | 178 940.10063989 + 193.73958047 $i$ |
| 1935 | 171 921.00612109 + 192.99335142 $i$ | 1975 | 175 474.88162257 + 193.26830537 $i$ | 2015 | 179 028.80720982 + 193.43461024 $i$ |
| 1936 | 172 009.92625528 + 192.80258531 $i$ | 1976 | 175 563.93826134 + 193.08232574 $i$ | 2016 | 179 117.99181797 + 193.51789068 $i$ |
| 1937 | 172 098.95204330 + 192.96709144 $i$ | 1977 | 175 652.85829837 + 193.42405076 $i$ | 2017 | 179 206.61084935 + 193.72351178 $i$ |
| 1938 | 172 187.55389193 + 193.09283600 $i$ | 1978 | 175 741.37234114 + 193.26413895 $i$ | 2018 | 179 295.40487010 + 193.49413926 $i$ |
| 1939 | 172 276.43430642 + 192.77263650 $i$ | 1979 | 175 830.54258812 + 193.13969563 $i$ | 2019 | 179 384.46160513 + 193.52943924 $i$ |
| 1940 | 172 365.55193218 + 193.02484669 $i$ | 1980 | 175 919.33351639 + 193.40904231 $i$ | 2020 | 179 473.23613360 + 193.71702021 $i$ |
| 1941 | 172 454.04183125 + 193.08810420 $i$ | 1981 | 176 008.00054613 + 193.31542670 $i$ | 2021 | 179 561.89114984 + 193.58013323 $i$ |
| 1942 | 172 543.03652153 + 192.82517248 $i$ | 1982 | 176 097.00620299 + 193.16606261 $i$ | 2022 | 179 651.02569967 + 193.47783860 $i$ |
| 1943 | 172 632.02725551 + 193.02853480 $i$ | 1983 | 176 185.93195735 + 193.39103125 $i$ | 2023 | 179 739.78935483 + 193.79929780 $i$ |
| 1944 | 172 720.65620280 + 193.10996014 $i$ | 1984 | 176 274.53454510 + 193.39557936 $i$ | 2024 | 179 828.41501658 + 193.56810841 $i$ |
| 1945 | 172 809.53423664 + 192.88560687 $i$ | 1985 | 176 363.51886618 + 193.14418842 $i$ | 2025 | 179 917.58990831 + 193.52118257 $i$ |
| 1946 | 172 898.57384038 + 192.99963444 $i$ | 1986 | 176 452.52850773 + 193.43456099 $i$ | 2026 | 180 006.30797062 + 193.81285202 $i$ |
| 1947 | 172 987.25157562 + 193.17913437 $i$ | 1987 | 176 541.02374951 + 193.41652128 $i$ | 2027 | 180 095.00441688 + 193.58610822 $i$ |
| 1948 | 173 076.00280224 + 192.88412802 $i$ | 1988 | 176 630.11569475 + 193.16591584 $i$ | 2028 | 180 184.07880861 + 193.57700793 $i$ |
| 1949 | 173 165.19202490 + 193.01991016 $i$ | 1989 | 176 719.01348974 + 193.46707602 $i$ | 2029 | 180 272.89378205 + 193.76366308 $i$ |
| 1950 | 173 253.72217330 + 193.21645951 $i$ | 1990 | 176 807.63572533 + 193.41391514 $i$ | 2030 | 180 361.53951962 + 193.70894817 $i$ |
| 1951 | 173 342.61651286 + 192.90380321 $i$ | 1991 | 176 896.60356027 + 193.24580840 $i$ | 2031 | 180 450.60082405 + 193.50494225 $i$ |
| 1952 | 173 431.67274574 + 193.05964283 $i$ | 1992 | 176 985.58363250 + 193.41345522 $i$ | 2032 | 180 539.47517577 + 193.85056507 $i$ |
| 1953 | 173 520.32615994 + 193.20064923 $i$ | 1993 | 177 074.20089601 + 193.50474491 $i$ | 2033 | 180 628.05126649 + 193.70178726 $i$ |
| 1954 | 173 609.13194523 + 192.97939450 $i$ | 1994 | 177 163.09222975 + 193.22597493 $i$ | 2034 | 180 717.16670432 + 193.53928895 $i$ |
| 1955 | 173 698.19404220 + 193.03642179 $i$ | 1995 | 177 252.18030035 + 193.45149083 $i$ | 2035 | 180 806.00527508 + 193.87390892 $i$ |
| 1956 | 173 786.93493550 + 193.24475620 $i$ | 1996 | 177 340.71251016 + 193.54068527 $i$ | 2036 | 180 894.61184058 + 193.69742763 $i$ |
| 1957 | 173 875.59903293 + 193.01882796 $i$ | 1997 | 177 429.65955916 + 193.22442955 $i$ | 2037 | 180 983.68830907 + 193.62582865 $i$ |
| 1958 | 173 964.80452208 + 193.02095544 $i$ | 1998 | 177 518.70241720 + 193.50956163 $i$ | 2038 | 181 072.55162985 + 193.81212134 $i$ |
| 1959 | 174 053.43028311 + 193.31108314 $i$ | 1999 | 177 607.27432004 + 193.50743572 $i$ | 2039 | 181 161.19434565 + 193.81452794 $i$ |
| 1960 | 174 142.18160768 + 193.00341999 $i$ | 2000 | 177 696.20406551 + 193.33497898 $i$ | 2040 | 181 250.17221666 + 193.57808279 $i$ |
| | | | | 2041 | 181 339.16205372 + 193.86149058 $i$ |
| | | | | 2042 | 181 427.68560646 + 193.83875131 $i$ |
| | | | | 2043 | 181 516.74688108 + 193.59020745 $i$ |
| | | | | 2044 | 181 605.68949411 + 193.90319955 $i$ |
| | | | | 2045 | 181 694.24618543 + 193.82743543 $i$ |
| | | | | 2046 | 181 783.28917919 + 193.66410165 $i$ |
| | | | | 2047 | 181 872.20387651 + 193.86206864 $i$ |
| | | | | 2048 | 181 960.85725881 + 193.89902540 $i$ |
| | | | | 2049 | 182 049.74751615 + 193.67089389 $i$ |
| | | | | 2050 | 182 138.81944177 + 193.86015695 $i$ |
| | | | | 2051 | 182 227.35313244 + 193.97297084 $i$ |
| | | | | 2052 | 182 316.32059153 + 193.64646494 $i$ |
| | | | | 2053 | 182 405.34655057 + 193.92418378 $i$ |
| | | | | 2054 | 182 493.91138921 + 193.95365913 $i$ |
| | | | | 2055 | 182 582.86600262 + 193.71289957 $i$ |
| | | | | 2056 | 182 671.85911732 + 193.91688305 $i$ |
| | | | | 2057 | 182 760.51850257 + 193.97157251 $i$ |
| | | | | 2058 | 182 849.33783721 + 193.78551644 $i$ |
| | | | | 2059 | 182 938.45865500 + 193.86065272 $i$ |
| | | | | 2060 | 183 027.04085907 + 194.07710138 $i$ |
| | | | | 2061 | 183 115.88551009 + 193.73755479 $i$ |
| | | | | 2062 | 183 205.00160243 + 193.93340449 $i$ |
| | | | | 2063 | 183 293.58253146 + 194.06115110 $i$ |
| | | | | 2064 | 183 382.44767710 + 193.79131471 $i$ |
| | | | | 2065 | 183 471.50755150 + 193.94550530 $i$ |
| | | | | 2066 | 183 560.17426324 + 194.05273961 $i$ |
| | | | | 2067 | 183 648.95638810 + 193.89023881 $i$ |
| | | | | 2068 | 183 738.05805484 + 193.87122361 $i$ |
| | | | | 2069 | 183 826.74167355 + 194.16152008 $i$ |
| | | | | 2070 | 183 915.46387567 + 193.85511275 $i$ |
| | | | | 2071 | 184 004.63577442 + 193.93063330 $i$ |
| | | | | 2072 | 184 093.25882107 + 194.16184124 $i$ |
| | | | | 2073 | 184 182.04598212 + 193.88686450 $i$ |
| | | | | 2074 | 184 271.13268639 + 193.96474739 $i$ |
| | | | | 2075 | 184 359.84481894 + 194.13757944 $i$ |
| | | | | 2076 | 184 448.57751398 + 193.97946685 $i$ |
| | | | | 2077 | 184 537.64661095 + 193.91878425 $i$ |
| | | | | 2078 | 184 626.45092738 + 194.20109386 $i$ |
| | | | | 2079 | 184 715.04980421 + 193.99321533 $i$ |
| | | | | 2080 | 184 804.25318116 + 193.93446644 $i$ |

| $k$ | $\sigma_k$ | $k$ | $\sigma_k$ | $k$ | $\sigma_k$ |
| --- | --- | --- | --- | --- | --- |
| 2081 | 184 892.94269854 + 194.23654754 $i$ | 2121 | 188 446.81196051 + 194.46886228 $i$ | 2161 | 192 000.67941792 + 194.75407502 $i$ |
| 2082 | 184 981.64792747 + 193.99930473 $i$ | 2122 | 188 535.61745965 + 194.26610879 $i$ | 2162 | 192 089.59708622 + 194.48453302 $i$ |
| 2083 | 185 070.74725559 + 193.98839976 $i$ | 2123 | 188 624.69495316 + 194.33149805 $i$ | 2163 | 192 178.66188491 + 194.72819950 $i$ |
| 2084 | 185 159.51757405 + 194.20776240 $i$ | 2124 | 188 713.39684716 + 194.52328054 $i$ | 2164 | 192 267.18477767 + 194.77431406 $i$ |
| 2085 | 185 248.20152920 + 194.07738453 $i$ | 2125 | 188 802.10500607 + 194.29460274 $i$ | 2165 | 192 356.16986225 + 194.49734062 $i$ |
| 2086 | 185 337.24287283 + 193.97780904 $i$ | 2126 | 188 891.29479865 + 194.31804323 $i$ | 2166 | 192 445.17771047 + 194.76655369 $i$ |
| 2087 | 185 426.13334903 + 194.22048710 $i$ | 2127 | 188 979.89554820 + 194.58623442 $i$ | 2167 | 192 533.75291248 + 194.75742267 $i$ |
| 2088 | 185 514.67721132 + 194.14198332 $i$ | 2128 | 189 068.69137011 + 194.27567876 $i$ | 2168 | 192 622.70959243 + 194.58985530 $i$ |
| 2089 | 185 603.83829715 + 193.94377835 $i$ | 2129 | 189 157.79726007 + 194.38504048 $i$ | 2169 | 192 711.69547087 + 194.69871227 $i$ |
| 2090 | 185 692.63441494 + 194.29858290 $i$ | 2130 | 189 246.46543287 + 194.54766023 $i$ | 2170 | 192 800.35550720 + 194.85402525 $i$ |
| 2091 | 185 781.26728054 + 194.11390516 $i$ | 2131 | 189 335.24180491 + 194.37115800 $i$ | 2171 | 192 889.18060984 + 194.57064900 $i$ |
| 2092 | 185 870.34345391 + 194.01610741 $i$ | 2132 | 189 424.29936128 + 194.34787751 $i$ | 2172 | 192 978.29926139 + 194.73095597 $i$ |
| 2093 | 185 959.18229231 + 194.27116280 $i$ | 2133 | 189 513.07600086 + 194.59278180 $i$ | 2173 | 193 066.87004001 + 194.88241982 $i$ |
| 2094 | 186 047.84443046 + 194.17445090 $i$ | 2134 | 189 601.72753261 + 194.40693260 $i$ | 2174 | 193 155.73380493 + 194.57750383 $i$ |
| 2095 | 186 136.82875421 + 194.03401140 $i$ | 2135 | 189 690.87641849 + 194.33212295 $i$ | 2175 | 193 244.84249989 + 194.77754045 $i$ |
| 2096 | 186 225.79298388 + 194.25039854 $i$ | 2136 | 189 779.61227733 + 194.65416199 $i$ | 2176 | 193 333.40753850 + 194.86106577 $i$ |
| 2097 | 186 314.34489011 + 194.26820467 $i$ | 2137 | 189 868.26473270 + 194.38845832 $i$ | 2177 | 193 422.31021586 + 194.65759517 $i$ |
| 2098 | 186 403.39126187 + 193.98639380 $i$ | 2138 | 189 957.43155032 + 194.40549211 $i$ | 2178 | 193 511.31737071 + 194.73960806 $i$ |
| 2099 | 186 492.33627634 + 194.33255713 $i$ | 2139 | 190 046.13181085 + 194.61648229 $i$ | 2179 | 193 600.04394371 + 194.91998363 $i$ |
| 2100 | 186 580.88656935 + 194.23580672 $i$ | 2140 | 190 134.86626156 + 194.46906798 $i$ | 2180 | 193 688.76700129 + 194.67938774 $i$ |
| 2101 | 186 669.94511093 + 194.06506415 $i$ | 2141 | 190 223.89336868 + 194.38925481 $i$ | 2181 | 193 777.92111141 + 194.73759265 $i$ |
| 2102 | 186 758.85117601 + 194.30491105 $i$ | 2142 | 190 312.75869708 + 194.63976405 $i$ | 2182 | 193 866.56721666 + 194.97025788 $i$ |
| 2103 | 186 847.48698044 + 194.28326456 $i$ | 2143 | 190 401.35741787 + 194.52002245 $i$ | 2183 | 193 955.31122824 + 194.67482624 $i$ |
| 2104 | 186 936.42786947 + 194.10220997 $i$ | 2144 | 190 490.45217968 + 194.37233411 $i$ | 2184 | 194 044.47567000 + 194.78206563 $i$ |
| 2105 | 187 025.43901378 + 194.27053677 $i$ | 2145 | 190 579.31906963 + 194.68821223 $i$ | 2185 | 194 133.08524653 + 194.95954226 $i$ |
| 2106 | 187 114.02028121 + 194.37179518 $i$ | 2146 | 190 667.86687005 + 194.52376831 $i$ | 2186 | 194 221.90557605 + 194.73619631 $i$ |
| 2107 | 187 202.94663856 + 194.07242047 $i$ | 2147 | 190 757.04108928 + 194.41683480 $i$ | 2187 | 194 310.94922791 + 194.77927465 $i$ |
| 2108 | 187 292.02521643 + 194.33038931 $i$ | 2148 | 190 845.79913218 + 194.68265231 $i$ | 2188 | 194 399.70778771 + 194.96377491 $i$ |
| 2109 | 187 380.52225566 + 194.36605969 $i$ | 2149 | 190 934.49885582 + 194.57077350 $i$ | 2189 | 194 488.38703310 + 194.81210457 $i$ |
| 2110 | 187 469.54241334 + 194.11900385 $i$ | 2150 | 191 023.49344371 + 194.44018243 $i$ | 2190 | 194 577.52070548 + 194.73153115 $i$ |
| 2111 | 187 558.49893002 + 194.33433448 $i$ | 2151 | 191 112.41780998 + 194.66834467 $i$ | 2191 | 194 666.25498769 + 195.04786245 $i$ |
| 2112 | 187 647.14839547 + 194.38595950 $i$ | 2152 | 191 201.01104986 + 194.64354330 $i$ | 2192 | 194 754.91851198 + 194.78925599 $i$ |
| 2113 | 187 736.01831390 + 194.17341073 $i$ | 2153 | 191 290.02604280 + 194.41781505 $i$ | 2193 | 194 844.07915230 + 194.78387359 $i$ |
| 2114 | 187 825.07066600 + 194.30812239 $i$ | 2154 | 191 378.99859210 + 194.71348261 $i$ | 2194 | 194 932.77779549 + 195.04415227 $i$ |
| 2115 | 187 913.71746562 + 194.45194707 $i$ | 2155 | 191 467.51476795 + 194.65668127 $i$ | 2195 | 195 021.50177928 + 194.82497436 $i$ |
| 2116 | 188 002.50838054 + 194.17166277 $i$ | 2156 | 191 556.61160311 + 194.44014055 $i$ | 2196 | 195 110.57329556 + 194.82088779 $i$ |
| 2117 | 188 091.67474405 + 194.32721956 $i$ | 2157 | 191 645.48780039 + 194.74094233 $i$ | 2197 | 195 199.36746746 + 195.00767384 $i$ |
| 2118 | 188 180.19929590 + 194.48910278 $i$ | 2158 | 191 734.12781452 + 194.65764610 $i$ | 2198 | 195 288.03183166 + 194.93065751 $i$ |
| 2119 | 188 269.12021312 + 194.18318302 $i$ | 2159 | 191 823.09232700 + 194.51236257 $i$ | 2199 | 195 377.09563327 + 194.75284758 $i$ |
| 2120 | 188 358.14897666 + 194.36650505 $i$ | 2160 | 191 912.06727870 + 194.68784462 $i$ | 2200 | 195 465.95386200 + 195.09480416 $i$ |
| 2201 | 195 554.53113777 + 194.90739915 $i$ | 2221 | 197 331.82391734 + 195.16779636 $i$ |  |  |
| 2202 | 195 643.66757464 + 194.80865259 $i$ | 2222 | 197 420.39580082 + 195.15016489 $i$ |  |  |
| 2203 | 195 732.47854515 + 195.09785811 $i$ | 2223 | 197 509.36022398 + 194.95114169 $i$ |  |  |
| 2204 | 195 821.10175819 + 194.92715546 $i$ | 2224 | 197 598.34770734 + 195.13944751 $i$ |  |  |
| 2205 | 195 910.18799152 + 194.86601554 $i$ | 2225 | 197 686.98775077 + 195.18846799 $i$ |  |  |
| 2206 | 195 999.02772395 + 195.04883150 $i$ | 2226 | 197 775.84328130 + 195.00361586 $i$ |  |  |
| 2207 | 196 087.68016811 + 195.03162853 $i$ | 2227 | 197 864.93992239 + 195.09614433 $i$ |  |  |
| 2208 | 196 176.66666742 + 194.81448694 $i$ | 2228 | 197 953.51541918 + 195.28820958 $i$ |  |  |
| 2209 | 196 265.64158213 + 195.11149056 $i$ | 2229 | 198 042.38919818 + 194.94931088 $i$ |  |  |
| 2210 | 196 354.16999420 + 195.04242885 $i$ | 2230 | 198 131.48215527 + 195.17730593 $i$ |  |  |
| 2211 | 196 443.25050710 + 194.83115517 $i$ | 2231 | 198 220.05848885 + 195.25749781 $i$ |  |  |
| 2212 | 196 532.15239199 + 195.14006415 $i$ | 2232 | 198 308.94906931 + 195.02059190 $i$ |  |  |
| 2213 | 196 620.74121989 + 195.04599935 $i$ | 2233 | 198 397.99274412 + 195.16504559 $i$ |  |  |
| 2214 | 196 709.78332907 + 194.89692775 $i$ | 2234 | 198 486.64888957 + 195.27089458 $i$ |  |  |
| 2215 | 196 798.68220395 + 195.09991027 $i$ | 2235 | 198 575.45711907 + 195.09444120 $i$ |  |  |
| 2216 | 196 887.34090024 + 195.11682825 $i$ | 2236 | 198 664.54416720 + 195.10542671 $i$ |  |  |
| 2217 | 196 976.24655758 + 194.89509059 $i$ | 2237 | 198 753.21966991 + 195.37060149 $i$ |  |  |
| 2218 | 197 065.29972481 + 195.10695360 $i$ | 2238 | 198 841.95999157 + 195.04974594 $i$ |  |  |
| 2219 | 197 153.83336325 + 195.17765618 $i$ | 2239 | 198 931.12149877 + 195.17420204 $i$ |  |  |
| 2220 | 197 242.82316268 + 194.87467035 $i$ | 2240 | 199 019.73271243 + 195.35719696 $i$ |  |  |

| k | $\sigma_k$ | k | $\sigma_k$ | k | $\sigma_k$ |
|---|---|---|---|---|---|
| 2241 | 199 108.54350470 + 195.09977603 $i$ | 2281 | 202 662.50713526 + 195.36399005 $i$ | 2321 | 206 216.53174500 + 195.59632511 $i$ |
| 2242 | 199 197.61989730 + 195.19157156 $i$ | 2282 | 202 751.56418816 + 195.50901789 $i$ | 2322 | 206 305.46690225 + 195.89508755 $i$ |
| 2243 | 199 286.32552765 + 195.34274130 $i$ | 2283 | 202 840.17993647 + 195.63116083 $i$ | 2323 | 206 394.00752752 + 195.80185232 $i$ |
| 2244 | 199 375.06444429 + 195.18183672 $i$ | 2284 | 202 929.02156880 + 195.35953632 $i$ | 2324 | 206 483.10665855 + 195.62116009 $i$ |
| 2245 | 199 464.14538454 + 195.14407145 $i$ | 2285 | 203 018.15169136 + 195.52979799 $i$ | 2325 | 206 571.96504631 + 195.91315889 $i$ |
| 2246 | 199 552.91837838 + 195.41222599 $i$ | 2286 | 203 106.67581111 + 195.66275743 $i$ | 2326 | 206 660.61483997 + 195.80716244 $i$ |
| 2247 | 199 641.55120988 + 195.18688673 $i$ | 2287 | 203 195.62162519 + 195.36994437 $i$ | 2327 | 206 749.58462780 + 195.68664073 $i$ |
| 2248 | 199 730.74591640 + 195.15903128 $i$ | 2288 | 203 284.62975878 + 195.57019451 $i$ | 2328 | 206 838.55212034 + 195.86121455 $i$ |
| 2249 | 199 819.40992189 + 195.44171500 $i$ | 2289 | 203 373.29206688 + 195.63945627 $i$ | 2329 | 206 927.15501655 + 195.90598444 $i$ |
| 2250 | 199 908.15092398 + 195.19337265 $i$ | 2290 | 203 462.11024818 + 195.45863819 $i$ | 2330 | 207 016.10070682 + 195.64919212 $i$ |
| 2251 | 199 997.23219119 + 195.21034364 $i$ | 2291 | 203 551.19487711 + 195.52059544 $i$ | 2331 | 207 105.13670159 + 195.90981623 $i$ |
| 2252 | 200 085.99659009 + 195.40912973 $i$ | 2292 | 203 639.85937015 + 195.70427224 $i$ | 2332 | 207 193.66235629 + 195.91980880 $i$ |
| 2253 | 200 174.68752528 + 195.27748436 $i$ | 2293 | 203 728.61229716 + 195.46543594 $i$ | 2333 | 207 282.68039970 + 195.67182079 $i$ |
| 2254 | 200 263.74338778 + 195.18569082 $i$ | 2294 | 203 817.77622369 + 195.52036028 $i$ | 2334 | 207 371.64876748 + 195.93249449 $i$ |
| 2255 | 200 352.59604852 + 195.43980814 $i$ | 2295 | 203 906.36946619 + 195.76076356 $i$ | 2335 | 207 460.24099984 + 195.91132124 $i$ |
| 2256 | 200 441.17913675 + 195.32769335 $i$ | 2296 | 203 995.19349418 + 195.44670230 $i$ | 2336 | 207 549.20318832 + 195.74547749 $i$ |
| 2257 | 200 530.32997974 + 195.15797953 $i$ | 2297 | 204 084.27592651 + 195.58787934 $i$ | 2337 | 207 638.17426444 + 195.87970501 $i$ |
| 2258 | 200 619.10598913 + 195.51271515 $i$ | 2298 | 204 172.94893704 + 195.72011634 $i$ | 2338 | 207 726.84228486 + 196.00307326 $i$ |
| 2259 | 200 707.76279192 + 195.29150647 $i$ | 2299 | 204 261.73156423 + 195.54420227 $i$ | 2339 | 207 815.67288514 + 195.72120972 $i$ |
| 2260 | 200 796.83239958 + 195.24276251 $i$ | 2300 | 204 350.78896717 + 195.54328245 $i$ | 2340 | 207 904.78770205 + 195.91561219 $i$ |
| 2261 | 200 885.67200454 + 195.46575982 $i$ | 2301 | 204 439.55458884 + 195.77548533 $i$ | 2341 | 207 993.34045629 + 196.02067372 $i$ |
| 2262 | 200 974.32302552 + 195.37000765 $i$ | 2302 | 204 528.21771170 + 195.56580025 $i$ | 2342 | 208 082.24411615 + 195.73897366 $i$ |
| 2263 | 201 063.33338875 + 195.23934751 $i$ | 2303 | 204 617.37252940 + 195.53620342 $i$ | 2343 | 208 171.31425343 + 195.94560723 $i$ |
| 2264 | 201 152.26693264 + 195.46084293 $i$ | 2304 | 204 706.08014613 + 195.82656047 $i$ | 2344 | 208 259.89299394 + 196.01141593 $i$ |
| 2265 | 201 240.82791689 + 195.44825847 $i$ | 2305 | 204 794.76932341 + 195.55515767 $i$ | 2345 | 208 348.80545642 + 195.80736107 $i$ |
| 2266 | 201 329.89435224 + 195.19757354 $i$ | 2306 | 204 883.91743653 + 195.59511605 $i$ | 2346 | 208 437.80207149 + 195.91226700 $i$ |
| 2267 | 201 418.81038085 + 195.54008783 $i$ | 2307 | 204 972.60667078 + 195.79152583 $i$ | 2347 | 208 526.52333226 + 196.06650811 $i$ |
| 2268 | 201 507.37251427 + 195.41003395 $i$ | 2308 | 205 061.36082494 + 195.63431054 $i$ | 2348 | 208 615.26003712 + 195.82774242 $i$ |
| 2269 | 201 596.44287054 + 195.28026840 $i$ | 2309 | 205 150.38331508 + 195.57530313 $i$ | 2349 | 208 704.41937180 + 195.90669187 $i$ |
| 2270 | 201 685.32990958 + 195.50102146 $i$ | 2310 | 205 239.23706522 + 195.81822191 $i$ | 2350 | 208 793.02445072 + 196.11267639 $i$ |
| 2271 | 201 773.97052127 + 195.46937860 $i$ | 2311 | 205 327.83991248 + 195.68108261 $i$ | 2351 | 208 881.82358175 + 195.82744650 $i$ |
| 2272 | 201 862.92124514 + 195.29802676 $i$ | 2312 | 205 416.96171444 + 195.55787942 $i$ | 2352 | 208 970.95640818 + 195.94640028 $i$ |
| 2273 | 201 951.92159417 + 195.48515260 $i$ | 2313 | 205 505.77369054 + 195.86511654 $i$ | 2353 | 209 059.56055661 + 196.10401404 $i$ |
| 2274 | 202 040.50091973 + 195.54966280 $i$ | 2314 | 205 594.37392257 + 195.68741988 $i$ | 2354 | 209 148.40452962 + 195.88431988 $i$ |
| 2275 | 202 129.45283163 + 195.26317495 $i$ | 2315 | 205 683.53182901 + 195.59532617 $i$ | 2355 | 209 237.43552689 + 195.93942712 $i$ |
| 2276 | 202 218.49177970 + 195.54379460 $i$ | 2316 | 205 772.27297484 + 195.86153925 $i$ | 2356 | 209 326.18483777 + 196.11606717 $i$ |
| 2277 | 202 307.01387188 + 195.54245005 $i$ | 2317 | 205 860.99463322 + 195.72324045 $i$ | 2357 | 209 414.88013642 + 195.94607857 $i$ |
| 2278 | 202 396.03952811 + 195.31083854 $i$ | 2318 | 205 949.98021489 + 195.62122742 $i$ | 2358 | 209 504.01224536 + 195.89880999 $i$ |
| 2279 | 202 484.97338878 + 195.54371184 $i$ | 2319 | 206 038.90452420 + 195.84430773 $i$ | 2359 | 209 592.72562475 + 196.19860056 $i$ |
| 2280 | 202 573.63987900 + 195.56144037 $i$ | 2320 | 206 127.48594411 + 195.79446284 $i$ | 2360 | 209 681.41780169 + 195.91899044 $i$ |

| k | $\sigma_k$ |
|---|---|
| 2361 | 209 770.56875347 + 195.95796932 $i$ |
| 2362 | 209 859.25461788 + 196.17881703 $i$ |
| 2363 | 209 947.99570705 + 195.96453051 $i$ |
| 2364 | 210 037.05852776 + 195.97608949 $i$ |
| 2365 | 210 125.84684249 + 196.16378279 $i$ |
| 2366 | 210 214.52235198 + 196.05531523 $i$ |
| 2367 | 210 303.58735397 + 195.91479708 $i$ |
| 2368 | 210 392.43068730 + 196.24699430 $i$ |
| 2369 | 210 481.01953631 + 196.03016854 $i$ |
| 2370 | 210 570.17120030 + 195.97382146 $i$ |
| 2371 | 210 658.94162137 + 196.23247208 $i$ |
| 2372 | 210 747.60220348 + 196.06845478 $i$ |
| 2373 | 210 836.67848669 + 196.00481042 $i$ |
| 2374 | 210 925.49944809 + 196.20506414 $i$ |
| 2375 | 211 014.17309792 + 196.15858021 $i$ |
| 2376 | 211 103.16256401 + 195.95795841 $i$ |
| 2377 | 211 192.11514605 + 196.26624965 $i$ |
| 2378 | 211 280.65457861 + 196.16286972 $i$ |
| 2379 | 211 369.75852474 + 195.98386485 $i$ |
| 2380 | 211 458.61910570 + 196.28020472 $i$ |
| 2381 | 211 547.23490398 + 196.16969318 $i$ |
| 2382 | 211 636.27249198 + 196.04128530 $i$ |
| 2383 | 211 725.16476250 + 196.24628826 $i$ |
| 2384 | 211 813.82074097 + 196.24000654 $i$ |
| 2385 | 211 902.74525530 + 196.03520297 $i$ |
| 2386 | 211 991.78351763 + 196.25895872 $i$ |
| 2387 | 212 080.30925549 + 196.28765727 $i$ |
| 2388 | 212 169.32291173 + 196.02075502 $i$ |
| 2389 | 212 258.30204849 + 196.31786402 $i$ |
| 2390 | 212 346.88014851 + 196.26207945 $i$ |
| 2391 | 212 435.85691256 + 196.09687341 $i$ |
| 2392 | 212 524.82896915 + 196.27513901 $i$ |
| 2393 | 212 613.46808662 + 196.32136807 $i$ |
| 2394 | 212 702.34827536 + 196.12049307 $i$ |
| 2395 | 212 791.41319758 + 196.24772458 $i$ |
| 2396 | 212 879.99639292 + 196.40833036 $i$ |
| 2397 | 212 968.89157622 + 196.07283399 $i$ |
| 2398 | 213 057.96165015 + 196.32898890 $i$ |
| 2399 | 213 146.53804626 + 196.36648640 $i$ |
| 2400 | 213 235.45191226 + 196.15384067 $i$ |

| $k$ | $\sigma_k$ | $k$ | $\sigma_k$ | $k$ | $\sigma_k$ |
|---|---|---|---|---|---|
| 2401 | 213 324.46942582 + 196.29511940 $i$ | 2441 | 216 878.39963662 + 196.61141060 $i$ | 2481 | 220 432.23528646 + 196.96113949 $i$ |
| 2402 | 213 413.12848356 + 196.40268103 $i$ | 2442 | 216 966.98111198 + 196.64149287 $i$ | 2482 | 220 520.88111426 + 196.75601224 $i$ |
| 2403 | 213 501.95483669 + 196.20550370 $i$ | 2443 | 217 055.95829769 + 196.37290058 $i$ | 2483 | 220 610.00881108 + 196.69789503 $i$ |
| 2404 | 213 591.02858275 + 196.25454612 $i$ | 2444 | 217 144.96172450 + 196.66764089 $i$ | 2484 | 220 698.75578860 + 196.95924130 $i$ |
| 2405 | 213 679.69489586 + 196.48733504 $i$ | 2445 | 217 233.50238271 + 196.62934029 $i$ | 2485 | 220 787.48503706 + 196.79149953 $i$ |
| 2406 | 213 768.45745542 + 196.16733701 $i$ | 2446 | 217 322.53455495 + 196.42569999 $i$ | 2486 | 220 876.46840817 + 196.72360563 $i$ |
| 2407 | 213 857.61326428 + 196.32337371 $i$ | 2447 | 217 411.45495863 + 196.65940064 $i$ | 2487 | 220 965.38713671 + 196.93762885 $i$ |
| 2408 | 213 946.20395642 + 196.46195621 $i$ | 2448 | 217 500.12002683 + 196.65110985 $i$ | 2488 | 221 053.96618791 + 196.86901325 $i$ |
| 2409 | 214 035.04521947 + 196.22614133 $i$ | 2449 | 217 589.00046493 + 196.47779537 $i$ | 2489 | 221 143.03887682 + 196.68857111 $i$ |
| 2410 | 214 124.10352674 + 196.32109082 $i$ | 2450 | 217 678.05723856 + 196.62098433 $i$ | 2490 | 221 231.92978713 + 196.99410790 $i$ |
| 2411 | 214 212.80032645 + 196.46493878 $i$ | 2451 | 217 766.64458425 + 196.72667230 $i$ | 2491 | 221 320.50424143 + 196.86894322 $i$ |
| 2412 | 214 301.55790526 + 196.29830569 $i$ | 2452 | 217 855.53359378 + 196.45978502 $i$ | 2492 | 221 409.60080642 + 196.71873061 $i$ |
| 2413 | 214 390.64208153 + 196.27679653 $i$ | 2453 | 217 944.62421526 + 196.65242010 $i$ | 2493 | 221 498.43832859 + 197.00205547 $i$ |
| 2414 | 214 479.38353213 + 196.53543952 $i$ | 2454 | 218 033.16110129 + 196.74953523 $i$ | 2494 | 221 587.10364240 + 196.88367110 $i$ |
| 2415 | 214 568.05344382 + 196.29291159 $i$ | 2455 | 218 122.11889117 + 196.47048059 $i$ | 2495 | 221 676.08376090 + 196.77450469 $i$ |
| 2416 | 214 657.23312587 + 196.29927006 $i$ | 2456 | 218 211.10997130 + 196.68859428 $i$ | 2496 | 221 765.02951056 + 196.95010169 $i$ |
| 2417 | 214 745.88228050 + 196.55915452 $i$ | 2457 | 218 299.77537846 + 196.72742600 $i$ | 2497 | 221 853.63380384 + 196.98238563 $i$ |
| 2418 | 214 834.65184138 + 196.30049121 $i$ | 2458 | 218 388.60513556 + 196.55701950 $i$ | 2498 | 221 942.61037672 + 196.72976884 $i$ |
| 2419 | 214 923.71535946 + 196.34615805 $i$ | 2459 | 218 477.67880105 + 196.62734659 $i$ | 2499 | 222 031.60839049 + 197.00438936 $i$ |
| 2420 | 215 012.47708049 + 196.52612148 $i$ | 2460 | 218 566.32898207 + 196.81186479 $i$ | 2500 | 222 120.14386929 + 196.98139066 $i$ |
| 2421 | 215 101.17713809 + 196.38838313 $i$ | 2461 | 218 655.12280416 + 196.54576981 $i$ | 2501 | 222 209.18378323 + 196.76060227 $i$ |
| 2422 | 215 190.23987893 + 196.30798115 $i$ | 2462 | 218 744.25037335 + 196.64315719 $i$ | 2502 | 222 298.11687218 + 197.02141911 $i$ |
| 2423 | 215 279.06514619 + 196.57168928 $i$ | 2463 | 218 832.84918910 + 196.85116561 $i$ | 2503 | 222 386.73365273 + 196.98322216 $i$ |
| 2424 | 215 367.67820899 + 196.41843301 $i$ | 2464 | 218 921.69440943 + 196.53711910 $i$ | 2504 | 222 475.69429069 + 196.82163141 $i$ |
| 2425 | 215 456.81719801 + 196.29433398 $i$ | 2465 | 219 010.75718953 + 196.70424512 $i$ | 2505 | 222 564.65968442 + 196.97622627 $i$ |
| 2426 | 215 545.58375891 + 196.63257069 $i$ | 2466 | 219 099.43036573 + 196.80866403 $i$ | 2506 | 222 653.32121023 + 197.06441954 $i$ |
| 2427 | 215 634.25411087 + 196.38606497 $i$ | 2467 | 219 188.22661243 + 196.63700622 $i$ | 2507 | 222 742.16735810 + 196.79975769 $i$ |
| 2428 | 215 723.32573061 + 196.37927664 $i$ | 2468 | 219 277.28187459 + 196.64548615 $i$ | 2508 | 222 831.27631478 + 197.01209689 $i$ |
| 2429 | 215 812.14910555 + 196.57169421 $i$ | 2469 | 219 366.02320008 + 196.87125639 $i$ | 2509 | 222 919.80966614 + 197.07940807 $i$ |
| 2430 | 215 900.80711722 + 196.48299519 $i$ | 2470 | 219 454.71283694 + 196.64901875 $i$ | 2510 | 223 008.75476000 + 196.81950762 $i$ |
| 2431 | 215 989.83522467 + 196.35139498 $i$ | 2471 | 219 543.87090816 + 196.64915608 $i$ | 2511 | 223 097.78225568 + 197.03240688 $i$ |
| 2432 | 216 078.73653986 + 196.59050947 $i$ | 2472 | 219 632.54288624 + 196.91160511 $i$ | 2512 | 223 186.38034759 + 197.08234487 $i$ |
| 2433 | 216 167.31957180 + 196.54062423 $i$ | 2473 | 219 721.27427485 + 196.64247666 $i$ | 2513 | 223 275.29989478 + 196.87956946 $i$ |
| 2434 | 216 256.39419315 + 196.31845497 $i$ | 2474 | 219 810.40179729 + 196.69993161 $i$ | 2514 | 223 364.29039104 + 197.00104144 $i$ |
| 2435 | 216 345.27959234 + 196.66018486 $i$ | 2475 | 219 899.08269535 + 196.88216855 $i$ | 2515 | 223 452.99683074 + 197.13171297 $i$ |
| 2436 | 216 433.86344826 + 196.50337662 $i$ | 2476 | 219 987.85381506 + 196.71814002 $i$ | 2516 | 223 541.75775054 + 196.89896801 $i$ |
| 2437 | 216 522.93982264 + 196.40233909 $i$ | 2477 | 220 076.87627011 + 196.67767857 $i$ | 2517 | 223 630.91053560 + 196.99368143 $i$ |
| 2438 | 216 611.80221615 + 196.61374169 $i$ | 2478 | 220 165.71288672 + 196.91126745 $i$ | 2518 | 223 719.48800770 + 197.18113813 $i$ |
| 2439 | 216 700.46038665 + 196.57421532 $i$ | 2479 | 220 254.32666324 + 196.75857680 $i$ | 2519 | 223 808.33890713 + 196.89252563 $i$ |
| 2440 | 216 789.41722825 + 196.40228462 $i$ | 2480 | 220 343.46503820 + 196.65977453 $i$ | 2520 | 223 897.42744126 + 197.03309454 $i$ |
| | | | | 2521 | 223 986.04528421 + 197.17150818 $i$ |
| | | | | 2522 | 224 074.90165664 + 196.94588303 $i$ |
| | | | | 2523 | 224 163.91765967 + 197.02129808 $i$ |
| | | | | 2524 | 224 252.66247567 + 197.18926188 $i$ |
| | | | | 2525 | 224 341.37460464 + 197.00150562 $i$ |
| | | | | 2526 | 224 430.50470540 + 196.98663935 $i$ |
| | | | | 2527 | 224 519.19581246 + 197.26484816 $i$ |
| | | | | 2528 | 224 607.91861549 + 196.97240459 $i$ |
| | | | | 2529 | 224 697.05553671 + 197.04599657 $i$ |
| | | | | 2530 | 224 785.72507975 + 197.23759630 $i$ |
| | | | | 2531 | 224 874.49496651 + 197.03229217 $i$ |
| | | | | 2532 | 224 963.54792080 + 197.04652720 $i$ |
| | | | | 2533 | 225 052.32136637 + 197.23616390 $i$ |
| | | | | 2534 | 225 141.01022348 + 197.10718163 $i$ |
| | | | | 2535 | 225 230.08615471 + 196.99869658 $i$ |
| | | | | 2536 | 225 318.90220643 + 197.31303714 $i$ |
| | | | | 2537 | 225 407.51269371 + 197.08092312 $i$ |
| | | | | 2538 | 225 496.67257269 + 197.05157136 $i$ |
| | | | | 2539 | 225 585.40069496 + 197.29224739 $i$ |
| | | | | 2540 | 225 674.10652129 + 197.13021364 $i$ |
| | | | | 2541 | 225 763.16408325 + 197.06803914 $i$ |
| | | | | 2542 | 225 851.98055133 + 197.28151205 $i$ |
| | | | | 2543 | 225 940.65955357 + 197.19860232 $i$ |
| | | | | 2544 | 226 029.65659396 + 197.03248304 $i$ |
| | | | | 2545 | 226 118.59302223 + 197.33803418 $i$ |
| | | | | 2546 | 226 207.13933634 + 197.20186290 $i$ |
| | | | | 2547 | 226 296.25942690 + 197.05904979 $i$ |
| | | | | 2548 | 226 385.08587335 + 197.34598419 $i$ |
| | | | | 2549 | 226 473.73223902 + 197.21649630 $i$ |
| | | | | 2550 | 226 562.76050165 + 197.10758446 $i$ |
| | | | | 2551 | 226 651.64816597 + 197.31412653 $i$ |
| | | | | 2552 | 226 740.30131188 + 197.28557617 $i$ |
| | | | | 2553 | 226 829.24566331 + 197.09540405 $i$ |
| | | | | 2554 | 226 918.25914775 + 197.33231516 $i$ |
| | | | | 2555 | 227 006.79163556 + 197.32862606 $i$ |
| | | | | 2556 | 227 095.83157895 + 197.08243903 $i$ |
| | | | | 2557 | 227 184.76986576 + 197.38006842 $i$ |
| | | | | 2558 | 227 273.36569620 + 197.30705573 $i$ |
| | | | | 2559 | 227 362.35977645 + 197.15697340 $i$ |
| | | | | 2560 | 227 451.30134835 + 197.33006156 $i$ |

| k | $\sigma_k$ | k | $\sigma_k$ | k | $\sigma_k$ |
|---|---|---|---|---|---|
| 2561 | 227 539.94902776 + 197.37974436 $i$ | 2601 | 231 093.81174138 + 197.55397004 $i$ | 2641 | 234 647.78281278 + 197.65588882 $i$ |
| 2562 | 227 628.85177319 + 197.16228806 $i$ | 2602 | 231 182.88943433 + 197.37004726 $i$ | 2642 | 234 736.88047094 + 197.73267986 $i$ |
| 2563 | 227 717.89083295 + 197.32284288 $i$ | 2603 | 231 271.75069089 + 197.70452540 $i$ | 2643 | 234 825.56288890 + 197.90377610 $i$ |
| 2564 | 227 806.47495488 + 197.44322582 $i$ | 2604 | 231 360.35635958 + 197.52437797 $i$ | 2644 | 234 914.34884030 + 197.72765642 $i$ |
| 2565 | 227 895.39168213 + 197.12809238 $i$ | 2605 | 231 449.43507125 + 197.44897813 $i$ | 2645 | 235 003.36641111 + 197.70520311 $i$ |
| 2566 | 227 984.44283810 + 197.39701175 $i$ | 2606 | 231 538.27582318 + 197.65324623 $i$ | 2646 | 235 092.18455529 + 197.93585926 $i$ |
| 2567 | 228 073.01202956 + 197.40085989 $i$ | 2607 | 231 626.95129041 + 197.60224424 $i$ | 2647 | 235 180.82036359 + 197.76492045 $i$ |
| 2568 | 228 161.95409345 + 197.21698462 $i$ | 2608 | 231 715.91117013 + 197.43441188 $i$ | 2648 | 235 269.96246096 + 197.68804328 $i$ |
| 2569 | 228 250.95132659 + 197.35169411 $i$ | 2609 | 231 804.87979320 + 197.66266004 $i$ | 2649 | 235 358.69757761 + 197.98705080 $i$ |
| 2570 | 228 339.61223494 + 197.45079220 $i$ | 2610 | 231 893.46104505 + 197.65556240 $i$ | 2650 | 235 447.38833867 + 197.75610362 $i$ |
| 2571 | 228 428.44624895 + 197.24148769 $i$ | 2611 | 231 982.46238240 + 197.41243801 $i$ | 2651 | 235 536.49156118 + 197.73120362 $i$ |
| 2572 | 228 517.51770830 + 197.33000953 $i$ | 2612 | 232 071.43190398 + 197.71314623 $i$ | 2652 | 235 625.23751970 + 197.97265267 $i$ |
| 2573 | 228 606.16923675 + 197.52096370 $i$ | 2613 | 232 159.98952374 + 197.64419623 $i$ | 2653 | 235 713.96701372 + 197.79893584 $i$ |
| 2574 | 228 694.95751381 + 197.21193390 $i$ | 2614 | 232 249.03227575 + 197.46879920 $i$ | 2654 | 235 802.96961011 + 197.75328551 $i$ |
| 2575 | 228 784.09961024 + 197.38801179 $i$ | 2615 | 232 337.93444896 + 197.69426622 $i$ | 2655 | 235 891.86156143 + 197.95133160 $i$ |
| 2576 | 228 872.67307665 + 197.49738589 $i$ | 2616 | 232 426.59794522 + 197.67440428 $i$ | 2656 | 235 980.44912281 + 197.87820245 $i$ |
| 2577 | 228 961.55189369 + 197.27500941 $i$ | 2617 | 232 515.50171301 + 197.51304786 $i$ | 2657 | 236 069.54826699 + 197.70629707 $i$ |
| 2578 | 229 050.58026447 + 197.37311650 $i$ | 2618 | 232 604.53927738 + 197.65720721 $i$ | 2658 | 236 158.39057741 + 198.01686377 $i$ |
| 2579 | 229 139.28009675 + 197.51758251 $i$ | 2619 | 232 693.11588658 + 197.75473368 $i$ | 2659 | 236 246.99976882 + 197.86571606 $i$ |
| 2580 | 229 228.05608124 + 197.33134862 $i$ | 2620 | 232 782.04682094 + 197.48060725 $i$ | 2660 | 236 336.09591636 + 197.74806499 $i$ |
| 2581 | 229 317.13053929 + 197.33310860 $i$ | 2621 | 232 871.09088564 + 197.70185929 $i$ | 2661 | 236 424.91564065 + 198.01243184 $i$ |
| 2582 | 229 405.85053387 + 197.58571540 $i$ | 2622 | 232 959.64850762 + 197.76356515 $i$ | 2662 | 236 513.58824156 + 197.88636202 $i$ |
| 2583 | 229 494.55901249 + 197.32188356 $i$ | 2623 | 233 048.61674359 + 197.50112971 $i$ | 2663 | 236 602.57908484 + 197.79221936 $i$ |
| 2584 | 229 583.71651501 + 197.36290210 $i$ | 2624 | 233 137.59015760 + 197.72682906 $i$ | 2664 | 236 691.50741942 + 197.97150106 $i$ |
| 2585 | 229 672.35542783 + 197.59914673 $i$ | 2625 | 233 226.25307843 + 197.74710184 $i$ | 2665 | 236 780.11627030 + 197.98217400 $i$ |
| 2586 | 229 761.15231972 + 197.33621939 $i$ | 2626 | 233 315.10530163 + 197.58532906 $i$ | 2666 | 236 869.11118524 + 197.74075154 $i$ |
| 2587 | 229 850.20084464 + 197.40289223 $i$ | 2627 | 233 404.16354981 + 197.66026865 $i$ | 2667 | 236 958.07924665 + 198.03526486 $i$ |
| 2588 | 229 938.95426623 + 197.56685471 $i$ | 2628 | 233 492.80129077 + 197.84055641 $i$ | 2668 | 237 046.63307711 + 197.97182028 $i$ |
| 2589 | 230 027.66898471 + 197.42601604 $i$ | 2629 | 233 581.62838718 + 197.55562046 $i$ | 2669 | 237 135.68357437 + 197.77897811 $i$ |
| 2590 | 230 116.73779548 + 197.35414714 $i$ | 2630 | 233 670.73054203 + 197.69438139 $i$ | 2670 | 237 224.58846493 + 198.03951563 $i$ |
| 2591 | 230 205.53366498 + 197.62228394 $i$ | 2631 | 233 759.32764034 + 197.85725062 $i$ | 2671 | 237 313.22554197 + 197.98070652 $i$ |
| 2592 | 230 294.17374818 + 197.43700934 $i$ | 2632 | 233 848.19120660 + 197.56338946 $i$ | 2672 | 237 402.18263384 + 197.83052480 $i$ |
| 2593 | 230 383.31026930 + 197.35760950 $i$ | 2633 | 233 937.24558859 + 197.74092747 $i$ | 2673 | 237 491.14513495 + 198.00018608 $i$ |
| 2594 | 230 472.05727205 + 197.66554864 $i$ | 2634 | 234 025.90525545 + 197.81929256 $i$ | 2674 | 237 579.79699796 + 198.06082606 $i$ |
| 2595 | 230 560.74314923 + 197.41709722 $i$ | 2635 | 234 114.71947346 + 197.66076546 $i$ | 2675 | 237 668.67405929 + 197.80864193 $i$ |
| 2596 | 230 649.82543893 + 197.43366245 $i$ | 2636 | 234 203.77398716 + 197.67641776 $i$ | 2676 | 237 757.75538735 + 198.02713294 $i$ |
| 2597 | 230 738.61935710 + 197.60218010 $i$ | 2637 | 234 292.49386950 + 197.89518184 $i$ | 2677 | 237 846.27848282 + 198.07743495 $i$ |
| 2598 | 230 827.29620332 + 197.52169325 $i$ | 2638 | 234 381.21154888 + 197.65594428 $i$ | 2678 | 237 935.27144491 + 197.82903259 $i$ |
| 2599 | 230 916.33203956 + 197.38644629 $i$ | 2639 | 234 470.35957394 + 197.68890590 $i$ | 2679 | 238 024.24924942 + 198.04373449 $i$ |
| 2600 | 231 005.20730122 + 197.64793221 $i$ | 2640 | 234 559.00843654 + 197.93003546 $i$ | 2680 | 238 112.86656909 + 198.08099714 $i$ |
| 2681 | 238 201.79288005 + 197.88416803 $i$ | | | | |
| 2682 | 238 290.77922230 + 198.01596800 $i$ | | | | |
| 2683 | 238 379.46610943 + 198.12947895 $i$ | | | | |
| 2684 | 238 468.26250301 + 197.90114509 $i$ | | | | |
| 2685 | 238 557.39614465 + 198.00728225 $i$ | | | | |
| 2686 | 238 645.95496491 + 198.17915317 $i$ | | | | |
| 2687 | 238 734.84893793 + 197.88773802 $i$ | | | | |
| 2688 | 238 823.89965670 + 198.05508921 $i$ | | | | |
| 2689 | 238 912.53290534 + 198.16248325 $i$ | | | | |
| 2690 | 239 001.39036021 + 197.94346218 $i$ | | | | |
| 2691 | 239 090.40690596 + 198.03784814 $i$ | | | | |
| 2692 | 239 179.13973756 + 198.18757207 $i$ | | | | |
| 2693 | 239 267.86986054 + 197.99079030 $i$ | | | | |
| 2694 | 239 356.99867047 + 198.00226575 $i$ | | | | |
| 2695 | 239 445.66154921 + 198.25754724 $i$ | | | | |
| 2696 | 239 534.42011913 + 197.96362761 $i$ | | | | |
| 2697 | 239 623.54279443 + 198.06213399 $i$ | | | | |
| 2698 | 239 712.19651745 + 198.22725375 $i$ | | | | |
| 2699 | 239 800.99457457 + 198.02971846 $i$ | | | | |
| 2700 | 239 890.03414564 + 198.04877787 $i$ | | | | |
| 2701 | 239 978.79853482 + 198.24132093 $i$ | | | | |
| 2702 | 240 067.50260696 + 198.08793245 $i$ | | | | |
| 2703 | 240 156.58364211 + 198.00828794 $i$ | | | | |
| 2704 | 240 245.36685506 + 198.30866593 $i$ | | | | |
| 2705 | 240 334.01029183 + 198.06657008 $i$ | | | | |
| 2706 | 240 423.16707635 + 198.05872545 $i$ | | | | |
| 2707 | 240 511.86513942 + 198.28880850 $i$ | | | | |
| 2708 | 240 600.61212454 + 198.11639022 $i$ | | | | |
| 2709 | 240 689.64384225 + 198.06655125 $i$ | | | | |
| 2710 | 240 778.46390375 + 198.28755167 $i$ | | | | |
| 2711 | 240 867.14130096 + 198.17516581 $i$ | | | | |
| 2712 | 240 956.16079582 + 198.04213942 $i$ | | | | |
| 2713 | 241 045.06600947 + 198.33080329 $i$ | | | | |
| 2714 | 241 133.62378204 + 198.17981770 $i$ | | | | |
| 2715 | 241 222.76317539 + 198.06462280 $i$ | | | | |
| 2716 | 241 311.55127295 + 198.33857626 $i$ | | | | |
| 2717 | 241 400.22752382 + 198.19816145 $i$ | | | | |
| 2718 | 241 489.25155635 + 198.10674466 $i$ | | | | |
| 2719 | 241 578.13027510 + 198.30834011 $i$ | | | | |
| 2720 | 241 666.77922216 + 198.26690335 $i$ | | | | |

| $k$ | $\sigma_k$ | $k$ | $\sigma_k$ | $k$ | $\sigma_k$ |
|---|---|---|---|---|---|
| 2721 | 241 755.75234585 + 198.08805879 $i$ | 2761 | 245 309.80329538 + 198.35427319 $i$ | 2801 | 248 863.73187820 + 198.70688512 $i$ |
| 2722 | 241 844.73086695 + 198.33060864 $i$ | 2762 | 245 398.52682262 + 198.63155977 $i$ | 2802 | 248 952.37772373 + 198.77444361 $i$ |
| 2723 | 241 933.27222589 + 198.30259638 $i$ | 2763 | 245 487.23785431 + 198.38800979 $i$ | 2803 | 249 041.22054755 + 198.61933717 $i$ |
| 2724 | 242 022.33704093 + 198.07990253 $i$ | 2764 | 245 576.31988757 + 198.41451110 $i$ | 2804 | 249 130.26120395 + 198.63976725 $i$ |
| 2725 | 242 111.23825036 + 198.37567216 $i$ | 2765 | 245 665.08700182 + 198.57829442 $i$ | 2805 | 249 218.96312440 + 198.85769789 $i$ |
| 2726 | 242 199.85478317 + 198.28415912 $i$ | 2766 | 245 753.79536688 + 198.49150545 $i$ | 2806 | 249 307.71292860 + 198.60128574 $i$ |
| 2727 | 242 288.85691058 + 198.14969688 $i$ | 2767 | 245 842.82404681 + 198.35467229 $i$ | 2807 | 249 396.84760774 + 198.66589877 $i$ |
| 2728 | 242 377.77829395 + 198.32400968 $i$ | 2768 | 245 931.67990617 + 198.63847395 $i$ | 2808 | 249 485.47780954 + 198.88042629 $i$ |
| 2729 | 242 466.43393511 + 198.36148599 $i$ | 2769 | 246 020.30154946 + 198.50144653 $i$ | 2809 | 249 574.28751824 + 198.60621656 $i$ |
| 2730 | 242 555.34706071 + 198.14043953 $i$ | 2770 | 246 109.38776161 + 198.36074176 $i$ | 2810 | 249 663.36003679 + 198.70346984 $i$ |
| 2731 | 242 644.37356172 + 198.33326467 $i$ | 2771 | 246 198.22158730 + 198.67457150 $i$ | 2811 | 249 752.04356548 + 198.85670123 $i$ |
| 2732 | 242 732.95259789 + 198.40552624 $i$ | 2772 | 246 286.84725878 + 198.48478365 $i$ | 2812 | 249 840.83859409 + 198.67818120 $i$ |
| 2733 | 242 821.89271771 + 198.12156526 $i$ | 2773 | 246 375.93317945 + 198.42892404 $i$ | 2813 | 249 929.86229717 + 198.67218290 $i$ |
| 2734 | 242 910.92268970 + 198.39236363 $i$ | 2774 | 246 464.74829492 + 198.62179692 $i$ | 2814 | 250 018.65592968 + 198.89076343 $i$ |
| 2735 | 242 999.48942571 + 198.37205408 $i$ | 2775 | 246 553.43554148 + 198.56629555 $i$ | 2815 | 250 107.31140545 + 198.70868730 $i$ |
| 2736 | 243 088.45810052 + 198.20572165 $i$ | 2776 | 246 642.40774558 + 198.40876233 $i$ | 2816 | 250 196.45964739 + 198.65891251 $i$ |
| 2737 | 243 177.42582601 + 198.33979576 $i$ | 2777 | 246 731.36136837 + 198.64300907 $i$ | 2817 | 250 285.16579749 + 198.94155162 $i$ |
| 2738 | 243 266.09540407 + 198.43551372 $i$ | 2778 | 246 819.93581876 + 198.60068673 $i$ | 2818 | 250 373.88652252 + 198.69351801 $i$ |
| 2739 | 243 354.94208363 + 198.21300699 $i$ | 2779 | 246 908.97061775 + 198.39483009 $i$ | 2819 | 250 462.97669932 + 198.70850060 $i$ |
| 2740 | 243 444.00546837 + 198.33122380 $i$ | 2780 | 246 997.90423324 + 198.68648236 $i$ | 2820 | 250 551.72129733 + 198.91790341 $i$ |
| 2741 | 243 532.63514282 + 198.49155755 $i$ | 2781 | 247 086.47470468 + 198.59556270 $i$ | 2821 | 250 640.45036113 + 198.74737207 $i$ |
| 2742 | 243 621.46732723 + 198.19723950 $i$ | 2782 | 247 175.53033128 + 198.44869187 $i$ | 2822 | 250 729.47232853 + 198.71236972 $i$ |
| 2743 | 243 710.58176082 + 198.37694247 $i$ | 2783 | 247 264.41335507 + 198.66268189 $i$ | 2823 | 250 818.32856488 + 198.90541411 $i$ |
| 2744 | 243 799.14260680 + 198.47233872 $i$ | 2784 | 247 353.07896226 + 198.63479723 $i$ | 2824 | 250 906.93977985 + 198.82329509 $i$ |
| 2745 | 243 888.06191682 + 198.25574772 $i$ | 2785 | 247 442.00256556 + 198.47793565 $i$ | 2825 | 250 996.04944625 + 198.66078842 $i$ |
| 2746 | 243 977.05634381 + 198.35690751 $i$ | 2786 | 247 531.01628903 + 198.63346300 $i$ | 2826 | 251 084.85518818 + 198.98040685 $i$ |
| 2747 | 244 065.75947548 + 198.49931043 $i$ | 2787 | 247 619.59503578 + 198.71426564 $i$ | 2827 | 251 173.49796988 + 198.79500025 $i$ |
| 2748 | 244 154.54838610 + 198.30140940 $i$ | 2788 | 247 708.55307999 + 198.43643584 $i$ | 2828 | 251 262.58377413 + 198.71509463 $i$ |
| 2749 | 244 243.62461521 + 198.32878945 $i$ | 2789 | 247 797.56027617 + 198.68914049 $i$ | 2829 | 251 351.39080812 + 198.96370765 $i$ |
| 2750 | 244 332.31963176 + 198.55915848 $i$ | 2790 | 247 886.13594448 + 198.70706491 $i$ | 2830 | 251 440.07682315 + 198.83291407 $i$ |
| 2751 | 244 421.05798615 + 198.28863937 $i$ | 2791 | 247 975.11066710 + 198.47344219 $i$ | 2831 | 251 529.07870969 + 198.74456465 $i$ |
| 2752 | 244 510.20392405 + 198.36423825 $i$ | 2792 | 248 064.07319433 + 198.69845608 $i$ | 2832 | 251 617.98341875 + 198.92823840 $i$ |
| 2753 | 244 598.83094976 + 198.56414049 $i$ | 2793 | 248 152.73131073 + 198.70368606 $i$ | 2833 | 251 706.59858456 + 198.91857607 $i$ |
| 2754 | 244 687.64691474 + 198.30987003 $i$ | 2794 | 248 241.60723713 + 198.54627854 $i$ | 2834 | 251 795.61260133 + 198.69336049 $i$ |
| 2755 | 244 776.68986734 + 198.39455667 $i$ | 2795 | 248 330.64521740 + 198.62929781 $i$ | 2835 | 251 884.55143632 + 198.99765440 $i$ |
| 2756 | 244 865.43323302 + 198.53816059 $i$ | 2796 | 248 419.27375587 + 198.80114103 $i$ | 2836 | 251 973.11866494 + 198.90004451 $i$ |
| 2757 | 244 954.15930624 + 198.39407992 $i$ | 2797 | 248 508.13048424 + 198.50727150 $i$ | 2837 | 252 062.18431808 + 198.74040081 $i$ |
| 2758 | 245 043.23152794 + 198.33722654 $i$ | 2798 | 248 597.21435354 + 198.68043045 $i$ | 2838 | 252 151.06180685 + 198.98779453 $i$ |
| 2759 | 245 132.00333049 + 198.60594512 $i$ | 2799 | 248 685.80228454 + 198.79662582 $i$ | 2839 | 252 239.70930755 + 198.92008506 $i$ |
| 2760 | 245 220.66983117 + 198.39161761 $i$ | 2800 | 248 774.69054043 + 198.53279601 $i$ | 2840 | 252 328.67902498 + 198.78524351 $i$ |
| 2841 | 252 417.63370015 + 198.95085783 $i$ | | | | |
| 2842 | 252 506.26431217 + 198.99527313 $i$ | | | | |
| 2843 | 252 595.18284231 + 198.76069281 $i$ | | | | |
| 2844 | 252 684.23229446 + 198.98101572 $i$ | | | | |
| 2845 | 252 772.75478344 + 199.01139052 $i$ | | | | |
| 2846 | 252 861.77937248 + 198.77511728 $i$ | | | | |
| 2847 | 252 950.71886528 + 199.00056314 $i$ | | | | |
| 2848 | 253 039.35822729 + 199.01205550 $i$ | | | | |
| 2849 | 253 128.28080741 + 198.82763787 $i$ | | | | |
| 2850 | 253 217.26963993 + 198.97052467 $i$ | | | | |
| 2851 | 253 305.93745313 + 199.06316190 $i$ | | | | |
| 2852 | 253 394.76485931 + 198.83781381 $i$ | | | | |
| 2853 | 253 483.87419164 + 198.96630073 $i$ | | | | |
| 2854 | 253 572.43229452 + 199.11542816 $i$ | | | | |
| 2855 | 253 661.35491567 + 198.81941141 $i$ | | | | |
| 2856 | 253 750.37243658 + 199.01634781 $i$ | | | | |
| 2857 | 253 839.01864009 + 199.09025775 $i$ | | | | |
| 2858 | 253 927.88159629 + 198.88636427 $i$ | | | | |
| 2859 | 254 016.90024965 + 198.98435843 $i$ | | | | |
| 2860 | 254 105.60718945 + 199.12401255 $i$ | | | | |
| 2861 | 254 194.36935109 + 198.92545268 $i$ | | | | |
| 2862 | 254 283.49121976 + 198.95372308 $i$ | | | | |
| 2863 | 254 372.12542553 + 199.19068062 $i$ | | | | |
| 2864 | 254 460.92658129 + 198.89780298 $i$ | | | | |
| 2865 | 254 550.02866335 + 199.01241678 $i$ | | | | |
| 2866 | 254 638.66760904 + 199.15598547 $i$ | | | | |
| 2867 | 254 727.49268788 + 198.96764809 $i$ | | | | |
| 2868 | 254 816.52257371 + 198.99237398 $i$ | | | | |
| 2869 | 254 905.27697146 + 199.17955951 $i$ | | | | |
| 2870 | 254 993.98928077 + 199.00872322 $i$ | | | | |
| 2871 | 255 083.08100610 + 198.96336274 $i$ | | | | |
| 2872 | 255 171.83445199 + 199.24141471 $i$ | | | | |
| 2873 | 255 260.50858448 + 198.99077346 $i$ | | | | |
| 2874 | 255 349.65674224 + 199.00904636 $i$ | | | | |
| 2875 | 255 438.33801939 + 199.22299218 $i$ | | | | |
| 2876 | 255 527.11049575 + 199.03880414 $i$ | | | | |
| 2877 | 255 616.12746998 + 199.01377333 $i$ | | | | |
| 2878 | 255 704.95266962 + 199.22258846 $i$ | | | | |
| 2879 | 255 793.61721123 + 199.09137881 $i$ | | | | |
| 2880 | 255 882.66693672 + 198.99065851 $i$ | | | | |

| $k$ | $\sigma_k$ | $k$ | $\sigma_k$ | $k$ | $\sigma_k$ |
|---|---|---|---|---|---|
| 2881 | 255 971.53144166 + 199.26184091 $i$ | 2921 | 259 525.30406313 + 199.47173949 $i$ | 2961 | 263 079.20345702 + 199.60273844 $i$ | 
| 2882 | 256 060.11427531 + 199.10255564 $i$ | 2922 | 259 614.14545023 + 199.22882576 $i$ | 2962 | 263 168.10621790 + 199.45457112 $i$ |
| 2883 | 256 149.26306054 + 199.00719094 $i$ | 2923 | 259 703.17732513 + 199.32070024 $i$ | 2963 | 263 257.13023751 + 199.54766775 $i$ |
| 2884 | 256 238.01766143 + 199.27344869 $i$ | 2924 | 259 791.90633721 + 199.45490984 $i$ | 2964 | 263 345.75080608 + 199.69750414 $i$ |
| 2885 | 256 326.72530762 + 199.11898997 $i$ | 2925 | 259 880.65334169 + 199.30769191 $i$ | 2965 | 263 434.62877443 + 199.40502999 $i$ |
| 2886 | 256 415.74069348 + 199.04452719 $i$ | 2926 | 259 969.72632656 + 199.26203760 $i$ | 2966 | 263 523.69873097 + 199.60805887 $i$ |
| 2887 | 256 504.60979728 + 199.24270661 $i$ | 2927 | 260 058.47133066 + 199.52803572 $i$ | 2967 | 263 612.27606392 + 199.67928035 $i$ |
| 2888 | 256 593.25798713 + 199.19065366 $i$ | 2928 | 260 147.16641207 + 199.29341886 $i$ | 2968 | 263 701.19150915 + 199.44665094 $i$ |
| 2889 | 256 682.25856270 + 199.02142233 $i$ | 2929 | 260 236.29888627 + 199.29159120 $i$ | 2969 | 263 790.21288664 + 199.61489684 $i$ |
| 2890 | 256 771.20012871 + 199.27169369 $i$ | 2930 | 260 324.99447553 + 199.53605928 $i$ | 2970 | 263 878.85528443 + 199.67756090 $i$ |
| 2891 | 256 859.75902643 + 199.21711657 $i$ | 2931 | 260 413.73391213 + 199.30392969 $i$ | 2971 | 263 967.71917222 + 199.51487538 $i$ |
| 2892 | 256 948.83990644 + 199.01757152 $i$ | 2932 | 260 502.81231527 + 199.33668728 $i$ | 2972 | 264 056.74325587 + 199.55611478 $i$ |
| 2893 | 257 037.70790285 + 199.31029178 $i$ | 2933 | 260 591.56035297 + 199.49572858 $i$ | 2973 | 264 145.44185076 + 199.76243414 $i$ |
| 2894 | 257 126.34395422 + 199.20046840 $i$ | 2934 | 260 680.28758039 + 199.39572828 $i$ | 2974 | 264 234.21278179 + 199.48648588 $i$ |
| 2895 | 257 215.35156618 + 199.08563232 $i$ | 2935 | 260 769.31389893 + 199.27896067 $i$ | 2975 | 264 323.33308835 + 199.58903427 $i$ |
| 2896 | 257 304.26023705 + 199.25703861 $i$ | 2936 | 260 858.16294925 + 199.56280540 $i$ | 2976 | 264 411.95023906 + 199.77054457 $i$ |
| 2897 | 257 392.91393212 + 199.27828215 $i$ | 2937 | 260 946.78302148 + 199.38786725 $i$ | 2977 | 264 500.78887078 + 199.50290135 $i$ |
| 2898 | 257 481.84201086 + 199.06781606 $i$ | 2938 | 261 035.88667090 + 199.30070523 $i$ | 2978 | 264 589.84218926 + 199.61817279 $i$ |
| 2899 | 257 570.85944747 + 199.27858590 $i$ | 2939 | 261 124.69311666 + 199.58158177 $i$ | 2979 | 264 678.52388806 + 199.74889869 $i$ |
| 2900 | 257 659.42530519 + 199.30981891 $i$ | 2940 | 261 213.33722380 + 199.39025437 $i$ | 2980 | 264 767.32543466 + 199.57492769 $i$ |
| 2901 | 257 748.40046945 + 199.06052611 $i$ | 2941 | 261 302.42807048 + 199.35275895 $i$ | 2981 | 264 856.35852150 + 199.58380329 $i$ |
| 2902 | 257 837.39958591 + 199.31976553 $i$ | 2942 | 261 391.22465293 + 199.53659501 $i$ | 2982 | 264 945.12612560 + 199.79031183 $i$ |
| 2903 | 257 925.96708766 + 199.28764319 $i$ | 2943 | 261 479.92312101 + 199.46912005 $i$ | 2983 | 265 033.80693615 + 199.59711103 $i$ |
| 2904 | 258 014.96090234 + 199.13378706 $i$ | 2944 | 261 568.90302519 + 199.32494485 $i$ | 2984 | 265 122.95485310 + 199.57238402 $i$ |
| 2905 | 258 103.90132745 + 199.27031562 $i$ | 2945 | 261 657.84104408 + 199.56489058 $i$ | 2985 | 265 211.63290605 + 199.83627767 $i$ |
| 2906 | 258 192.57831415 + 199.35637400 $i$ | 2946 | 261 746.41446749 + 199.50018977 $i$ | 2986 | 265 300.38294880 + 199.58113652 $i$ |
| 2907 | 258 281.43554451 + 199.13002984 $i$ | 2947 | 261 835.47650803 + 199.31610675 $i$ | 2987 | 265 389.46568563 + 199.62782166 $i$ |
| 2908 | 258 370.49581026 + 199.27286016 $i$ | 2948 | 261 924.36997493 + 199.60369673 $i$ | 2988 | 265 478.20017985 + 199.80481381 $i$ |
| 2909 | 258 459.10165202 + 199.39988579 $i$ | 2949 | 262 012.96692209 + 199.49500544 $i$ | 2989 | 265 566.93610856 + 199.64528525 $i$ |
| 2910 | 258 547.97586873 + 199.12281389 $i$ | 2950 | 262 102.02833059 + 199.36512214 $i$ | 2990 | 265 655.97305040 + 199.61259797 $i$ |
| 2911 | 258 637.05767434 + 199.30989648 $i$ | 2951 | 262 190.88727093 + 199.57422360 $i$ | 2991 | 265 744.79545164 + 199.80957069 $i$ |
| 2912 | 258 725.62077181 + 199.38860326 $i$ | 2952 | 262 279.56250651 + 199.54071758 $i$ | 2992 | 265 833.43509566 + 199.70688121 $i$ |
| 2913 | 258 814.56306493 + 199.17433997 $i$ | 2953 | 262 368.50247392 + 199.38439703 $i$ | 2993 | 265 922.54300582 + 199.56532219 $i$ |
| 2914 | 258 903.53408934 + 199.29293934 $i$ | 2954 | 262 457.49293518 + 199.55525634 $i$ | 2994 | 266 011.32630043 + 199.88643763 $i$ |
| 2915 | 258 992.24680393 + 199.41613914 $i$ | 2955 | 262 546.07516768 + 199.61191685 $i$ | 2995 | 266 099.99149433 + 199.66893070 $i$ |
| 2916 | 259 081.03601214 + 199.21136320 $i$ | 2956 | 262 635.05600603 + 199.34345916 $i$ | 2996 | 266 189.07517291 + 199.63373969 $i$ |
| 2917 | 259 170.11922504 + 199.26835674 $i$ | 2957 | 262 724.03638790 + 199.61489359 $i$ | 2997 | 266 277.87213759 + 199.85096332 $i$ |
| 2918 | 259 258.78880183 + 199.47055675 $i$ | 2958 | 262 812.61716013 + 199.59102866 $i$ | 2998 | 266 366.55970501 + 199.72018368 $i$ |
| 2919 | 259 347.55878619 + 199.19950040 $i$ | 2959 | 262 901.60775862 + 199.39755370 $i$ | 2999 | 266 455.57349290 + 199.64522552 $i$ |
| 2920 | 259 436.68921331 + 199.30284576 $i$ | 2960 | 262 990.55979644 + 199.60404139 $i$ | 3000 | 266 544.46165009 + 199.83312439 $i$ |

| $k$ | $\sigma_k$ |
|---|---|
| 3001 | 266 633.08286521 + 199.79675174 $i$ |
| 3002 | 266 722.11364725 + 199.59489006 $i$ |
| 3003 | 266 811.02578761 + 199.89948878 $i$ |
| 3004 | 266 899.60238363 + 199.77228827 $i$ |
| 3005 | 266 988.68301849 + 199.64947782 $i$ |
| 3006 | 267 077.53658954 + 199.88013330 $i$ |
| 3007 | 267 166.19566021 + 199.80441393 $i$ |
| 3008 | 267 255.17416349 + 199.68011668 $i$ |
| 3009 | 267 344.11483912 + 199.84923520 $i$ |
| 3010 | 267 432.73802429 + 199.87927001 $i$ |
| 3011 | 267 521.69158768 + 199.65222979 $i$ |
| 3012 | 267 610.70216220 + 199.88476283 $i$ |
| 3013 | 267 699.24009006 + 199.89007931 $i$ |
| 3014 | 267 788.28058317 + 199.66583964 $i$ |
| 3015 | 267 877.19114684 + 199.90492095 $i$ |
| 3016 | 267 965.84836886 + 199.86642791 $i$ |
| 3017 | 268 054.76894649 + 199.72134075 $i$ |
| 3018 | 268 143.76016981 + 199.86787122 $i$ |
| 3019 | 268 232.40703812 + 199.94558419 $i$ |
| 3020 | 268 321.27252160 + 199.72242357 $i$ |
| 3021 | 268 410.35342949 + 199.86789749 $i$ |
| 3022 | 268 498.90899769 + 199.98980582 $i$ |
| 3023 | 268 587.85523061 + 199.70494222 $i$ |
| 3024 | 268 676.85386212 + 199.91978186 $i$ |
| 3025 | 268 765.49793350 + 199.95692452 $i$ |
| 3026 | 268 854.37042648 + 199.78216039 $i$ |
| 3027 | 268 943.39492260 + 199.87627346 $i$ |
| 3028 | 269 032.07601696 + 200.00746547 $i$ |
| 3029 | 269 120.87069328 + 199.80284636 $i$ |
| 3030 | 269 209.97855398 + 199.84475574 $i$ |
| 3031 | 269 298.59716159 + 200.06737767 $i$ |
| 3032 | 269 387.42808333 + 199.77312374 $i$ |
| 3033 | 269 476.50819985 + 199.91577840 $i$ |
| 3034 | 269 565.14960886 + 200.03191092 $i$ |
| 3035 | 269 653.98764665 + 199.84374171 $i$ |
| 3036 | 269 743.00540684 + 199.88671356 $i$ |
| 3037 | 269 831.75859147 + 200.06577915 $i$ |
| 3038 | 269 920.47898960 + 199.87497962 $i$ |
| 3039 | 270 009.57686923 + 199.86313127 $i$ |
| 3040 | 270 098.30111342 + 200.11910092 $i$ |

| $k$ | $\sigma_k$ | $k$ | $\sigma_k$ | $k$ | $\sigma_k$ |
|---|---|---|---|---|---|
| 3041 | 270 187.00917533 + 199.86193530 $i$ | 3081 | 273 741.05610542 + 200.04200499 $i$ | 3121 | 277 294.99714890 + 200.24285937 $i$ |
| 3042 | 270 276.14288455 + 199.90503706 $i$ | 3082 | 273 830.01772749 + 200.18001385 $i$ | 3122 | 277 383.97427454 + 200.42828992 $i$ |
| 3043 | 270 364.81249916 + 200.10150685 $i$ | 3083 | 273 918.73147331 + 200.27396621 $i$ | 3123 | 277 472.55459064 + 200.45083554 $i$ |
| 3044 | 270 453.60509991 + 199.91103979 $i$ | 3084 | 274 007.52244669 + 200.07535869 $i$ | 3124 | 277 561.55435817 + 200.20755176 $i$ |
| 3045 | 270 542.61644666 + 199.90724940 $i$ | 3085 | 274 096.61568275 + 200.15351175 $i$ | 3125 | 277 650.51677998 + 200.48756317 $i$ |
| 3046 | 270 631.43306523 + 200.09918917 $i$ | 3086 | 274 185.25461571 + 200.32924251 $i$ | 3126 | 277 739.09929035 + 200.42329215 $i$ |
| 3047 | 270 720.09458905 + 199.96580342 $i$ | 3087 | 274 274.06345926 + 200.06286760 $i$ | 3127 | 277 828.10358075 + 200.26765765 $i$ |
| 3048 | 270 809.17849891 + 199.88227508 $i$ | 3088 | 274 363.17333222 + 200.18613332 $i$ | 3128 | 277 917.03755840 + 200.46129598 $i$ |
| 3049 | 270 897.99197827 + 200.13955278 $i$ | 3089 | 274 451.77876989 + 200.32673529 $i$ | 3129 | 278 005.68361991 + 200.45956502 $i$ |
| 3050 | 270 986.61074293 + 199.97285324 $i$ | 3090 | 274 540.64295705 + 200.09493876 $i$ | 3130 | 278 094.60745392 + 200.30398889 $i$ |
| 3051 | 271 075.75750374 + 199.89711835 $i$ | 3091 | 274 629.66155973 + 200.19622923 $i$ | 3131 | 278 183.60879663 + 200.41837221 $i$ |
| 3052 | 271 164.49016667 + 200.15458103 $i$ | 3092 | 274 718.38400317 + 200.32127239 $i$ | 3132 | 278 272.23074655 + 200.54358615 $i$ |
| 3053 | 271 253.21951415 + 199.98280638 $i$ | 3093 | 274 807.14662512 + 200.16390997 $i$ | 3133 | 278 361.12799296 + 200.25422801 $i$ |
| 3054 | 271 342.22844026 + 199.93467793 $i$ | 3094 | 274 896.21679358 + 200.14014159 $i$ | 3134 | 278 450.17892862 + 200.48274577 $i$ |
| 3055 | 271 431.09445303 + 200.12206120 $i$ | 3095 | 274 984.94379096 + 200.39587260 $i$ | 3135 | 278 538.75234674 + 200.51800754 $i$ |
| 3056 | 271 519.73660432 + 200.05613279 $i$ | 3096 | 275 073.66128606 + 200.14181381 $i$ | 3136 | 278 627.69825837 + 200.30684993 $i$ |
| 3057 | 271 608.76091594 + 199.90373009 $i$ | 3097 | 275 162.78924207 + 200.18028415 $i$ | 3137 | 278 716.68788347 + 200.46966692 $i$ |
| 3058 | 271 697.67187264 + 200.15979992 $i$ | 3098 | 275 251.46658499 + 200.39142826 $i$ | 3138 | 278 805.33655129 + 200.53274020 $i$ |
| 3059 | 271 786.24672793 + 200.07215050 $i$ | 3099 | 275 340.23168634 + 200.16454139 $i$ | 3139 | 278 894.21614464 + 200.35772273 $i$ |
| 3060 | 271 875.33462604 + 199.90735878 $i$ | 3100 | 275 429.29684968 + 200.20985442 $i$ | 3140 | 278 983.22800627 + 200.42400921 $i$ |
| 3061 | 271 964.18290036 + 200.19309626 $i$ | 3101 | 275 518.04054305 + 200.36599715 $i$ | 3141 | 279 071.91738429 + 200.60690394 $i$ |
| 3062 | 272 052.83303069 + 200.06054698 $i$ | 3102 | 275 606.77867545 + 200.24369990 $i$ | 3142 | 279 160.70520188 + 200.33094684 $i$ |
| 3063 | 272 141.84181920 + 199.97124749 $i$ | 3103 | 275 695.80420626 + 200.15619712 $i$ | 3143 | 279 249.82405776 + 200.46579238 $i$ |
| 3064 | 272 230.74388416 + 200.13964802 $i$ | 3104 | 275 784.64308581 + 200.43006352 $i$ | 3144 | 279 338.42314102 + 200.60472269 $i$ |
| 3065 | 272 319.39540206 + 200.14256193 $i$ | 3105 | 275 873.26517156 + 200.23202968 $i$ | 3145 | 279 427.28664229 + 200.35575279 $i$ |
| 3066 | 272 408.34114277 + 199.94250416 $i$ | 3106 | 275 962.39133002 + 200.18525413 $i$ | 3146 | 279 516.32768050 + 200.48392390 $i$ |
| 3067 | 272 497.34023057 + 200.16555334 $i$ | 3107 | 276 051.15830476 + 200.43447565 $i$ | 3147 | 279 605.00571658 + 200.58951515 $i$ |
| 3068 | 272 585.89877876 + 200.16600619 $i$ | 3108 | 276 139.83293829 + 200.24924066 $i$ | 3148 | 279 693.81284192 + 200.42139612 $i$ |
| 3069 | 272 674.90929797 + 199.94369229 $i$ | 3109 | 276 228.92052741 + 200.21913041 $i$ | 3149 | 279 782.85303611 + 200.44515701 $i$ |
| 3070 | 272 763.86740595 + 200.19686578 $i$ | 3110 | 276 317.69844386 + 200.40192734 $i$ | 3150 | 279 871.59709193 + 200.63990897 $i$ |
| 3071 | 272 852.45378589 + 200.15710038 $i$ | 3111 | 276 406.41069044 + 200.32199863 $i$ | 3151 | 279 960.30290901 + 200.43308083 $i$ |
| 3072 | 272 941.46187734 + 200.00439765 $i$ | 3112 | 276 495.40007083 + 200.19216021 $i$ | 3152 | 280 049.44411747 + 200.44093634 $i$ |
| 3073 | 273 030.37712379 + 200.15151992 $i$ | 3113 | 276 584.31917843 + 200.43297994 $i$ | 3153 | 280 138.10718628 + 200.68185393 $i$ |
| 3074 | 273 119.06176172 + 200.22128831 $i$ | 3114 | 276 672.89480447 + 200.34440719 $i$ | 3154 | 280 226.88188040 + 200.41561744 $i$ |
| 3075 | 273 207.92853724 + 199.99808821 $i$ | 3115 | 276 761.98282789 + 200.18758803 $i$ | 3155 | 280 315.94890009 + 200.49300992 $i$ |
| 3076 | 273 296.98611930 + 200.16077940 $i$ | 3116 | 276 850.83785390 + 200.46647543 $i$ | 3156 | 280 404.67712346 + 200.64762781 $i$ |
| 3077 | 273 385.56986300 + 200.25528616 $i$ | 3117 | 276 939.45860156 + 200.33882089 $i$ | 3157 | 280 493.42814859 + 200.49182385 $i$ |
| 3078 | 273 474.48456729 + 199.99531414 $i$ | 3118 | 277 028.51780830 + 200.23425608 $i$ | 3158 | 280 582.46676975 + 200.46142453 $i$ |
| 3079 | 273 563.53107176 + 200.19147916 $i$ | 3119 | 277 117.36706678 + 200.43983185 $i$ | 3159 | 280 671.26610357 + 200.66918007 $i$ |
| 3080 | 273 652.10402716 + 200.24852996 $i$ | 3120 | 277 206.04823026 + 200.38945762 $i$ | 3160 | 280 759.93136932 + 200.53370310 $i$ |
| 3161 | 280 849.03144637 + 200.42781683 $i$ | | | | |
| 3162 | 280 937.80262237 + 200.73824981 $i$ | | | | |
| 3163 | 281 026.48241591 + 200.49598955 $i$ | | | | |
| 3164 | 281 115.56910730 + 200.50156497 $i$ | | | | |
| 3165 | 281 204.34831748 + 200.68642384 $i$ | | | | |
| 3166 | 281 293.04490289 + 200.56466828 $i$ | | | | |
| 3167 | 281 382.07134157 + 200.49334821 $i$ | | | | |
| 3168 | 281 470.93693519 + 200.68657257 $i$ | | | | |
| 3169 | 281 559.56595022 + 200.62597890 $i$ | | | | |
| 3170 | 281 648.61406446 + 200.45148457 $i$ | | | | |
| 3171 | 281 737.49906573 + 200.75031309 $i$ | | | | |
| 3172 | 281 826.08916607 + 200.59945506 $i$ | | | | |
| 3173 | 281 915.18203515 + 200.50712938 $i$ | | | | |
| 3174 | 282 004.00858388 + 200.72282805 $i$ | | | | |
| 3175 | 282 092.68371434 + 200.64246811 $i$ | | | | |
| 3176 | 282 181.66875017 + 200.52570392 $i$ | | | | |
| 3177 | 282 270.59691357 + 200.70166343 $i$ | | | | |
| 3178 | 282 359.21745912 + 200.70904939 $i$ | | | | |
| 3179 | 282 448.19409162 + 200.49359662 $i$ | | | | |
| 3180 | 282 537.17341820 + 200.74361337 $i$ | | | | |
| 3181 | 282 625.72962840 + 200.71195700 $i$ | | | | |
| 3182 | 282 714.77296170 + 200.51166650 $i$ | | | | |
| 3183 | 282 803.66847505 + 200.76028159 $i$ | | | | |
| 3184 | 282 892.33520533 + 200.71067839 $i$ | | | | |
| 3185 | 282 981.25906805 + 200.56969045 $i$ | | | | |
| 3186 | 283 070.24932229 + 200.71413939 $i$ | | | | |
| 3187 | 283 158.87801038 + 200.78015754 $i$ | | | | |
| 3188 | 283 247.77871060 + 200.55661445 $i$ | | | | |
| 3189 | 283 336.82849467 + 200.72571139 $i$ | | | | |
| 3190 | 283 425.39172447 + 200.81479259 $i$ | | | | |
| 3191 | 283 514.35443596 + 200.54317083 $i$ | | | | |
| 3192 | 283 603.33377444 + 200.77149974 $i$ | | | | |
| 3193 | 283 691.97695217 + 200.78060371 $i$ | | | | |
| 3194 | 283 780.86596533 + 200.62638990 $i$ | | | | |
| 3195 | 283 869.88201873 + 200.71664123 $i$ | | | | |
| 3196 | 283 958.54850416 + 200.84705217 $i$ | | | | |
| 3197 | 284 047.37244104 + 200.62529277 $i$ | | | | |
| 3198 | 284 136.45724788 + 200.71286273 $i$ | | | | |
| 3199 | 284 225.07503872 + 200.89747631 $i$ | | | | |
| 3200 | 284 313.92991184 + 200.60122562 $i$ | | | | |

| k | $\sigma_k$ | k | $\sigma_k$ | k | $\sigma_k$ |
|---|---|---|---|---|---|
| 3201 | 284 402.98808137 + 200.76953749 $i$ | 3241 | 287 956.85948255 + 200.99266074 $i$ | 3281 | 291 510.79627731 + 201.25709073 $i$ |
| 3202 | 284 491.63089688 + 200.85724712 $i$ | 3242 | 288 045.54663739 + 201.03192245 $i$ | 3282 | 291 599.38023918 + 201.14214020 $i$ |
| 3203 | 284 580.48012926 + 200.67638197 $i$ | 3243 | 288 134.41901187 + 200.82161955 $i$ | 3283 | 291 688.48343798 + 201.01187643 $i$ |
| 3204 | 284 669.49353419 + 200.73294497 $i$ | 3244 | 288 223.47449336 + 201.00251130 $i$ | 3284 | 291 777.30388370 + 201.28754666 $i$ |
| 3205 | 284 758.23656690 + 200.89818722 $i$ | 3245 | 288 312.04225890 + 201.06372668 $i$ | 3285 | 291 865.95499324 + 201.13937946 $i$ |
| 3206 | 284 846.96709989 + 200.69873581 $i$ | 3246 | 288 400.99036597 + 200.81874196 $i$ | 3286 | 291 955.00689994 + 201.05764151 $i$ |
| 3207 | 284 936.07447044 + 200.71600113 $i$ | 3247 | 288 490.00059562 + 201.03117289 $i$ | 3287 | 292 043.84839361 + 201.25802869 $i$ |
| 3208 | 285 024.76732968 + 200.94589818 $i$ | 3248 | 288 578.59498928 + 201.06167551 $i$ | 3288 | 292 132.53252346 + 201.19223535 $i$ |
| 3209 | 285 113.50935624 + 200.68838391 $i$ | 3249 | 288 667.54796330 + 200.86131740 $i$ | 3289 | 292 221.49309015 + 201.05725015 $i$ |
| 3210 | 285 202.62932269 + 200.75573794 $i$ | 3250 | 288 756.50140312 + 201.01640979 $i$ | 3290 | 292 310.45316902 + 201.25319279 $i$ |
| 3211 | 285 291.29185977 + 200.92773737 $i$ | 3251 | 288 845.21039794 + 201.08717011 $i$ | 3291 | 292 399.03357537 + 201.24705527 $i$ |
| 3212 | 285 380.09568976 + 200.73529735 $i$ | 3252 | 288 934.01597172 + 200.89541552 $i$ | 3292 | 292 488.05725254 + 201.02765762 $i$ |
| 3213 | 285 469.10676819 + 200.75381391 $i$ | 3253 | 289 023.10594536 + 200.98622391 $i$ | 3293 | 292 576.99067484 + 201.30781810 $i$ |
| 3214 | 285 557.91309439 + 200.92800297 $i$ | 3254 | 289 111.71952518 + 201.14782032 $i$ | 3294 | 292 665.58146415 + 201.22086259 $i$ |
| 3215 | 285 646.57851920 + 200.78871368 $i$ | 3255 | 289 200.57384550 + 200.87850554 $i$ | 3295 | 292 754.60673174 + 201.09030117 $i$ |
| 3216 | 285 735.68008437 + 200.72194841 $i$ | 3256 | 289 289.64937947 + 201.02153779 $i$ | 3296 | 292 843.51307808 + 201.27016006 $i$ |
| 3217 | 285 824.45461824 + 200.97799405 $i$ | 3257 | 289 378.25696053 + 201.14104204 $i$ | 3297 | 292 932.16778411 + 201.26881223 $i$ |
| 3218 | 285 913.11489042 + 200.78643271 $i$ | 3258 | 289 467.14318651 + 200.91449736 $i$ | 3298 | 293 021.10306196 + 201.10457919 $i$ |
| 3219 | 286 002.23925949 + 200.74035779 $i$ | 3259 | 289 556.14601965 + 201.02438045 $i$ | 3299 | 293 110.08701384 + 201.24848660 $i$ |
| 3220 | 286 090.96797404 + 200.99122877 $i$ | 3260 | 289 644.86313685 + 201.13855180 $i$ | 3300 | 293 198.71226360 + 201.33960926 $i$ |
| 3221 | 286 179.71238239 + 200.79739111 $i$ | 3261 | 289 733.63648953 + 200.97279427 $i$ | 3301 | 293 287.62626403 + 201.06247172 $i$ |
| 3222 | 286 268.71460961 + 200.77833743 $i$ | 3262 | 289 822.70741621 + 200.97705314 $i$ | 3302 | 293 376.65942254 + 201.30882231 $i$ |
| 3223 | 286 357.57583517 + 200.95598106 $i$ | 3263 | 289 911.41897548 + 201.21268007 $i$ | 3303 | 293 465.22961340 + 201.31028025 $i$ |
| 3224 | 286 446.21934713 + 200.87816269 $i$ | 3264 | 290 000.15600828 + 200.94537358 $i$ | 3304 | 293 554.19990130 + 201.12138831 $i$ |
| 3225 | 286 535.26408193 + 200.73738978 $i$ | 3265 | 290 089.27867443 + 201.02142314 $i$ | 3305 | 293 643.16204254 + 201.28594840 $i$ |
| 3226 | 286 624.14059537 + 201.00057842 $i$ | 3266 | 290 177.93956297 + 201.19865524 $i$ | 3306 | 293 731.82463404 + 201.33884338 $i$ |
| 3227 | 286 712.73796610 + 200.88288829 $i$ | 3267 | 290 266.72893290 + 200.97907200 $i$ | 3307 | 293 820.70494275 + 201.15497379 $i$ |
| 3228 | 286 801.83164492 + 200.75002803 $i$ | 3268 | 290 355.78066117 + 201.03739781 $i$ | 3308 | 293 909.71522268 + 201.25457580 $i$ |
| 3229 | 286 890.65675373 + 201.02203569 $i$ | 3269 | 290 444.52246748 + 201.18613361 $i$ | 3309 | 293 998.39563194 + 201.40260701 $i$ |
| 3230 | 286 979.31999955 + 200.87817618 $i$ | 3270 | 290 533.26430977 + 201.04648122 $i$ | 3310 | 294 087.20060773 + 201.13240241 $i$ |
| 3231 | 287 068.33700997 + 200.80836882 $i$ | 3271 | 290 622.29802869 + 200.99112115 $i$ | 3311 | 294 176.31129392 + 201.29206373 $i$ |
| 3232 | 287 157.22384287 + 200.96912994 $i$ | 3272 | 290 711.11924215 + 201.24683801 $i$ | 3312 | 294 264.89493673 + 201.39777404 $i$ |
| 3233 | 287 245.87321678 + 200.96297673 $i$ | 3273 | 290 799.75157461 + 201.03603486 $i$ | 3313 | 294 353.78885335 + 201.16319925 $i$ |
| 3234 | 287 334.84317142 + 200.77261487 $i$ | 3274 | 290 888.89285590 + 201.01808379 $i$ | 3314 | 294 442.80571714 + 201.30048518 $i$ |
| 3235 | 287 423.81754204 + 201.00418341 $i$ | 3275 | 290 977.61979362 + 201.24747845 $i$ | 3315 | 294 531.48686783 + 201.39438001 $i$ |
| 3236 | 287 512.37684970 + 200.97720536 $i$ | 3276 | 291 066.33595285 + 201.06240176 $i$ | 3316 | 294 620.30759673 + 201.22180299 $i$ |
| 3237 | 287 601.41559242 + 200.77816368 $i$ | 3277 | 291 155.40611067 + 201.03918359 $i$ | 3317 | 294 709.34209748 + 201.25907451 $i$ |
| 3238 | 287 690.33509407 + 201.03048314 $i$ | 3278 | 291 244.17792993 + 201.22605581 $i$ | 3318 | 294 798.06824549 + 201.45001117 $i$ |
| 3239 | 287 778.94629592 + 200.97407320 $i$ | 3279 | 291 332.89951962 + 201.12309616 $i$ | 3319 | 294 886.80584175 + 201.22323450 $i$ |
| 3240 | 287 867.95214856 + 200.82931217 $i$ | 3280 | 291 421.89320925 + 201.01532327 $i$ | 3320 | 294 975.92954365 + 201.26094500 $i$ |
| 3321 | 295 064.57822922 + 201.48159666 $i$ | | | | |
| 3322 | 295 153.38133764 + 201.21323785 $i$ | | | | |
| 3323 | 295 242.43661028 + 201.31123961 $i$ | | | | |
| 3324 | 295 331.15400159 + 201.44231458 $i$ | | | | |
| 3325 | 295 419.91769058 + 201.29206203 $i$ | | | | |
| 3326 | 295 508.95734699 + 201.27047602 $i$ | | | | |
| 3327 | 295 597.74095975 + 201.48240335 $i$ | | | | |
| 3328 | 295 686.42311460 + 201.31566692 $i$ | | | | |
| 3329 | 295 775.52403118 + 201.25345050 $i$ | | | | |
| 3330 | 295 864.28017528 + 201.53440397 $i$ | | | | |
| 3331 | 295 952.96812833 + 201.28480790 $i$ | | | | |
| 3332 | 296 042.06417062 + 201.32479050 $i$ | | | | |
| 3333 | 296 130.82308130 + 201.48065591 $i$ | | | | |
| 3334 | 296 219.53675202 + 201.36416035 $i$ | | | | |
| 3335 | 296 308.56355521 + 201.29367674 $i$ | | | | |
| 3336 | 296 397.41200450 + 201.50169556 $i$ | | | | |
| 3337 | 296 486.05529484 + 201.40790147 $i$ | | | | |
| 3338 | 296 575.11093564 + 201.26191663 $i$ | | | | |
| 3339 | 296 663.96868423 + 201.55733697 $i$ | | | | |
| 3340 | 296 752.58119662 + 201.38524561 $i$ | | | | |
| 3341 | 296 841.67711250 + 201.31693442 $i$ | | | | |
| 3342 | 296 930.47923230 + 201.52621738 $i$ | | | | |
| 3343 | 297 019.17838972 + 201.43591266 $i$ | | | | |
| 3344 | 297 108.16034292 + 201.32241795 $i$ | | | | |
| 3345 | 297 197.07489458 + 201.51268071 $i$ | | | | |
| 3346 | 297 285.69917117 + 201.49562091 $i$ | | | | |
| 3347 | 297 374.69649697 + 201.29353134 $i$ | | | | |
| 3348 | 297 463.64458164 + 201.55587609 $i$ | | | | |
| 3349 | 297 552.21702112 + 201.49099713 $i$ | | | | |
| 3350 | 297 641.26794575 + 201.32040738 $i$ | | | | |
| 3351 | 297 730.14846448 + 201.56359664 $i$ | | | | |
| 3352 | 297 818.81731340 + 201.49557853 $i$ | | | | |
| 3353 | 297 907.75646982 + 201.37377289 $i$ | | | | |
| 3354 | 297 996.73320885 + 201.51200163 $i$ | | | | |
| 3355 | 298 085.35079855 + 201.57310884 $i$ | | | | |
| 3356 | 298 174.28461974 + 201.34508024 $i$ | | | | |
| 3357 | 298 263.30058332 + 201.54043466 $i$ | | | | |
| 3358 | 298 351.87470630 + 201.59396589 $i$ | | | | |
| 3359 | 298 440.85046613 + 201.34367808 $i$ | | | | |
| 3360 | 298 529.81481767 + 201.57971108 $i$ | | | | |

| $k$ | $\sigma_k$ | $k$ | $\sigma_k$ | $k$ | $\sigma_k$ |
|---|---|---|---|---|---|
| 3361 | 298 618.45945216 + 201.56380167 $i$ | 3401 | 302 172.35811982 + 201.74021692 $i$ | 3441 | 305 726.24514042 + 201.80446753 $i$ |
| 3362 | 298 707.36291186 + 201.42169885 $i$ | 3402 | 302 261.34205470 + 201.55511922 $i$ | 3442 | 305 815.39151298 + 201.80308567 $i$ |
| 3363 | 298 796.36257836 + 201.51710793 $i$ | 3403 | 302 350.28886367 + 201.80486098 $i$ | 3443 | 305 904.08214177 + 202.02041368 $i$ |
| 3364 | 298 885.02757779 + 201.64383972 $i$ | 3404 | 302 438.86394876 + 201.74692405 $i$ | 3444 | 305 992.83941316 + 201.82943939 $i$ |
| 3365 | 298 973.87139640 + 201.40184748 $i$ | 3405 | 302 527.91779710 + 201.56629199 $i$ | 3445 | 306 081.88599046 + 201.82117958 $i$ |
| 3366 | 299 062.93633026 + 201.53225092 $i$ | 3406 | 302 616.80187398 + 201.82413062 $i$ | 3446 | 306 170.66171695 + 202.00734646 $i$ |
| 3367 | 299 151.55610766 + 201.67695130 $i$ | 3407 | 302 705.44205326 + 201.74668171 $i$ | 3447 | 306 259.38369039 + 201.88265334 $i$ |
| 3368 | 299 240.42671102 + 201.38948848 $i$ | 3408 | 302 794.44001290 + 201.61272249 $i$ | 3448 | 306 348.38959504 + 201.80133864 $i$ |
| 3369 | 299 329.46826332 + 201.57979867 $i$ | 3409 | 302 883.34242287 + 201.78801941 $i$ | 3449 | 306 437.27098339 + 202.03512206 $i$ |
| 3370 | 299 418.11159241 + 201.64043857 $i$ | 3410 | 302 972.02566016 + 201.80313197 $i$ | 3450 | 306 525.86506236 + 201.90206167 $i$ |
| 3371 | 299 506.97506096 + 201.46841477 $i$ | 3411 | 303 060.91701526 + 201.60645061 $i$ | 3451 | 306 614.98368616 + 201.79897435 $i$ |
| 3372 | 299 595.98210026 + 201.53251961 $i$ | 3412 | 303 149.95612368 + 201.79528228 $i$ | 3452 | 306 703.77385814 + 202.06549659 $i$ |
| 3373 | 299 684.71110882 + 201.68842263 $i$ | 3413 | 303 238.51308328 + 201.83766316 $i$ | 3453 | 306 792.45165876 + 201.89693851 $i$ |
| 3374 | 299 773.46038113 + 201.48105848 $i$ | 3414 | 303 327.50099260 + 201.60371292 $i$ | 3454 | 306 881.49375676 + 201.84020005 $i$ |
| 3375 | 299 862.56821333 + 201.52262902 $i$ | 3415 | 303 416.47121885 + 201.82486000 $i$ | 3455 | 306 970.32813182 + 202.03377112 $i$ |
| 3376 | 299 951.23377094 + 201.73196477 $i$ | 3416 | 303 505.08371008 + 201.83119161 $i$ | 3456 | 307 059.01522748 + 201.95785569 $i$ |
| 3377 | 300 040.01407700 + 201.47092637 $i$ | 3417 | 303 594.03928160 + 201.64485249 $i$ | 3457 | 307 147.99237711 + 201.83218561 $i$ |
| 3378 | 300 129.11045166 + 201.56003172 $i$ | 3418 | 303 682.98835610 + 201.80653350 $i$ | 3458 | 307 236.92911334 + 202.03581328 $i$ |
| 3379 | 300 217.77030156 + 201.71322773 $i$ | 3419 | 303 771.68462718 + 201.86014749 $i$ | 3459 | 307 325.51582413 + 202.00478526 $i$ |
| 3380 | 300 306.58834942 + 201.51943172 $i$ | 3420 | 303 860.51488616 + 201.67531047 $i$ | 3460 | 307 414.56093967 + 201.80528400 $i$ |
| 3381 | 300 395.59641673 + 201.55363190 $i$ | 3421 | 303 949.59219281 + 201.77596886 $i$ | 3461 | 307 503.45945048 + 202.08657813 $i$ |
| 3382 | 300 484.38851180 + 201.71835828 $i$ | 3422 | 304 038.18941713 + 201.92411027 $i$ | 3462 | 307 592.07012244 + 201.98112238 $i$ |
| 3383 | 300 573.07018430 + 201.56933836 $i$ | 3423 | 304 127.08095435 + 201.64954732 $i$ | 3463 | 307 681.10532534 + 201.86535728 $i$ |
| 3384 | 300 662.17590966 + 201.51681303 $i$ | 3424 | 304 216.12353827 + 201.82060207 $i$ | 3464 | 307 769.98283828 + 202.04584364 $i$ |
| 3385 | 300 750.91838096 + 201.77579843 $i$ | 3425 | 304 304.74322658 + 201.90810909 $i$ | 3465 | 307 858.65845267 + 202.03928946 $i$ |
| 3386 | 300 839.61823096 + 201.55745536 $i$ | 3426 | 304 393.63477671 + 201.68883160 $i$ | 3466 | 307 947.59768738 + 201.86522258 $i$ |
| 3387 | 300 928.72215839 + 201.54745364 $i$ | 3427 | 304 482.62757689 + 201.81724671 $i$ | 3467 | 308 036.56816024 + 202.03784222 $i$ |
| 3388 | 301 017.45033698 + 201.77674179 $i$ | 3428 | 304 571.34582100 + 201.91551983 $i$ | 3468 | 308 125.19163728 + 202.09122213 $i$ |
| 3389 | 301 106.19881427 + 201.57403996 $i$ | 3429 | 304 660.12741816 + 201.73977017 $i$ | 3469 | 308 214.12576431 + 201.83537254 $i$ |
| 3390 | 301 195.20831010 + 201.58008195 $i$ | 3430 | 304 749.19700313 + 201.77293392 $i$ | 3470 | 308 303.13773593 + 202.09007567 $i$ |
| 3391 | 301 284.05408658 + 201.74129842 $i$ | 3431 | 304 837.89260823 + 201.98694631 $i$ | 3471 | 308 391.70872926 + 202.06595155 $i$ |
| 3392 | 301 372.70211846 + 201.65903914 $i$ | 3432 | 304 926.65389808 + 201.71125567 $i$ | 3472 | 308 480.70310426 + 201.89386821 $i$ |
| 3393 | 301 461.76659108 + 201.52894829 $i$ | 3433 | 305 015.76461203 + 201.81780651 $i$ | 3473 | 308 569.63509889 + 202.06004728 $i$ |
| 3394 | 301 550.60894936 + 201.79684229 $i$ | 3434 | 305 104.41182743 + 201.96912374 $i$ | 3474 | 308 658.31123897 + 202.10353786 $i$ |
| 3395 | 301 639.22942932 + 201.65037123 $i$ | 3435 | 305 193.23213053 + 201.75415169 $i$ | 3475 | 308 747.19511073 + 201.91664379 $i$ |
| 3396 | 301 728.32571690 + 201.55281797 $i$ | 3436 | 305 282.26402120 + 201.81842851 $i$ | 3476 | 308 836.20624328 + 202.03885782 $i$ |
| 3397 | 301 817.13160041 + 201.80814209 $i$ | 3437 | 305 371.00091905 + 201.96507774 $i$ | 3477 | 308 924.86616693 + 202.15610210 $i$ |
| 3398 | 301 905.81003556 + 201.65482427 $i$ | 3438 | 305 459.75389436 + 201.81244558 $i$ | 3478 | 309 013.70183253 + 201.90005435 $i$ |
| 3399 | 301 994.83040007 + 201.59908445 $i$ | 3439 | 305 548.79633646 + 201.77774248 $i$ | 3479 | 309 102.79563156 + 202.07131926 $i$ |
| 3400 | 302 083.70077794 + 201.76035880 $i$ | 3440 | 305 637.58600674 + 202.01932451 $i$ | 3480 | 309 191.36752893 + 202.15421979 $i$ |
| $k$ | $\sigma_k$ | | | | |
| 3481 | 309 280.29462754 + 201.92849900 $i$ | | | | |
| 3482 | 309 369.28071023 + 202.07367935 $i$ | | | | |
| 3483 | 309 457.96927900 + 202.15880791 $i$ | | | | |
| 3484 | 309 546.80109652 + 201.98027478 $i$ | | | | |
| 3485 | 309 635.83080840 + 202.03482784 $i$ | | | | |
| 3486 | 309 724.54092809 + 202.21437580 $i$ | | | | |
| 3487 | 309 813.30308311 + 201.97355959 $i$ | | | | |
| 3488 | 309 902.41308463 + 202.04759960 $i$ | | | | |
| 3489 | 309 991.05382449 + 202.23866602 $i$ | | | | |
| 3490 | 310 079.87899672 + 201.97248724 $i$ | | | | |
| 3491 | 310 168.92227564 + 202.08986759 $i$ | | | | |
| 3492 | 310 257.63223059 + 202.20040829 $i$ | | | | |
| 3493 | 310 346.40980079 + 202.05165474 $i$ | | | | |
| 3494 | 310 435.44731072 + 202.04100868 $i$ | | | | |
| 3495 | 310 524.21738167 + 202.25277416 $i$ | | | | |
| 3496 | 310 612.91451762 + 202.05828749 $i$ | | | | |
| 3497 | 310 702.01633201 + 202.03779512 $i$ | | | | |
| 3498 | 310 790.75100277 + 202.29078030 $i$ | | | | |
| 3499 | 310 879.46389993 + 202.04310027 $i$ | | | | |
| 3500 | 310 968.55950406 + 202.09413750 $i$ | | | | |
| 3501 | 311 057.28863775 + 202.23953885 $i$ | | | | |
| 3502 | 311 146.03453631 + 202.12676191 $i$ | | | | |
| 3503 | 311 235.05427188 + 202.05241213 $i$ | | | | |
| 3504 | 311 323.88634788 + 202.27643712 $i$ | | | | |
| 3505 | 311 412.54382775 + 202.15129020 $i$ | | | | |
| 3506 | 311 501.61013867 + 202.03603023 $i$ | | | | |
| 3507 | 311 590.43738357 + 202.31902138 $i$ | | | | |
| 3508 | 311 679.07276034 + 202.13502283 $i$ | | | | |
| 3509 | 311 768.17293900 + 202.08838574 $i$ | | | | |
| 3510 | 311 856.95297526 + 202.28627636 $i$ | | | | |
| 3511 | 311 945.66854732 + 202.18716994 $i$ | | | | |
| 3512 | 312 034.65071913 + 202.08574064 $i$ | | | | |
| 3513 | 312 123.55705868 + 202.28196754 $i$ | | | | |
| 3514 | 312 212.17794831 + 202.23906975 $i$ | | | | |
| 3515 | 312 301.19741103 + 202.06046477 $i$ | | | | |
| 3516 | 312 390.11806193 + 202.32646033 $i$ | | | | |
| 3517 | 312 478.70551324 + 202.23075384 $i$ | | | | |
| 3518 | 312 567.76401081 + 202.08966130 $i$ | | | | |
| 3519 | 312 656.62510513 + 202.32322917 $i$ | | | | |
| 3520 | 312 745.30040126 + 202.24620209 $i$ | | | | |